\documentclass[11pt,a4paper]{article}
\usepackage[T1]{fontenc}
\usepackage[utf8]{inputenc}
\usepackage[english]{babel}
\usepackage{lmodern,microtype}
\usepackage[margin=2.5cm]{geometry}
\usepackage{amsmath,amssymb,amsthm,amscd,dsfont}
\usepackage{enumerate,enumitem}
\usepackage[table,svgnames]{xcolor}
\usepackage{array,multicol,float,graphicx,lscape}
\usepackage{tikz,tikz-cd}
\usetikzlibrary{positioning,matrix,calc,shapes.arrows,shadows.blur,shadows,shapes.geometric,shapes.symbols}
\usepackage[square,numbers,sort&compress]{natbib}
\usepackage[hidelinks,hypertexnames=false,hyperfootnotes=false]{hyperref}
\hypersetup{pdftitle={Values of differentials of plane branches: block decomposition and semiroots},
pdfauthor={Marcelo Osnar Rodrigues de Abreu}}
\providecommand{\doi}[1]{\href{https://doi.org/#1}{\nolinkurl{doi:#1}}}

\newtheorem{teorema}{Theorem}[section]
\newtheorem{lema}[teorema]{Lemma}
\newtheorem{proposicao}[teorema]{Proposition}
\newtheorem{corolario}[teorema]{Corollary}
\newtheorem{definicao}[teorema]{Definition}
\newtheorem{exemplo}[teorema]{Example}
\newtheorem{obs}[teorema]{Remark}

\newcommand{\teo}{\begin{teorema}}
\newcommand{\fteo}{\end{teorema}}
\newcommand{\lem}{\begin{lema}}
\newcommand{\flem}{\end{lema}}
\newcommand{\prop}{\begin{proposicao}}
\newcommand{\fprop}{\end{proposicao}}
\newcommand{\cor}{\begin{corolario}}
\newcommand{\fcor}{\end{corolario}}
\newcommand{\defn}{\begin{definicao}}
\newcommand{\fdefn}{\end{definicao}}
\newcommand{\ex}{\begin{exemplo}}
\newcommand{\fex}{\end{exemplo}}
\newcommand{\ob}{\begin{obs}}
\newcommand{\fob}{\end{obs}}
\newcommand{\dem}{\noindent\textsc{Proof:}\enspace}
\newcommand{\f}{\hfill $\rule{2mm}{2mm}$\par\medskip}

\title{Values of differentials of plane branches: block decomposition and semiroots}
\author{Marcelo Osnar Rodrigues de Abreu\\[3pt]
\small Universidade Estadual de Maring\'a\\
\small Maring\'a, Paran\'a, Brazil\\
\small\href{mailto:morabreu@uem.br}{\texttt{morabreu@uem.br}}}
\date{}
\begin{document}
\maketitle
\begin{abstract}
We study the passage of special gaps from a Zariski semiroot to a plane branch. A decomposition into rectangular blocks allows us to compute the maximal valuations of the coefficients of differentials and describe the image of the transition map. We treat vertical and horizontal decompositions, as well as decompositions with alternations, obtaining criteria for the preservation of special gaps. For towers whose ramification indices after the first are all equal to two, the value semigroup of the branch and the special gaps of the first semiroot determine the special gaps of the branch.
\end{abstract}

\section{Introduction}

Let $C_f$ be an irreducible plane branch. Its value semigroup $\Gamma_f$ determines the topological class of $C_f$, whereas the set $\Lambda_f$ of values of differentials is an analytic invariant. The relationship between these two sets plays a central role in the analytic classification of plane branches; for an account of this setting, see \citet{hefez2011analytic,de2022value}.

This paper compares $\Gamma_f$ and $\Lambda_f$ with the corresponding invariants of the Zariski semiroots of $f$, obtained by truncating a Puiseux parametrization previously chosen in the normal form established by \citet{hefez2011analytic}. The truncations follow the construction described in Section~2 of \citet{de2022analytic}. In particular, we study how elements of $\Lambda_k\setminus\Gamma_k$ give rise to elements of $\Lambda_f\setminus\Gamma_f$ and how these special gaps propagate along a tower of semiroots. The combinatorial decomposition considered below has two basic configurations, treated separately in Cases~1 and~2, and a general configuration in which alternations between them may occur.

After the introduction, Section~\ref{sec:preliminares} collects the definitions and auxiliary results. Section~\ref{sec:transicao-semirraizes} studies propagation along the tower and gives a preliminary decomposition of the transition from $f_1$ to $f$; Section~\ref{sec:caso-vertical} treats the vertical decomposition, Section~\ref{sec:caso-horizontal} the horizontal decomposition, and Section~\ref{sec:caso-geral} their simultaneous occurrence and the alternations between the two cases.

\section{Preliminaries}\label{sec:preliminares}

\subsection{Plane branches and value semigroups}

We write $\mathbb N=\{0,1,2,\ldots\}$ and $\mathbb N^*=\mathbb N\setminus\{0\}$. For sets of integers, we use the notation $aS=\{as:s\in S\}$, $a+S=\{a+s:s\in S\}$, and $S+T=\{s+t:s\in S,\ t\in T\}$. The symbols $\dot\cup$ and $\sharp S$ (or $\#S$) denote disjoint union and cardinality, respectively.

Let $\mathbb C\{x,y\}$ be the ring of convergent power series, with maximal ideal $\mathfrak m=\langle x,y\rangle$. A reduced analytic plane curve in $(\mathbb C^2,0)$ is the germ of the zero set of an element $f\in\mathfrak m$ with no repeated irreducible factors. A germ identifies sets that agree in some neighborhood of the origin. If $f$ is irreducible, $C_f=\{f=0\}$ is called a \emph{plane branch}. Its multiplicity is $n=\operatorname{ord}_{(x,y)}f$, that is, the least total degree of a nonzero monomial of $f$, and its local ring is
\[
 \mathcal O_f=\mathbb C\{x,y\}/\langle f\rangle.
\]
The branch is singular if $n>1$ and smooth if $n=1$. Throughout this paper, $C_f$ is singular, although some of its semiroots may be smooth.

Two branches are analytically equivalent if a germ of a biholomorphism of $(\mathbb C^2,0)$ maps one onto the other; they are topologically equivalent if there is a germ of a homeomorphism with this property. An analytic invariant, respectively a topological invariant, is preserved under the corresponding equivalence.

After choosing coordinates in which $x=0$ is transverse to the branch and $y=0$ is tangent to it, the Weierstrass Preparation Theorem allows us to multiply $f$ by a unit and write
\[
 f(x,y)=y^n+c_1(x)y^{n-1}+\cdots+c_n(x),
 \qquad c_i(x)\in x\mathbb C\{x\}.
\]
A monic polynomial of this form is called a \emph{Weierstrass polynomial} in $y$. By the Newton--Puiseux Theorem, $C_f$ admits a convergent primitive parametrization
\[
 \varphi(t)=(x(t),y(t))
 =\left(t^n,\sum_{j>n}a_jt^j\right).
\]
The primitivity condition means that $\operatorname{gcd}\bigl(\{n\}\cup\{j:a_j\neq0\}\bigr)=1$; it rules out reparametrization by a proper power of $t$. By removing the terms that are powers of $x$ preceding the first exponent not divisible by $n$ and normalizing the leading coefficient, we may further assume that
\[
 y(t)=t^{v_1}+\sum_{j>v_1}a_jt^j,\qquad n\nmid v_1.
\]
These coordinate choices and the results used in this subsection are described in \citet[Section~2]{de2022analytic}.

For $z(t)=\sum_jc_jt^j\neq0$, set $\operatorname{ord}_t z=\min\{j:c_j\neq0\}$ and $\operatorname{ord}_t0=\infty$. The homomorphism $\varphi^*(h)=h\circ\varphi$ has kernel $\langle f\rangle$ and identifies $\mathcal O_f$ with $\mathbb C\{x(t),y(t)\}\subset\mathbb C\{t\}$. The valuation of functions is defined by
\[
 \nu_f(h)=\operatorname{ord}_t(h\circ\varphi),
 \qquad h\in\mathbb C\{x,y\}.
\]
Thus, $\nu_f(h)=\infty$ if and only if $h\in\langle f\rangle$. The \emph{value semigroup} of the branch is
\[
 \Gamma_f=\{\nu_f(h):h\in\mathbb C\{x,y\}\setminus\langle f\rangle\}
          =\langle v_0,\ldots,v_g\rangle,
 \qquad v_0=n<v_1<\cdots<v_g,
\]
where $\langle v_0,\ldots,v_g\rangle$ denotes the set of linear combinations with coefficients in $\mathbb N$, and the $v_i$ are its minimal generators. In particular, $\nu_f(x)=v_0$ and $\nu_f(y)=v_1$.

A numerical semigroup is an additive submonoid of $\mathbb N$ with finite complement. The semigroup $\Gamma_f$ is numerical; the elements of $\mathbb N\setminus\Gamma_f$ are its \emph{gaps}, and its \emph{conductor} is
\[
 \mu_f=\min\{c\in\mathbb N:c+\mathbb N\subseteq\Gamma_f\}.
\]
To fix the Puiseux data, set $\beta_0=\varepsilon_0=v_0$, $n_0=1$, and
\[
 \beta_i=\min\{j:a_j\neq0,\ \varepsilon_{i-1}\nmid j\},\qquad
 \varepsilon_i=\operatorname{gcd}(\varepsilon_{i-1},\beta_i),\qquad
 n_i=\frac{\varepsilon_{i-1}}{\varepsilon_i},
 \quad 1\leq i\leq g,
\]
where $g$ is the first index for which $\varepsilon_g=1$. The $\beta_i$ are the \emph{characteristic exponents}; $\beta_1=v_1$, $n_i\geq2$, and
\[
 v_i=n_{i-1}v_{i-1}+\beta_i-\beta_{i-1},
 \qquad \varepsilon_i=\operatorname{gcd}(v_0,\ldots,v_i),
 \qquad 1\leq i\leq g.
\]
The case $g=1$ corresponds to a branch with a single Puiseux pair and semigroup $\langle v_0,v_1\rangle$, where $\operatorname{gcd}(v_0,v_1)=1$.

The following properties will be used below:
\[
 \varepsilon_k\nmid\beta_{k+1}\quad(0\leq k<g),\qquad
 \mu_f=\sum_{i=1}^g(n_i-1)v_i-v_0+1.
\]
Every $z\in\mathbb Z$ has a unique expression
\[
 z=a_0v_0+\cdots+a_gv_g,\qquad
 a_0\in\mathbb Z,\quad 0\leq a_i<n_i\ (1\leq i\leq g),
\]
and $z\in\Gamma_f$ if and only if $a_0\geq0$. Moreover, $\Gamma_f$ is symmetric, that is,
\[
 z\in\Gamma_f\ \Longleftrightarrow\
 \mu_f-1-z\notin\Gamma_f,
 \qquad
 \sharp(\mathbb N\setminus\Gamma_f)=\frac{\mu_f}{2}.
\]
The conductor coincides with the Milnor number $\dim_{\mathbb C}\mathbb C\{x,y\}/\langle f_x,f_y\rangle$, where $f_x$ and $f_y$ are the partial derivatives of $f$. For these properties and their relationship with the topological class of the branch, see \citet[Section~2]{de2022analytic,de2022value}.

\subsection{Values of differentials and special gaps}

Let
\[
 \Omega^1=\mathbb C\{x,y\}\,dx+\mathbb C\{x,y\}\,dy
\]
be the module of holomorphic differential forms. For $\omega=A\,dx+B\,dy\in\Omega^1$, its pullback under the parametrization is
\[
 \varphi^*(\omega)=
 \bigl(A(\varphi(t))x'(t)+B(\varphi(t))y'(t)\bigr)\,dt.
\]
We define
\[
 \nu_f(\omega)=
 \operatorname{ord}_t\bigl(A(\varphi(t))x'(t)
                         +B(\varphi(t))y'(t)\bigr)+1,
\]
with $\nu_f(\omega)=\infty$ when $\varphi^*(\omega)=0$. Thus, the term $ct^{\ell-1}\,dt$, with $c\neq0$, has differential value $\ell$.

The module of Kähler differentials of $\mathcal O_f$ is
\[
 \Omega_{\mathcal O_f/\mathbb C}
 =\frac{\Omega^1}{f\Omega^1+\mathbb C\{x,y\}\,df},
 \qquad df=f_x\,dx+f_y\,dy.
\]
The kernel of the pullback induced on this quotient is its torsion submodule, consisting of the classes annihilated by some nonzero element of $\mathcal O_f$. The classes with nonzero pullback determine the \emph{set of values of differentials}
\[
 \Lambda_f=\{\nu_f(\omega):\omega\in\Omega^1,\ 
                             \varphi^*(\omega)\neq0\}
 \subset\mathbb N^*.
\]
These definitions are independent of the chosen primitive parametrization. For their interpretation in terms of the module of differentials, see \citet[Section~3]{de2022analytic}.

For $h\in\mathfrak m\setminus\langle f\rangle$, we have $\nu_f(dh)=\nu_f(h)$. For $h$ and $\omega$ with finite values,
\[
 \nu_f(h\omega)=\nu_f(h)+\nu_f(\omega).
\]
Consequently,
\[
 \Gamma_f\setminus\{0\}\subseteq\Lambda_f,
 \qquad \Gamma_f+\Lambda_f\subseteq\Lambda_f.
\]
The latter property means that $\Lambda_f$ is a $\Gamma_f$-\emph{semimodule}. It admits a finite system of generators $\lambda_1,\ldots,\lambda_r$:
\[
 \Lambda_f=\bigcup_{i=1}^r(\lambda_i+\Gamma_f).
\]
The system is minimal if none of its elements belongs to the semimodule generated by the others. Finiteness is established in Lemma~4.6 of \citet{hefez2007standard}.

In this paper, the elements of the finite set $\Lambda_f\setminus\Gamma_f$ are called the \emph{special gaps} of $C_f$. The sum of a special gap and an element of $\Gamma_f$ belongs to $\Lambda_f$, but it may belong to $\Gamma_f$. In particular,
\[
 \mathbb N\setminus\Gamma_f
 =(\Lambda_f\setminus\Gamma_f)
  \mathbin{\dot\cup}(\mathbb N^*\setminus\Lambda_f).
\]
If $\Lambda_f\setminus\Gamma_f\neq\varnothing$, then $\min(\Lambda_f\setminus\Gamma_f)>v_0+v_1$ \citep[Remark~2.3]{de2022value}.

\teo[Hefez--Hernandes]\label{teo:forma-normal-hh}
Let $C_f$ be a plane branch with semigroup $\Gamma_f=\langle v_0,\ldots,v_g\rangle$ and set of values of differentials $\Lambda_f$. If $\Lambda_f\setminus\Gamma_f=\varnothing$, then $C_f$ is analytically equivalent to the branch parametrized by $(t^{v_0},t^{v_1})$. Otherwise, set
\[
 \lambda=\min(\Lambda_f\setminus\Gamma_f)-v_0.
\]
Then $C_f$ is analytically equivalent to a branch with parametrization
\begin{equation}\label{eq:forma-normal-hh}
 \varphi(t)=
 \left(t^{v_0},\,
       t^{v_1}+t^\lambda+
       \sum_{\substack{j>\lambda\\j\notin\Lambda_f-v_0}}a_jt^j
 \right).
\end{equation}
\fteo
\dem The existence of this normal form is established in Theorem~2.1 of \citet{hefez2011analytic}; see also Theorem~2.2 of \citet{de2022value}. \f

We first fix a representative of $C_f$ in the normal form of Theorem~\ref{teo:forma-normal-hh}. The semiroots are constructed only after this choice, in the same ambient coordinates. In the nonmonomial case, we retain the notation $\lambda$ from \eqref{eq:forma-normal-hh}.

\subsection{Semiroots and scaling of values}

\defn Given $0\leq k<g$, a monic polynomial $f_k\in\mathbb C\{x\}[y]$ is a \emph{$k$-semiroot} of $f$ if
\[
 \operatorname{deg}_y(f_k)=n_0n_1\cdots n_k=\frac{v_0}{\varepsilon_k},
 \qquad \nu_f(f_k)=v_{k+1}.
\]
We also adopt the convention $f_g=f$, so that $\nu_f(f_g)=v_{g+1}=\infty$. \fdefn

Each $f_k$ is irreducible and defines a branch $C_k:=C_{f_k}$. Its primitive parametrizations can be written in the form $\psi_k(u)=(u^{v_0/\varepsilon_k},y_k(u))$. We denote by $\nu_k:=\nu_{f_k}$ the valuation of functions and forms defined by $\psi_k$, according to the preceding conventions. The sets $\Gamma_k$ and $\Lambda_k$ are, respectively, the value semigroup of functions and the set of values of differentials of $C_k$. The characteristic sequence of $C_k$ is $(\beta_0/\varepsilon_k,\ldots,\beta_k/\varepsilon_k)$ and
\[
 \Gamma_k=
 \left\langle\frac{v_0}{\varepsilon_k},\ldots,\frac{v_k}{\varepsilon_k}\right\rangle,
 \qquad
 \mu_k=\sum_{i=1}^k(n_i-1)\frac{v_i}{\varepsilon_k}
                  -\frac{v_0}{\varepsilon_k}+1.
\]
In particular, $C_0$ is smooth, $\Gamma_0=\mathbb N$, $\Lambda_0=\mathbb N^*$, and $\mu_0=0$. In the chosen coordinates, we may take $f_{-1}=x$ and $f_0=y$. These properties are described in Definition~2.3 and the subsequent discussion in \citet{de2022analytic}.

Throughout this paper, we fix the \emph{Zariski semiroots}, obtained by truncating the parametrization of $f$ previously chosen in the normal form of Theorem~\ref{teo:forma-normal-hh}. This construction by truncations is recalled in Section~2 of \citet{de2022analytic}, immediately before Definition~2.3. For $0\leq k<g$, we consider
\begin{equation}\label{eq:semirraizes-zariski}
 \psi_k(u)=
 \left(u^{v_0/\varepsilon_k},\ 
       \sum_{\substack{v_1\leq j<\beta_{k+1}\\ \varepsilon_k\mid j}}
           a_j u^{j/\varepsilon_k}\right),
 \qquad a_{v_1}=1,
\end{equation}
and take $f_k$ to be the monic minimal equation of this branch. The resulting tower is compatible: for $j<k$, the truncation defining $f_j$ is also a truncation of $f_k$, so that $f_j$ is a $j$-semiroot of $f_k$. The general definition allows other choices; in this paper, we use the Zariski tower fixed by \eqref{eq:semirraizes-zariski}.

\prop[Proposition~2.6 of \citet{de2022analytic}]
\label{prop:escala-funcoes}
If $f_k$ is a $k$-semiroot of $f$, then, for every $h\in\mathbb C\{x\}[y]\setminus\{0\}$ with $\operatorname{deg}_y(h)<\operatorname{deg}_y(f_k)$,
\[
 \nu_f(h)=\varepsilon_k\nu_k(h).
\]
\fprop

We compute the values of $\nu_k$ and the elements of $\Gamma_k$ and $\Lambda_k$ using the primitive parametrization of $C_k$. Reparametrizing by $u=t^{\varepsilon_k}$ multiplies both function and differential values by $\varepsilon_k$: if $\psi_k^*(\omega)=c\,u^{d-1}\,du+\text{higher-order terms}$, then
\[
 \bigl(\psi_k(t^{\varepsilon_k})\bigr)^*(\omega)
 =\varepsilon_kc\,t^{\varepsilon_kd-1}\,dt+\text{higher-order terms}.
\]
This identity concerns the reparametrization of $C_k$; comparison with the parametrization of $C_f$ requires the transition results. For $1\leq k<g$, the conductor formula also gives
\begin{equation}\label{eq:condutor-escala-semirraiz}
 \varepsilon_k(\mu_k-1)+v_0=n_kv_k-\beta_k
                 =v_{k+1}-\beta_{k+1}.
\end{equation}

\subsection{Representatives of differentials and the transition map}

For an irreducible Weierstrass polynomial $q\in\mathbb C\{x\}[y]$ of degree $m$ in $y$, define
\[
 E(q)=\{A\,dx+B\,dy:
 A,B\in\mathbb C\{x\}[y],\
 \operatorname{deg}_y A<m,\
 \operatorname{deg}_y B<m-1\},
\]
with $\operatorname{deg}_y0=-\infty$. By Weierstrass division, every $\omega\in\Omega^1$ admits a unique decomposition
\begin{equation}\label{eq:representante-reduzido}
 \omega=\omega_0+Q\,dq+Pq\,dx,\qquad
 \omega_0\in E(q),\quad P,Q\in\mathbb C\{x,y\};
\end{equation}
in particular, $\omega$ and $\omega_0$ have the same pullback on $C_q$. This is Lemma~3.3 of \citet{de2022analytic}.

A meromorphic form $W\in q^{-1}\Omega^1$ is \emph{logarithmic along $C_q$} if there exist $P,Q\in\mathbb C\{x,y\}$ and $\eta\in\Omega^1$, with $Q\notin\langle q\rangle$, such that $QW=P\,dq/q+\eta$. We denote the module of these forms by $\Omega(\log C_q)$. Thus,
\[
 q\Omega(\log C_q)
 =\{\omega\in\Omega^1:Q\omega=P\,dq+q\eta
       \text{ for some }P,Q,\eta,\ Q\notin\langle q\rangle\}.
\]
For $\omega=A\,dx-B\,dy$, the following equivalences hold:
\[
 \omega\in q\Omega(\log C_q)
 \ \Longleftrightarrow\
 q\mid Aq_y+Bq_x
 \ \Longleftrightarrow\
 \nu_q(\omega)=\infty.
\]
They follow from $q_y\omega=-B\,dq+(Aq_y+Bq_x)\,dx$ and the irreducibility of $q$; see \citet[Section~3]{de2022analytic}, which follows Saito's definition.

Fix $1\leq k<g$ and $\delta_k\in\Lambda_k\setminus\Gamma_k$. Choose a representative
\[
 w_{\delta_k}=A_{\delta_k}\,dx-B_{\delta_k}\,dy
\]
with
\begin{equation}\label{eq:coeficiente-maximal}
 \nu_k(w_{\delta_k})=\delta_k,\qquad
 \nu_k(B_{\delta_k})
 =\max\{\nu_k(B):A,B\in\mathbb C\{x,y\},\
                      \nu_k(A\,dx-B\,dy)=\delta_k\}.
\end{equation}
The values $\nu_k(A)$ and $\nu_k(B)$ are finite: if either pullback were zero, the value of the differential would belong to $\Gamma_k$ or would be infinite. Since the values of $A\,dx$ and $B\,dy$ belong to $\Gamma_k$, a special gap can arise only from cancellation of their leading terms. Therefore,
\[
 \nu_k(A)+\frac{v_0}{\varepsilon_k}
 =\nu_k(B)+\frac{v_1}{\varepsilon_k}<\delta_k.
\]
The values of $B$ considered in \eqref{eq:coeficiente-maximal} thus belong to a nonempty finite set of integers, which guarantees the existence of the maximum.

The representative can be chosen in $E(f_k)$. Indeed, \eqref{eq:representante-reduzido} preserves the pullback, and reducing the coefficient of $dy$ subtracts a multiple of $(f_k)_y$. The identity
\[
 \nu_k((f_k)_y)=\mu_k-1+\frac{v_0}{\varepsilon_k},
\]
recorded in Section~4.2 of \citet{de2022analytic}, together with the inequalities $\nu_k(B)<\delta_k<\mu_k$, shows that this reduction also preserves $\nu_k(B)$: the multiple being subtracted has value greater than $\nu_k(B)$. We therefore fix $w_{\delta_k}\in E(f_k)$ satisfying \eqref{eq:coeficiente-maximal}. The maximum, denoted by $\Theta_k(\delta_k)$ in \citet{de2022analytic}, will be written as $\nu_k(B_{\delta_k})$.

By Lemma~4.3, Remark~4.5, and Proposition~4.7 of \citet{de2022analytic}, we obtain the map
\begin{equation}\label{eq:definicao-rho}
\begin{aligned}
 \rho_k:\Lambda_k\setminus\Gamma_k&\longrightarrow
                    \Lambda_f\setminus\Gamma_f,\\
 \delta_k&\longmapsto\nu_f(w_{\delta_k})=\min\{\varepsilon_k\delta_k,\
       \beta_{k+1}+\varepsilon_k\nu_k(B_{\delta_k})\}=\min\{\varepsilon_k\delta_k,\
       \beta_{k+1}+\nu_f(B_{\delta_k})\}.
\end{aligned}
\end{equation}
The last equality follows from Proposition~\ref{prop:escala-funcoes}. The two quantities being compared are distinct, since $\varepsilon_k\nmid\beta_{k+1}$. The value of $\rho_k$ depends only on $\delta_k$ and the maximum in \eqref{eq:coeficiente-maximal}, although the representative $w_{\delta_k}$ need not be unique.

\lem\label{lem:injetividade-rho}
For $1\leq k<g$, the map $\delta_k\mapsto\nu_k(B_{\delta_k})$ is injective on $\Lambda_k\setminus\Gamma_k$. Consequently, $\rho_k$ is injective.
\flem
\dem Let $\delta'<\delta''$ be two special gaps such that $\nu_k(B_{\delta'})=\nu_k(B_{\delta''})$. There exists $c\in\mathbb C^*$ that cancels the leading coefficients of the pullbacks of these two functions, so that
\[
 \nu_k(B_{\delta'}-cB_{\delta''})>\nu_k(B_{\delta'}).
\]
However, $\nu_k(w_{\delta'}-cw_{\delta''})=\delta'$, contradicting maximality in \eqref{eq:coeficiente-maximal}. Hence the coefficient values are distinct. The images of the two alternatives in \eqref{eq:definicao-rho} belong to distinct residue classes modulo $\varepsilon_k$. In each alternative, injectivity follows, respectively, from multiplication by $\varepsilon_k$ and from the injectivity of $\delta_k\mapsto\nu_k(B_{\delta_k})$. \f

\subsection{Previously established results}

The following results describe the special gaps arising from a semiroot. The first concerns the last extension, with ramification index $n_g=2$.

\teo[Theorem~5.2 of \citet{de2022analytic}]\label{t5.2}
If $g\geq2$ and $n_g=2$, then
\[
\Lambda_f\setminus\Gamma_f=\rho_{g-1}(\Lambda_{g-1}\setminus\Gamma_{g-1})\mathbin{\dot\cup}
\{v_g-2\delta:\delta\in\mathbb N^*\setminus\Lambda_{g-1}\}\mathbin{\dot\cup}\{v_g+2\delta:\delta\in\mathbb N\setminus\Gamma_{g-1}\},
\]
and $\sharp(\Lambda_f\setminus\Gamma_f)=\mu_{g-1}$.
\fteo

For values smaller than $v_{k+1}$, the following description does not require the assumption $n_g=2$.

\prop[Proposition~5.1 of \citet{de2022analytic}]\label{2.2.15}
For $1\leq k<g$, we have
\[
\{\gamma\in\Lambda_f\setminus\Gamma_f:\gamma<v_{k+1}\}=\rho_k(\Lambda_k\setminus\Gamma_k)
 \mathbin{\dot\cup} \{v_{k+1}-\varepsilon_k\delta:\delta\in\mathbb N^*\setminus\Lambda_k\}.
\]
\fprop

\subsection{Comparison of the parametrizations}\label{subsec:comparacao-parametrizacoes}

Assume $g\geq2$. To compare the terms preceding and following $\beta_2$, set
\[
 \psi(t)=\sum_{\substack{v_1\leq j<\beta_2\\\varepsilon_1\mid j}}a_jt^j,
 \qquad
 \varphi(t)=(t^{v_0},y(t)),\qquad
 \varphi_1(t)=(t^{v_0},\psi(t)).
\]
By the definition of $\beta_2$,
\[
 y(t)-\psi(t)=a_{\beta_2}t^{\beta_2}
                +O(t^{\beta_2+1}),\qquad a_{\beta_2}\neq0.
\]
The notation $O(t^N)$ denotes a series of order at least $N$. By convention~\eqref{eq:semirraizes-zariski}, $f_1$ is the Zariski semiroot defined by this truncation; its primitive parametrization is
\[
 \bar\varphi_1(u)=
 \left(u^{v_0/\varepsilon_1},\sum_{\substack{v_1\leq j<\beta_2\\\varepsilon_1\mid j}}
                      a_j u^{j/\varepsilon_1}\right),
 \qquad
 \varphi_1(t)=\bar\varphi_1(t^{\varepsilon_1}).
\]

For $\varphi^*(w)=\sum_{\ell\geq1}c_\ell t^{\ell-1}\,dt$, define
\begin{equation}\label{eq:definicao-E2}
 E_2(w)=\min\{\ell\geq1:c_\ell\neq0,\ \varepsilon_1\nmid\ell\},
\end{equation}
with $E_2(w)=\infty$ if the set is empty. Divisibility refers to the value $\ell$, rather than to the exponent $\ell-1$ of the coefficient of the pullback.

If $H=\sum_{a,b}h_{ab}x^ay^b\in\mathbb C\{x\}[y]$ satisfies $\operatorname{deg}_y H<n_1=v_0/\varepsilon_1$, the weights $av_0+bv_1$ of its nonzero monomials are distinct: equality of two such weights would imply $n_1\mid(b-b')$, with $0\leq b,b'<n_1$. For $b\geq1$, the binomial expansion gives
\[
 (x^ay^b)(\varphi)-(x^ay^b)(\varphi_1)
 =b\,a_{\beta_2}t^{av_0+(b-1)v_1+\beta_2}
       +O(t^{av_0+(b-1)v_1+\beta_2+1}).
\]
For $b=0$, the difference is zero. Therefore, if $H\notin\mathbb C\{x\}$,
\begin{equation}\label{eq:ordem-diferenca-funcoes}
 \operatorname{ord}_t(H\circ\varphi-H\circ\varphi_1)
 =\beta_2-v_1+
    \min_{\substack{h_{ab}\neq0\\b\geq1}}(av_0+bv_1).
\end{equation}
The minimum is taken only over monomials that depend on $y$.

\lem\label{lem:primeiro-valor-nao-divisivel}
Let $w=A\,dx+B\,dy\in E(f_1)$, with $A,B\neq0$, satisfy
\[
 \nu_f(A)+v_0=\nu_f(B)+v_1<\nu_f(w).
\]
Then
\[
 E_2(w)=\beta_2+\nu_f(B).
\]
In particular, the conclusion applies to $w\in E(f_1)$ with $\nu_1(w)\in\Lambda_1\setminus\Gamma_1$.
\flem
\dem By the degree bounds in $E(f_1)$, the leading monomials of $A$ and $B$ are unique. Equality of their differential values forces these monomials to be, respectively,
\[
 c v_1x^ay^b
 \quad\text{and}\quad
 -c v_0x^{a+1}y^{b-1},
 \qquad c\in\mathbb C^*,\quad a\geq0,\quad 1\leq b<n_1.
\]
Indeed, if $x^{a'}y^{b'}$ is the leading monomial of $B$, the equality $(a+1)v_0+bv_1=a'v_0+(b'+1)v_1$, with $0\leq b<n_1$ and $0\leq b'<n_1-1$, implies $b=b'+1$ and $a'=a+1$. The cancellation assumption determines the ratio of their coefficients.

After writing the pullbacks as series multiplied by $dt$, the difference of their coefficients, multiplied by $t$, satisfies
\[
\begin{aligned}
 \frac{t\,(\varphi^*w-\varphi_1^*w)}{dt}
 &=c v_0a_{\beta_2}
    \bigl(v_1b-(b-1)v_1-\beta_2\bigr)
       t^{(a+1)v_0+(b-1)v_1+\beta_2}\\
 &\quad+O(t^{(a+1)v_0+(b-1)v_1+\beta_2+1})\\
 &=c v_0a_{\beta_2}(v_1-\beta_2)
       t^{\nu_f(B)+\beta_2}
       +O(t^{\nu_f(B)+\beta_2+1}).
\end{aligned}
\]
The remaining monomials contribute only at higher orders, by the same expansion and the uniqueness of the weights. The displayed coefficient is nonzero, since $\beta_2>v_1$. On the other hand, $t\,\varphi_1^*w/dt$ contains only powers divisible by $\varepsilon_1$, whereas $\varepsilon_1\nmid\nu_f(B)+\beta_2$. This proves the equality for $E_2(w)$. Finally, if $\nu_1(w)\in\Lambda_1\setminus\Gamma_1$, its leading terms cancel on $C_1$; Proposition~\ref{prop:escala-funcoes} and equality of the leading terms of the parametrizations ensure the same cancellation on $C_f$. \f

\prop\label{prop:lambda-semirraiz-normalizada}
Assume $g\geq2$, with $f$ in the normal form \eqref{eq:forma-normal-hh} and $f_1$ its Zariski semiroot defined by \eqref{eq:semirraizes-zariski}. Then
\[
 \Lambda_1\setminus\Gamma_1\neq\varnothing
 \quad\Longleftrightarrow\quad \lambda<\beta_2.
\]
In this case,
\begin{equation}\label{eq:lambda-escalas}
 \varepsilon_1\min(\Lambda_1\setminus\Gamma_1)
 =\lambda+v_0
 =\min(\Lambda_f\setminus\Gamma_f).
\end{equation}
\fprop
\dem Since $g\geq2$, the branch $C_f$ is not analytically equivalent to a monomial branch; by Theorem~\ref{teo:forma-normal-hh}, $\Lambda_f\setminus\Gamma_f\neq\varnothing$. In the normal form, $\lambda$ is the first exponent after $v_1$ with nonzero coefficient. Consequently, $\lambda\leq\beta_2$. If $\lambda=\beta_2$, the truncation defining $f_1$ is $(u^{v_0/\varepsilon_1},u^{v_1/\varepsilon_1})$ and $\Lambda_1\setminus\Gamma_1=\varnothing$.

Assume $\lambda<\beta_2$. By the definition of $\beta_2$, $\varepsilon_1\mid\lambda$. The differential $w_0=v_1y\,dx-v_0x\,dy$ satisfies
\[
 \nu_1(w_0)=\frac{\lambda+v_0}{\varepsilon_1},
 \qquad
 \nu_f(w_0)=\lambda+v_0.
\]
Moreover, $(\lambda+v_0)/\varepsilon_1\notin\Gamma_1$, since $\varepsilon_1\Gamma_1=\langle v_0,v_1\rangle\subseteq\Gamma_f$ and $\lambda+v_0\notin\Gamma_f$. Therefore, $\Lambda_1\setminus\Gamma_1\neq\varnothing$. For every $\eta\in\Lambda_1\setminus\Gamma_1$, the map \eqref{eq:definicao-rho} gives
\[
 \lambda+v_0
 =\min(\Lambda_f\setminus\Gamma_f)
 \leq\rho_1(\eta)\leq \varepsilon_1\eta.
\]
Since $(\lambda+v_0)/\varepsilon_1\in\Lambda_1\setminus\Gamma_1$, \eqref{eq:lambda-escalas} follows. \f

Thus, when $\Lambda_1\setminus\Gamma_1\neq\varnothing$, $\lambda/\varepsilon_1$ is the Zariski invariant of $C_1$ and $\lambda<\beta_2$. The parametrization of $f$ can then be written as
\begin{equation}\label{eq:normalizacao-lambda}
 y(t)=t^{v_1}+t^\lambda+
        \sum_{\lambda<j<\beta_2}a_jt^j
        +a_{\beta_2}t^{\beta_2}+\cdots,
\end{equation}
where the coefficients absent from the normal form are zero. The abbreviations $(h.o.t.)$ and $(i.o.t.)$ denote, respectively, higher-order terms and intermediate-order terms, specified in each expression.

\section{Transition between semiroots}\label{sec:transicao-semirraizes}

We retain the representatives $w_{\delta_k}\in E(f_k)$ fixed in \eqref{eq:coeficiente-maximal} and the map
\[
 \rho_k(\delta_k)
 =\min\{\varepsilon_k\delta_k,\,
           \beta_{k+1}+\nu_f(B_{\delta_k})\}.
\]
The passage from $f_k$ to $f_{k+1}$ is described by the same formula, divided by $\varepsilon_{k+1}$. Thus, the values of $\rho_k$ are on the scale of $f$, whereas $\rho_k(\delta_k)/\varepsilon_{k+1}$ is on the scale of $f_{k+1}$.

\subsection{Propagation along the tower}\label{subsec:propagacao-torre}

The results in this subsection describe the propagation of a known special gap of a semiroot to the subsequent stages of the tower. Here, the index $k$ is arbitrary, and the semigroup $\Gamma_k$ may have more than two generators.

\teo\label{teo:propagacao-primeira-alternativa}
Let $g\geq3$, $1\leq k\leq g-2$, and $\delta_k\in\Lambda_k\setminus\Gamma_k$. If
\[
 \delta_{k+1}:=\frac{\rho_k(\delta_k)}{\varepsilon_{k+1}}
 =\frac{\varepsilon_k}{\varepsilon_{k+1}}\delta_k,
\]
then $\delta_{k+1}\in\Lambda_{k+1}\setminus\Gamma_{k+1}$ and
\[
 \rho_{k+1}(\delta_{k+1})
 =\varepsilon_{k+1}\delta_{k+1}=\varepsilon_k\delta_k.
\]
\fteo

\dem Applying \eqref{eq:definicao-rho} to the branch $f_{k+1}$ ensures that $\delta_{k+1}$ is a special gap and that $\nu_{k+1}(w_{\delta_k})=\delta_{k+1}$. The two terms in the minimum are distinct, since $\varepsilon_k\nmid\beta_{k+1}$. Hence,
\[
 \varepsilon_{k+1}\delta_{k+1}
 <\beta_{k+1}+\varepsilon_k\nu_k(B_{\delta_k}).
\]
Since $w_{\delta_k}\in E(f_k)\subset E(f_{k+1})$, the maximality of $\nu_{k+1}(B_{\delta_{k+1}})$ and Proposition~\ref{prop:escala-funcoes} give
\[
 \varepsilon_k\nu_k(B_{\delta_k})
 =\varepsilon_{k+1}\nu_{k+1}(B_{\delta_k})
 \leq\varepsilon_{k+1}\nu_{k+1}(B_{\delta_{k+1}}).
\]
Therefore,
\[
 \varepsilon_{k+1}\delta_{k+1}
 <\beta_{k+2}+
       \varepsilon_{k+1}\nu_{k+1}(B_{\delta_{k+1}}),
\]
and the formula for $\rho_{k+1}$ completes the proof.\f

\teo\label{teo:propagacao-segunda-alternativa}
Let $g\geq3$, $1\leq k\leq g-2$, and $\delta_k\in\Lambda_k\setminus\Gamma_k$. If
\[
 \delta_{k+1}:=\frac{\rho_k(\delta_k)}{\varepsilon_{k+1}}
 =\frac{\beta_{k+1}+\varepsilon_k\nu_k(B_{\delta_k})}
        {\varepsilon_{k+1}},
\]
then $\delta_{k+1}\in\Lambda_{k+1}\setminus\Gamma_{k+1}$ and
\[
 \rho_{k+1}(\delta_{k+1})
 =\varepsilon_{k+1}\delta_{k+1}.
\]
\fteo

\dem Again, $\nu_{k+1}(w_{\delta_k})=\delta_{k+1}$ and $\delta_{k+1}\in\Lambda_{k+1}\setminus\Gamma_{k+1}$ by the transition formula applied to $f_{k+1}$. Comparing the coefficients of $dy$ as in the preceding proof gives
\[
 \begin{aligned}
 \varepsilon_{k+1}\delta_{k+1}
 &=\beta_{k+1}+\varepsilon_k\nu_k(B_{\delta_k})\\
 &\leq\beta_{k+1}
       +\varepsilon_{k+1}\nu_{k+1}(B_{\delta_{k+1}})\\
 &<\beta_{k+2}
       +\varepsilon_{k+1}\nu_{k+1}(B_{\delta_{k+1}}).
 \end{aligned}
\]
The conclusion follows from \eqref{eq:definicao-rho}.\f

\teo\label{teo:propagacao-lacunas-novas}
Let $g\geq3$, $1\leq k\leq g-2$, and $\delta_{k+1}\in\Lambda_{k+1}\setminus\Gamma_{k+1}$ satisfy
\[
 \delta_{k+1}\notin
 \frac{1}{\varepsilon_{k+1}}\rho_k(\Lambda_k\setminus\Gamma_k),
 \qquad
 \delta_{k+1}<\frac{v_{k+1}}{\varepsilon_{k+1}}.
\]
Then $\rho_{k+1}(\delta_{k+1})=\varepsilon_{k+1}\delta_{k+1}$.
\fteo

\dem Proposition~\ref{2.2.15}, applied to $f_{k+1}$, gives
\[
 \delta_{k+1}
 =\frac{v_{k+1}}{\varepsilon_{k+1}}
  -\frac{\varepsilon_k}{\varepsilon_{k+1}}\alpha,
 \qquad \alpha\in\mathbb N^*\setminus\Lambda_k.
\]
By Corollary~4.2 of \citet{de2022analytic}, there exists $\omega=A\,dx-B\,dy\in E(f_k)\cap f_k\Omega(\log C_k)$ with $\nu_{k+1}(\omega)=\delta_{k+1}$. Proposition~4.1 of the same paper and \eqref{eq:condutor-escala-semirraiz} imply
\[
 \begin{aligned}
 \varepsilon_{k+1}\delta_{k+1}
 &=v_{k+1}-\varepsilon_k(\mu_k-1)-v_0
       +\varepsilon_k\nu_k(B)\\
 &=\beta_{k+1}+\varepsilon_{k+1}\nu_{k+1}(B)\\
 &\leq\beta_{k+1}
       +\varepsilon_{k+1}\nu_{k+1}(B_{\delta_{k+1}})\\
 &<\beta_{k+2}
       +\varepsilon_{k+1}\nu_{k+1}(B_{\delta_{k+1}}).
 \end{aligned}
\]
The change of scale for the coefficient $B$ is justified by Proposition~\ref{prop:escala-funcoes}. The definition of $\rho_{k+1}$ completes the proof.\f

\cor\label{cor:propagacao-lacunas-novas}
Let $g\geq3$, $1\leq k\leq g-2$, and $\delta_{k+1}\in\Lambda_{k+1}\setminus\Gamma_{k+1}$, with $\delta_{k+1}<v_{k+1}/\varepsilon_{k+1}$. Then $\varepsilon_{k+1}\delta_{k+1}\in\Lambda_f\setminus\Gamma_f$.
\fcor

\dem If $\delta_{k+1}$ belongs to the image of the transition from $f_k$ to $f_{k+1}$, apply Theorems~\ref{teo:propagacao-primeira-alternativa} and \ref{teo:propagacao-segunda-alternativa}. Otherwise, apply Theorem~\ref{teo:propagacao-lacunas-novas}. In both cases, $\rho_{k+1}(\delta_{k+1})=\varepsilon_{k+1}\delta_{k+1}$.\f

\lem\label{lem:propagacao-total-binaria}
If $2\leq k<g$ and $n_k=2$, then
\[
 \rho_k(\delta)=\varepsilon_k\delta
 \qquad\text{for every }\delta\in\Lambda_k\setminus\Gamma_k.
\]
\flem
\dem By Theorem~\ref{t5.2}, applied to $f_k$, each $\delta\in\Lambda_k\setminus\Gamma_k$ belongs to the image of the transition from $f_{k-1}$ to $f_k$ or has the form
\[
 \delta=\frac{v_k}{\varepsilon_k}-2\alpha,
 \qquad
 \alpha\in(\mathbb N^*\setminus\Lambda_{k-1})
       \cup\bigl(-(\mathbb N\setminus\Gamma_{k-1})\bigr).
\]
In the first case, the conclusion follows from Theorems~\ref{teo:propagacao-primeira-alternativa} and \ref{teo:propagacao-segunda-alternativa}, with $k-1$ in place of $k$. In the second, Corollary~4.2 of \citet{de2022analytic} provides a differential
\[
 \omega=A\,dx-B\,dy
 \in E(f_{k-1})\cap f_{k-1}\Omega(\log C_{k-1}),
 \qquad \nu_k(\omega)=\delta.
\]
By Proposition~4.1 of the same paper and \eqref{eq:condutor-escala-semirraiz},
\[
 \begin{aligned}
 \varepsilon_k\delta
 &=v_k-\varepsilon_{k-1}(\mu_{k-1}-1)-v_0
       +\varepsilon_{k-1}\nu_{k-1}(B)\\
 &=\beta_k+\varepsilon_k\nu_k(B)\\
 &\leq\beta_k+\varepsilon_k\nu_k(B_\delta)
 <\beta_{k+1}+\varepsilon_k\nu_k(B_\delta).
 \end{aligned}
\]
Therefore, $\rho_k(\delta)=\varepsilon_k\delta$. The argument includes $\alpha<0$, that is, values greater than $v_k/\varepsilon_k$.\f

\subsection{The transition from the first semiroot to the branch}
\label{subsec:primeira-transicao}

We now turn to the central problem of this paper: determining what happens to the special gaps of the first semiroot in the transition to a branch with more generators. Since
\[
 \Gamma_1=\langle n_1,m\rangle,\qquad
 n_1=\frac{v_0}{\varepsilon_1},\qquad
 m=\frac{v_1}{\varepsilon_1},
\]
$C_1$ is a branch with a single Puiseux pair; on the scale of $f$, the corresponding semigroup is $\varepsilon_1\Gamma_1=\langle v_0,v_1\rangle$. For $\delta\in\Lambda_1\setminus\Gamma_1$, the transition formula is
\begin{equation}\label{eq:rho-primeira-transicao}
 \rho_1(\delta)=
 \min\{\varepsilon_1\delta,\,
       \beta_2+\varepsilon_1\nu_1(B_\delta)\}.
\end{equation}
The first alternative means that the gap $\delta$ remains special after rescaling to $f$. In the second, its value is replaced by a gap congruent to $\beta_2$ modulo $\varepsilon_1$. Thus, the problem is to decide which alternative holds simultaneously for all elements of $\Lambda_1\setminus\Gamma_1$.

Proposition~\ref{2.2.15} separates the special gaps of $f$ smaller than $v_2$ into two families. The family
\[
 \{v_2-\varepsilon_1\alpha:
       \alpha\in\mathbb N^*\setminus\Lambda_1\}
\]
is determined by $\Lambda_1$; it therefore remains to compute $\rho_1(\Lambda_1\setminus\Gamma_1)$. To avoid analyzing each gap individually, we decompose $\Lambda_1\setminus\Gamma_1$ into blocks $A_i$, defined from the minimal generators of this $\Gamma_1$-semimodule. Multiplication of differentials by monomials in $x$ and $y$ allows us to carry the minimum comparison in~\eqref{eq:rho-primeira-transicao} over to the other elements of the same block.

If
\[
 \delta_i=s_i m-r_i n_1
\]
is a minimal generator, the geometry of the corresponding block depends on the ordering of the pairs $(s_i,r_i)$. Simultaneous increase of both coordinates yields the vertical decomposition of Case~1; simultaneous decrease yields the horizontal decomposition of Case~2; changes between these behaviors give rise to the general case. For each configuration, we proceed in four steps: describe the blocks $A_i$, compute $\nu_1(B_\delta)$ in each block, determine $\rho_1(A_i)$, and derive preservation and propagation criteria.

All sets $A_i$ are considered on the scale of $f_1$. The sets $\varepsilon_1A_i$ and the alternative blocks $B_i$ are considered on the scale of $f$. We retain this convention in the following sections.

Every special gap of $f_1$ is greater than $v_1/\varepsilon_1$. Moreover, the least such gap always remains special after rescaling. Indeed, by Proposition~\ref{prop:lambda-semirraiz-normalizada}, if $\Lambda_1\setminus\Gamma_1\neq\varnothing$, then
\[
 \nu_1(v_1y\,dx-v_0x\,dy)=\frac{\lambda+v_0}{\varepsilon_1},
 \qquad
 \nu_f(v_1y\,dx-v_0x\,dy)=\lambda+v_0.
\]
Consequently,
\[
 [\varepsilon_1(\Lambda_1\setminus\Gamma_1)]
 \cap(\Lambda_f\setminus\Gamma_f)=\varnothing
 \quad\Longleftrightarrow\quad
 \Lambda_1\setminus\Gamma_1=\varnothing.
\]
By Theorem~\ref{teo:forma-normal-hh}, the last condition is equivalent to the analytic equivalence of $f_1$ with the monomial branch $(t^{v_0/\varepsilon_1},t^{v_1/\varepsilon_1})$.

\subsection{Motivating examples}

For the following examples, consider $f_1=(t^4,t^9+t^{10})$, whose special gaps are $\{14,19,23\}$.

\begin{center}
	\begin{tabular}{r|cccccc}
	27 & \cellcolor{green!25}23 & \cellcolor{green!25}19 & 15 & 11 & 7 & 3 \\
	18 & \cellcolor{green!25}14 & 10 & 6 & 2 \\
	9 & 5 & 1 \\
	0 &  \\
\end{tabular}
\end{center}

\ex\label{ex1} Consider $f=(t^8,t^{18}+t^{20}+t^{23})$.
\begin{center}
\begin{tabular}{r|cccccccccccccccc}
	131 & \cellcolor{green!25}123 & \cellcolor{green!25}115 & \cellcolor{green!25}107 & \cellcolor{green!25}99 & \cellcolor{green!25}91 & \cellcolor{green!25}83 & \cellcolor{green!25}75 & \cellcolor{green!25}67 & 59 & 51 & 43 & 35 & 27 & 19 & 11 & 3 \\
	113 & \cellcolor{green!25}105 & \cellcolor{green!25}97 & \cellcolor{green!25}89 & \cellcolor{green!25}81 & \cellcolor{green!25}73 & \cellcolor{green!25}65 & \cellcolor{green!25}57 & 49 & 41 & 33 & 25 & 17 & 9 & 1 \\
	95 & \cellcolor{green!25}87 & \cellcolor{green!25}79 & \cellcolor{green!25}71 & \cellcolor{green!25}63 & \cellcolor{green!25}55 & \cellcolor{green!25}47 & 39 & 31 & 23 & 15 & 7 \\
	77 & 69 & 61 & 53 & 45 & 37 & 29 & 21 & 13 & 5 \\
	54 & \cellcolor{green!25}46 & \cellcolor{green!25}38 & 30 & 22 & 14 & 6 \\
	36 & \cellcolor{green!25}28 & 20 & 12 & 4 \\
	18 & 10 & 2 \\
	0 &  \\
\end{tabular}
\end{center}

The values $2\cdot14$, $2\cdot19$, and $2\cdot23$ are special gaps of $f$. Since $n_2=2$, Theorem~\ref{t5.2} gives $\sharp(\Lambda_f\setminus\Gamma_f)=\mu_1=(4-1)(9-1)=24$.
\fex

\ex Consider $f=(t^8,t^{18}+t^{20}+t^{21})$.
\begin{center}
\begin{tabular}{r|cccccccccccccccc}
	129 & \cellcolor{green!25}121 & \cellcolor{green!25}113 & \cellcolor{green!25}105 & \cellcolor{green!25}97 & \cellcolor{green!25}89 & \cellcolor{green!25}81 & \cellcolor{green!25}73 & \cellcolor{green!25}65 & 57 & 49 & 41 & 33 & 25 & 17 & 9 & 1 \\
	111 & \cellcolor{green!25}103 & \cellcolor{green!25}95 & \cellcolor{green!25}87 & \cellcolor{green!25}79 & \cellcolor{green!25}71 & \cellcolor{green!25}63 & \cellcolor{green!25}55 & 47 & 39 & 31 & 23 & 15 & 7 \\
	93 & \cellcolor{green!25}85 & \cellcolor{green!25}77 & \cellcolor{green!25}69 & \cellcolor{green!25}61 & \cellcolor{green!25}53 & \cellcolor{green!25}45 & \cellcolor{green!25}37 & 29 & 21 & 13 & 5 \\
	75 & 67 & 59 & 51 & 43 & 35 & 27 & 19 & 11 & 3 \\
	54 & \cellcolor{green!25}46 & 38 & 30 & 22 & 14 & 6 \\
	36 & \cellcolor{green!25}28 & 20 & 12 & 4 \\
	18 & 10 & 2 \\
	0 &  \\
\end{tabular}
\end{center}

\fex

The difference between the examples occurs at the minimal generator $19$:
\[
 \begin{array}{c|c|c}
  \beta_2&2\cdot14\ \text{and}\ \beta_2+v_0
         &2\cdot19\ \text{and}\ \beta_2+2v_0\\ \hline
  23&28<31&38<39\\
  21&28<29&38>37
 \end{array}
\]
Thus, the first comparison remains unchanged, whereas the second reverses direction. The following results justify using these comparisons to determine the blocks in the image of $\rho_1$.

\subsection{Common tools for the block decomposition}

Before distinguishing the vertical and horizontal configurations, we establish two common properties: the duality that determines the set of maximal values of the coefficients of $dy$, and a formula for these values in any partition into rectangular blocks.

\lem\label{lem:dualidade-coeficientes}
Write $n_1=v_0/\varepsilon_1$, $m=v_1/\varepsilon_1$, and $K=(n_1-1)m$. Then
\begin{equation}\label{eq:dualidade-coeficientes}
 \{\nu_1(B_\delta):\delta\in\Lambda_1\setminus\Gamma_1\}
 =K-(\Lambda_1\setminus\Gamma_1).
\end{equation}
Moreover,
\begin{equation}\label{eq:monotonia-coeficientes}
 \nu_1(B_{\delta+\gamma})\geq\nu_1(B_\delta)+\gamma
\end{equation}
whenever $\gamma\in\Gamma_1$ and $\delta,\delta+\gamma\in\Lambda_1\setminus\Gamma_1$.
\flem

\dem All values in this proof are on the scale of $f_1$. Let $F=\mu_1-1=n_1m-n_1-m$, so that $K=F+n_1$. Fix $\delta\in\Lambda_1\setminus\Gamma_1$ and a representative $\omega=P\,dx-Q\,dy$ with $\nu_1(\omega)=\delta$ and $\nu_1(Q)=\nu_1(B_\delta)=\theta$. Since $\delta\notin\Gamma_1$, the leading terms of the two summands cancel and
\[
 \nu_1(P)+n_1=\theta+m<\delta\leq F.
\]
Write $\theta=an_1+bm$, with $a\geq0$ and $0\leq b<n_1$. If $a=0$, the inequality $\theta+m<F$ implies $b+1<n_1$, but $\nu_1(P)=(b+1)m-n_1\notin\Gamma_1$, a contradiction. Thus, $\theta-n_1\in\Gamma_1$. By the symmetry of the semigroup,
\[
 \alpha:=K-\theta=F-(\theta-n_1)\notin\Gamma_1,
 \qquad \alpha>0.
\]

The duality of logarithmic residues, presented in Theorem~3.8 of \citet{de2022analytic}, states that
\[
 \ell\in\Lambda_1\quad\Longleftrightarrow\quad-\ell\notin\Delta_1,
\]
where
\[
 \Delta_1=
 \left\{\nu_1(Q_0)-K:
 P_0\,dx-Q_0\,dy\in f_1\Omega(\log C_1),\
 Q_0\not\equiv0\pmod{f_1}\right\}.
\]
Indeed, $\nu_1(\partial f_1/\partial y)=K$, and the residue of $P_0\,dx-Q_0\,dy$ is $-Q_0/(\partial f_1/\partial y)$. If $\alpha\notin\Lambda_1$, this duality gives $\eta=P_0\,dx-Q_0\,dy\in f_1\Omega(\log C_1)$ with $\nu_1(Q_0)=K-\alpha=\theta$. Its pullback to $C_1$ is zero. Hence there exists $c\in\mathbb C$ such that
\[
 \nu_1(\omega-c\eta)=\delta,\qquad
 \nu_1(Q-cQ_0)>\theta,
\]
contradicting the maximality of $\nu_1(B_\delta)$. Therefore, $\alpha\in\Lambda_1\setminus\Gamma_1$, which proves one inclusion in \eqref{eq:dualidade-coeficientes}. Equality follows from the finiteness of $\Lambda_1\setminus\Gamma_1$ and the injectivity of $\delta\mapsto\nu_1(B_\delta)$, established in Lemma~\ref{lem:injetividade-rho}.

To prove \eqref{eq:monotonia-coeficientes}, choose a function $h$ with $\nu_1(h)=\gamma$. The differential $h\,w_\delta$ has value $\delta+\gamma$, and its coefficient of $dy$ has value $\nu_1(B_\delta)+\gamma$. The definition of the maximum completes the proof.\f

To define the blocks, assume $\Lambda_1\setminus\Gamma_1\neq\emptyset$. By Lemma~4.6 of \citet{hefez2007standard}, the $\Gamma_1$-semimodule $\Lambda_1$ admits a finite system of minimal generators. Since $\Gamma_1=\langle n_1,m\rangle$, this system has the form
\[
 \{n_1,m,\delta_0,\ldots,\delta_t\},
 \qquad \delta_0<\cdots<\delta_t,\qquad
 \delta_i\in\Lambda_1\setminus\Gamma_1.
\]
Consequently, every special gap belongs to some $\delta_i+\Gamma_1$. Minimality implies $\delta_j-\delta_i\notin\Gamma_1$ for $i<j$; in particular, the $\delta_i$ belong to distinct residue classes modulo $n_1$ and modulo $m$. Write
\[
 \delta_i=s_i m-r_i n_1,\qquad
 2\leq s_i<n_1,\qquad r_i\geq1.
\]
The integers $s_i$ and $r_i$ are pairwise distinct within each family. Moreover, by Proposition~\ref{prop:lambda-semirraiz-normalizada},
\[
 \delta_0=\min(\Lambda_1\setminus\Gamma_1)
          =\frac{\lambda+v_0}{\varepsilon_1}.
\]

The ordering of the generators defines the partition
\begin{equation}\label{eq:particao-geradores}
 A_i=\bigl[(\Lambda_1\setminus\Gamma_1)\cap
                (\delta_i+\Gamma_1)\bigr]
       \setminus\bigcup_{h<i}A_h,\qquad 0\leq i\leq t.
\end{equation}
Thus, $\min A_i=\delta_i$ and $\Lambda_1\setminus\Gamma_1=\dot\bigcup_{i=0}^t A_i$. These sets are on the scale of $f_1$; on the scale of $f$, the corresponding blocks are $\varepsilon_1A_i$, with minimum $\varepsilon_1\delta_i$. The geometric shape of each block depends on the ordering of the pairs $(s_i,r_i)$ and will be described in the following sections.

The next lemma determines the maximal coefficients in a rectangular partition. Its statement also applies to subblocks, provided they are ordered by their respective minima.

\lem\label{lem:coeficientes-particao-retangular}
Let $A_0,\ldots,A_t$ be a partition of $\Lambda_1\setminus\Gamma_1$ of the form
\[
 A_i=\{\delta_i+a m+b n_1:0\leq a\leq p_i,\ 0\leq b\leq q_i\},
 \qquad p_i,q_i\in\mathbb N.
\]
Assume that $\delta_i=\min A_i$ is strictly increasing and that $M_i=\max A_i$ is strictly decreasing. Then, for $\delta\in A_i$,
\begin{equation}\label{eq:coeficientes-particao-retangular}
 \nu_1(B_\delta)=K-M_i+\delta-\delta_i,
 \qquad K=(n_1-1)m.
\end{equation}
\flem

\dem The rectangular shape gives
\[
 M_i=\delta_i+p_i m+q_i n_1,\qquad
 A_i=\delta_i+M_i-A_i.
\]
In particular, $\delta-\delta_i,M_i-\delta\in\Gamma_1$ for $\delta\in A_i$.

The proof is by induction on $i$. By the duality~\eqref{eq:dualidade-coeficientes}, the coefficient values form the set $\dot\bigcup_{h=0}^t(K-A_h)$. Suppose that \eqref{eq:coeficientes-particao-retangular} has already been proved for $h<i$. The symmetry of $A_h$ shows that the image of this block is $K-A_h$. By the injectivity in Lemma~\ref{lem:injetividade-rho}, the coefficient values on the remaining blocks therefore form $\dot\bigcup_{h=i}^t(K-A_h)$. Its minimum is $K-M_i$ and its maximum is $K-\delta_i$, by the ordering of the endpoints. Monotonicity~\eqref{eq:monotonia-coeficientes} implies
\[
 K-\delta_i
 \geq\nu_1(B_{M_i})
 \geq\nu_1(B_{\delta_i})+M_i-\delta_i
 \geq K-M_i+M_i-\delta_i=K-\delta_i.
\]
All inequalities are equalities. Applying \eqref{eq:monotonia-coeficientes} once again, now between $\delta_i$, $\delta$, and $M_i$, we obtain
\[
 K-M_i+\delta-\delta_i
 \leq\nu_1(B_\delta)
 \leq\nu_1(B_{M_i})-(M_i-\delta)
 =K-M_i+\delta-\delta_i.
\]
This proves the formula and completes the induction. \f

By Proposition~\ref{2.2.15},
$$
\left\{\gamma \in \Lambda_f\setminus\Gamma_f ; \gamma<v_{2}\right\}= \rho_1\left(\Lambda_1 \backslash \Gamma_1\right) \dot{\cup}\left\{v_{2}-\varepsilon_{1} \delta_1 ; \delta_1 \in \mathbb{N}^* \backslash \Lambda_1\right\} .
$$

For $1\leq s<n_1$, let $\ell_s=s m-r_s n_1$ be the least element of $\Lambda_1$ in the class $s m$ modulo $n_1$. Since $s m\in\Gamma_1\setminus\{0\}$, we have $r_s\geq0$, and this class contains exactly $r_s$ special gaps. If $\alpha_s=\ell_s-n_1>0$, then $\alpha_s$ is the greatest element of $\mathbb N^*\setminus\Lambda_1$ in this class. The reflection $\alpha\mapsto v_2-\varepsilon_1\alpha$ associates it with the class $v_2-sv_1$ modulo $v_0$, whose least element arising from this family is
\[
 v_2-\varepsilon_1\alpha_s
 =v_2-sv_1+(r_s+1)v_0
 =v_2+(n_1-s)v_1-(m-r_s-1)v_0.
\]
Distinct classes modulo $n_1$ give rise to distinct classes modulo $v_0$: the congruence $\varepsilon_1\alpha\equiv\varepsilon_1\alpha'\pmod{v_0}$ is equivalent to $\alpha\equiv\alpha'\pmod{n_1}$.

The decomposition also gives the number of special gaps below $v_2$. By the injectivity of $\rho_1$ and the decomposition in Proposition~\ref{2.2.15},
\[
 \begin{aligned}
 \#\{\gamma\in\Lambda_f\setminus\Gamma_f:\gamma<v_2\}
 &=\#(\Lambda_1\setminus\Gamma_1)
   +\#(\mathbb N^*\setminus\Lambda_1)\\
 &=\#(\mathbb N\setminus\Gamma_1)=\frac{\mu_1}{2},
 \end{aligned}
\]
where $\mu_1=(v_0/\varepsilon_1-1)(v_1/\varepsilon_1-1)$ is the conductor of $\Gamma_1$. If $g=2$ and $n_2=2$, Theorem~\ref{t5.2} gives $\#(\Lambda_f\setminus\Gamma_f)=\mu_1$; the other $\mu_1/2$ elements belong to the family $\{v_2+2\eta:\eta\in\mathbb N\setminus\Gamma_1\}$, which lies above $v_2$.

\section{Case~1: vertical decomposition}\label{sec:caso-vertical}

\subsection{Construction of the vertical blocks}

Throughout this section, we assume $g\geq2$. If $\Lambda_1\setminus\Gamma_1=\emptyset$, then $\rho_1(\Lambda_1\setminus\Gamma_1)=\emptyset$, and the inclusion $\varepsilon_1(\Lambda_1\setminus\Gamma_1)\subset\Lambda_f\setminus\Gamma_f$ is immediate. We therefore assume that $\Lambda_1\setminus\Gamma_1\neq\emptyset$ and that its minimal generators $\delta_0<\cdots<\delta_t$ satisfy $s_i<s_j$ and $r_i<r_j$ whenever $i<j$. By Proposition~\ref{prop:lambda-semirraiz-normalizada}, we have $\lambda<\beta_2$ and use the normalized parametrization \eqref{eq:normalizacao-lambda}.

Set $r_{-1}=0$. In this case, the partition~\eqref{eq:particao-geradores} admits the description
\begin{equation}\label{eq:blocos-verticais}
 A_i=\left\{\delta_i+a m+b n_1:
 0\leq a\leq n_1-1-s_i,\quad
 0\leq b\leq r_i-r_{i-1}-1\right\}.
\end{equation}
Indeed, the gaps in $\delta_i+\Gamma_1$ are represented by $s m-r n_1$, with $s_i\leq s<n_1$ and $1\leq r\leq r_i$. Since $s_h<s_i$ and $r_h\leq r_{i-1}$ for $h<i$, the union of the preceding blocks removes precisely the elements with $r\leq r_{i-1}$. In particular,
\[
 \min A_i=\delta_i,\qquad
 M_i:=\max A_i=(n_1-1)m-(r_{i-1}+1)n_1.
\]
The sequence $(M_i)$ is strictly decreasing.

On the scale of $f$, define
\begin{equation}\label{eq:blocos-B-verticais}
 \delta_i'=\beta_2+(r_{i-1}+1)v_0,\qquad
 B_i=\delta_i'+\varepsilon_1(A_i-\delta_i).
\end{equation}
Since $v_2=(n_1-1)v_1+\beta_2$, the symmetry of each rectangle gives
\[
 B_i=v_2-\varepsilon_1A_i,\qquad \min B_i=\delta_i'.
\]
The blocks $\varepsilon_1A_i$ and $B_i$ consist of gaps of $\Gamma_f$ smaller than $v_2$. Indeed, $\varepsilon_1\max A_i\leq\varepsilon_1(\mu_1-1)=v_2-\beta_2-v_0<v_2$, and every element of $\Gamma_f$ below $v_2$ belongs to $\varepsilon_1\Gamma_1$. The elements of $B_i$ are positive, smaller than $v_2$, and congruent to $\beta_2\not\equiv0\pmod{\varepsilon_1}$. The blocks in each family are pairwise disjoint; the two families are disjoint by congruence modulo $\varepsilon_1$.

\lem\label{lem:coeficientes-blocos-verticais}
For $\delta\in A_i$, we have
\begin{equation}\label{eq:coeficientes-blocos-verticais}
 \nu_1(B_\delta)=(r_{i-1}+1)n_1+\delta-\delta_i.
\end{equation}
Moreover, the sequence $\bigl(\varepsilon_1\delta_i-(r_{i-1}+1)v_0\bigr)_{i=0}^t$ is strictly increasing, and its first term is $\lambda$.
\flem

\dem The formula for the coefficients follows from Lemma~\ref{lem:coeficientes-particao-retangular}, applied to \eqref{eq:blocos-verticais}, since $K-M_i=(r_{i-1}+1)n_1$.

For $i<t$, the value
\[
 \eta=\delta_i+(s_{i+1}-s_i)m
     =\delta_{i+1}+(r_{i+1}-r_i)n_1
\]
belongs to $A_i$. By the formula just proved and \eqref{eq:monotonia-coeficientes},
\[
 (r_{i-1}+1)n_1+(s_{i+1}-s_i)m
 =\nu_1(B_\eta)
 \geq\nu_1(B_{\delta_{i+1}})+(r_{i+1}-r_i)n_1
 =(r_{i+1}+1)n_1.
\]
The inequality is strict: equality would imply $n_1\mid(s_{i+1}-s_i)$, although $0<s_{i+1}-s_i<n_1$. Multiplying by $\varepsilon_1$ and rearranging, we obtain
\[
 \varepsilon_1\delta_{i+1}-(r_i+1)v_0
 >
 \varepsilon_1\delta_i-(r_{i-1}+1)v_0.
\]
Finally, $\varepsilon_1\delta_0-v_0=\lambda$. \f

\subsection{Examples}

The following example illustrates the decomposition of $\Lambda_1\setminus\Gamma_1$ into blocks $A_i$, $i=0,\ldots,5$.

\ex\label{ex2} Consider $f_1=(t^8,t^{25}+t^{26})$, whose special gaps are highlighted in green in the following figure.

\begin{figure}[H]\centering
	\setlength{\tabcolsep}{2pt}\footnotesize\noindent
	\begin{tabular}{r|ccccccccccccccccccccc}
		175 & \cellcolor{green!25}167 & \cellcolor{green!25}159 & \cellcolor{green!25}151 & \cellcolor{green!25}143 & \cellcolor{green!25}135 & \cellcolor{green!25}127 & \cellcolor{green!25}119 & \cellcolor{green!25}111 & \cellcolor{green!25}103 & 95 & 87 & 79 & 71 & 63 & 55 & 47 & 39 & 31 & 23 & 15 & 7 \\
		150 & \cellcolor{green!25}142 & \cellcolor{green!25}134 & \cellcolor{green!25}126 & \cellcolor{green!25}118 & \cellcolor{green!25}110 & \cellcolor{green!25}102 & \cellcolor{green!25}94 & \cellcolor{green!25}86 & 78 & 70 & 62 & 54 & 46 & 38 & 30 & 22 & 14 & 6 \\
		125 & \cellcolor{green!25}117 & \cellcolor{green!25}109 & \cellcolor{green!25}101 & \cellcolor{green!25}93 & \cellcolor{green!25}85 & \cellcolor{green!25}77 & 69 & 61 & 53 & 45 & 37 & 29 & 21 & 13 & 5 \\
		100 & \cellcolor{green!25}92 & \cellcolor{green!25}84 & \cellcolor{green!25}76 & \cellcolor{green!25}68 & \cellcolor{green!25}60 & 52 & 44 & 36 & 28 & 20 & 12 & 4 \\
		75 & \cellcolor{green!25}67 & \cellcolor{green!25}59 & \cellcolor{green!25}51 & 43 & 35 & 27 & 19 & 11 & 3 \\
		50 & \cellcolor{green!25}42 & \cellcolor{green!25}34 & 26 & 18 & 10 & 2 \\
		25 & 17 & 9 & 1 \\
		0 &  \\
	\end{tabular}
\end{figure}

In this case, $\delta_0=34$, $\delta_1=51$, $\delta_2=60$, $\delta_3=77$, $\delta_4=86$, and $\delta_5=103$. The sets $A_i$, $i=0,\ldots,5$, are shown in different colors in the following figure.

\begin{figure}[H]\centering
	\setlength{\tabcolsep}{2pt}\footnotesize\noindent
	\begin{tabular}{r|ccccccccccccccccccccc}
		175 & \cellcolor{brown!60}167 & \cellcolor{brown!60}159 & \cellcolor{blue!25}151 & \cellcolor{red!25}143 & \cellcolor{red!25}135 & \cellcolor{yellow!15}127 & \cellcolor{gray!25}119 & \cellcolor{gray!25}111 & \cellcolor{orange!55}103 & 95 & 87 & 79 & 71 & 63 & 55 & 47 & 39 & 31 & 23 & 15 & 7 \\
		150 & \cellcolor{brown!60}142 & \cellcolor{brown!60}134 & \cellcolor{blue!25}126 & \cellcolor{red!25}118 & \cellcolor{red!25}110 & \cellcolor{yellow!15}102 & \cellcolor{gray!25}94 & \cellcolor{gray!55}86 & 78 & 70 & 62 & 54 & 46 & 38 & 30 & 22 & 14 & 6 \\
		125 & \cellcolor{brown!60}117 & \cellcolor{brown!60}109 & \cellcolor{blue!25}101 & \cellcolor{red!25}93 & \cellcolor{red!25}85 & \cellcolor{yellow!55}77 & 69 & 61 & 53 & 45 & 37 & 29 & 21 & 13 & 5 \\
		100 & \cellcolor{brown!60}92 & \cellcolor{brown!60}84 & \cellcolor{blue!25}76 & \cellcolor{red!25}68 & \cellcolor{red!55}60 & 52 & 44 & 36 & 28 & 20 & 12 & 4 \\
		75 & \cellcolor{brown!60}67 & \cellcolor{brown!60}59 & \cellcolor{blue!55}51 & 43 & 35 & 27 & 19 & 11 & 3 \\
		50 & \cellcolor{brown!60}42 & \cellcolor{brown}34 & 26 & 18 & 10 & 2 \\
		25 & 17 & 9 & 1 \\
		0 &  \\
	\end{tabular}
\end{figure}

\fex

The next example illustrates how the gaps of $f_1$ determine the special gaps of $f$.

\ex Returning to Example~\ref{ex2}, consider $f=(t^{16},t^{50}+t^{52}+t^{53})$. In the following figure, the special gaps of the forms $v_2+\varepsilon_1\delta$ and $v_2-\varepsilon_1\delta$ are highlighted in green. The other colored gaps are obtained by multiplying the gaps of $f_1$ by $\varepsilon_1=2$ and are not necessarily special gaps of $f$.

\begin{figure}[H]\centering
	\setlength{\tabcolsep}{1pt}\tiny\noindent
	
	\begin{tabular}{r|ccccccccccccccccccccccccccccccccccccccccccccccc}
		753 & \cellcolor{green!25}737 & \cellcolor{green!25}721 & \cellcolor{green!25}705 & \cellcolor{green!25}689 & \cellcolor{green!25}673 & \cellcolor{green!25}657 & \cellcolor{green!25}641 & \cellcolor{green!25}625 & \cellcolor{green!25}609 & \cellcolor{green!25}593 & \cellcolor{green!25}577 & \cellcolor{green!25}561 & \cellcolor{green!25}545 & \cellcolor{green!25}529 & \cellcolor{green!25}513 & \cellcolor{green!25}497 & \cellcolor{green!25}481 & \cellcolor{green!25}465 & \cellcolor{green!25}449 & \cellcolor{green!25}433 & \cellcolor{green!25}417 & \cellcolor{green!25}401 & \cellcolor{green!25}385 & \cellcolor{green!25}369 & 353 & 337 & 321 & 305 & 289 & 273 & 257 & 241 & 225\\
		703 & \cellcolor{green!25}687 & \cellcolor{green!25}671 & \cellcolor{green!25}655 & \cellcolor{green!25}639 & \cellcolor{green!25}623 & \cellcolor{green!25}607 & \cellcolor{green!25}591 & \cellcolor{green!25}575 & \cellcolor{green!25}559 & \cellcolor{green!25}543 & \cellcolor{green!25}527 & \cellcolor{green!25}511 & \cellcolor{green!25}495 & \cellcolor{green!25}479 & \cellcolor{green!25}463 & \cellcolor{green!25}447 & \cellcolor{green!25}431 & \cellcolor{green!25}415 & \cellcolor{green!25}399 & \cellcolor{green!25}383 & \cellcolor{green!25}367 & \cellcolor{green!25}351 & 335 & 319 & 303 & 287 & 271 & 255 & 239 & 223 & 207 & 191 & 175\\
		653 & \cellcolor{green!25}637 & \cellcolor{green!25}621 & \cellcolor{green!25}605 & \cellcolor{green!25}589 & \cellcolor{green!25}573 & \cellcolor{green!25}557 & \cellcolor{green!25}541 & \cellcolor{green!25}525 & \cellcolor{green!25}509 & \cellcolor{green!25}493 & \cellcolor{green!25}477 & \cellcolor{green!25}461 & \cellcolor{green!25}445 & \cellcolor{green!25}429 & \cellcolor{green!25}413 & \cellcolor{green!25}397 & \cellcolor{green!25}381 & \cellcolor{green!25}365 & \cellcolor{green!25}349 & \cellcolor{green!25}333 & \cellcolor{green!25}317 & 301 & 285 & 269 & 253 & 237 & 221 & 205 & 189 & 173 & 157 & 141 & 125\\
		603 & \cellcolor{green!25}587 & \cellcolor{green!25}571 & \cellcolor{green!25}555 & \cellcolor{green!25}539 & \cellcolor{green!25}523 & \cellcolor{green!25}507 & \cellcolor{green!25}491 & \cellcolor{green!25}475 & \cellcolor{green!25}459 & \cellcolor{green!25}443 & \cellcolor{green!25}427 & \cellcolor{green!25}411 & \cellcolor{green!25}395 & \cellcolor{green!25}379 & \cellcolor{green!25}363 & \cellcolor{green!25}347 & \cellcolor{green!25}331 & \cellcolor{green!25}315 & \cellcolor{green!25}299 & 283 & 267 & 251 & 235 & 219 & 203 & 187 & 171 & 155 & 139 & 123 & 107 & 91 & 75\\
		553 & \cellcolor{green!25}537 & \cellcolor{green!25}521 & \cellcolor{green!25}505 & \cellcolor{green!25}489 & \cellcolor{green!25}473 & \cellcolor{green!25}457 & \cellcolor{green!25}441 & \cellcolor{green!25}425 & \cellcolor{green!25}409 & \cellcolor{green!25}393 & \cellcolor{green!25}377 & \cellcolor{green!25}361 & \cellcolor{green!25}345 & \cellcolor{green!25}329 & \cellcolor{green!25}313 & \cellcolor{green!25}297 & \cellcolor{green!25}281 & \cellcolor{green!25}265 & 249 & 233 & 217 & 201 & 185 & 169 & 153 & 137 & 121 & 105 & 89 & 73 & 57 & 41 & 25\\
		503 & \cellcolor{green!25}487 & \cellcolor{green!25}471 & \cellcolor{green!25}455 & \cellcolor{green!25}439 & \cellcolor{green!25}423 & \cellcolor{green!25}407 & \cellcolor{green!25}391 & \cellcolor{green!25}375 & \cellcolor{green!25}359 & \cellcolor{green!25}343 & \cellcolor{green!25}327 & \cellcolor{green!25}311 & \cellcolor{green!25}295 & \cellcolor{green!25}279 & \cellcolor{green!25}263 & \cellcolor{green!25}247 & 231 & 215 & 199 & 183 & 167 & 151 & 135 & 119 & 103 & 87 & 71 & 55 & 39 & 23 & 7 \\
		453 & \cellcolor{green!25}437 & \cellcolor{green!25}421 & \cellcolor{green!25}405 & \cellcolor{green!25}389 & \cellcolor{green!25}373 & \cellcolor{green!25}357 & \cellcolor{green!25}341 & \cellcolor{green!25}325 & \cellcolor{green!25}309 & \cellcolor{green!25}293 & \cellcolor{green!25}277 & \cellcolor{green!25}261 & \cellcolor{green!25}245 & \cellcolor{green!25}229 & \cellcolor{green!25}213 & 197 & 181 & 165 & 149 & 133 & 117 & 101 & 85  & 69 & 53 & 37 & 21 & 5 \\
		403 & 387 & 371 & 355 & 339 & 323 & 307 & 291 & 275 & 259 & 243 & 227 & 211 & 195 & 179 & 163 & 147 & 131 & 115 & 99 & 83 & 67 & 51 & 35 & 19 & 3 \\
		350 & \cellcolor{brown!60}334 & \cellcolor{brown!60}318 & \cellcolor{blue!25}302 & \cellcolor{red!25}286 & \cellcolor{red!25}270 & \cellcolor{yellow!15}254 & \cellcolor{gray!25}238 & \cellcolor{gray!25}222 & \cellcolor{orange!55}206 & 190 & 174 & 158 & 142 & 126 & 110 & 94 & 78 & 62 & 46 & 30 & 14 \\
		300 & \cellcolor{brown!60}284 & \cellcolor{brown!60}268 & \cellcolor{blue!25}252 & \cellcolor{red!25}236 & \cellcolor{red!25}220 & \cellcolor{yellow!15}204 & \cellcolor{gray!25}188 & \cellcolor{gray!55}172 & 156 & 140 & 124 & 108 & 92 & 76 & 60 & 44 & 28 & 12 \\
		250 & \cellcolor{brown!60}234 & \cellcolor{brown!60}218 & \cellcolor{blue!25}202 & \cellcolor{red!25}186 & \cellcolor{red!25}170 & \cellcolor{yellow!55}154 & 138 & 122 & 106 & 90 & 74 & 58 & 42 & 26 & 10 \\
		200 & \cellcolor{brown!60}184 & \cellcolor{brown!60}168 & \cellcolor{blue!25}152 & \cellcolor{red!25}136 & \cellcolor{red!55}120 & 104 & 88 & 72 & 56 & 40 & 24 & 8 \\
		150 & \cellcolor{brown!60}134 & \cellcolor{brown!60}118 & \cellcolor{blue!55}102 & 86 & 70 & 54 & 38 & 22 & 6 \\
		100 & \cellcolor{brown!60}84 & \cellcolor{brown}68 & 52 & 36 & 20 & 4 \\
		50 & 34 & 18 & 2 \\
		0 &  \\
	\end{tabular}
\end{figure}

We have $\delta_0=34$, $2\delta_0\in\Lambda_f\setminus\Gamma_f$, and $\beta_2=53$. For $w_{\delta_0}=v_1y\,dx-v_0x\,dy$, we have $\nu_{f_1}(w_{\delta_0})=\delta_0$. Denoting the coefficient of $dy$ by $B_{\delta_0}=v_0x$, we obtain $\nu_1(B_{\delta_0})=v_0/\varepsilon_1$ and hence
$$\varepsilon_1\delta_0=68 < 69 = \beta_2 + \varepsilon_1\nu_{f_1}(B_{\delta_0})= \beta_2 + \nu_{f}(B_{\delta_0}) = \beta_2+v_0=\beta_2+(r_{-1}+1)v_0,$$
where $r_{-1}=0$.

For $i=1,\ldots,5$, comparing $\varepsilon_1\delta_i$ with $\beta_2+(r_{i-1}+1)v_0$ gives
$$
\begin{aligned}
\varepsilon_1\delta_1&=102 > 101 = \beta_2 + 3v_0 = \beta_2+(r_0+1)v_0\\
\varepsilon_1\delta_2&=120 > 117 = \beta_2 + 4v_0 = \beta_2+(r_1+1)v_0\\
\varepsilon_1\delta_3&=154 > 149 = \beta_2 + 6v_0 = \beta_2+(r_2+1)v_0\\
\varepsilon_1\delta_4&=172 > 165 = \beta_2 + 7v_0 = \beta_2+(r_3+1)v_0\\
\varepsilon_1\delta_5&=206 > 197 = \beta_2 + 9v_0 = \beta_2+(r_4+1)v_0.
\end{aligned}
$$

By identity~\eqref{eq:coeficientes-blocos-verticais}, $\nu_f(B_{\delta_i})=(r_{i-1}+1)v_0$. We thus obtain the special gaps of $f$ highlighted in the following table:

\begin{figure}[H]\centering
	\setlength{\tabcolsep}{1pt}\tiny\noindent
	
	\begin{tabular}{r|ccccccccccccccccccccccccccccccccccccccccccccccc}
753 & \cellcolor{green!25}737 & \cellcolor{green!25}721 & \cellcolor{green!25}705 & \cellcolor{green!25}689 & \cellcolor{green!25}673 & \cellcolor{green!25}657 & \cellcolor{green!25}641 & \cellcolor{green!25}625 & \cellcolor{green!25}609 & \cellcolor{green!25}593 & \cellcolor{green!25}577 & \cellcolor{green!25}561 & \cellcolor{green!25}545 & \cellcolor{green!25}529 & \cellcolor{green!25}513 & \cellcolor{green!25}497 & \cellcolor{green!25}481 & \cellcolor{green!25}465 & \cellcolor{green!25}449 & \cellcolor{green!25}433 & \cellcolor{green!25}417 & \cellcolor{green!25}401 & \cellcolor{green!25}385 & \cellcolor{green!25}369 & 353 & 337 & 321 & 305 & 289 & 273 & 257 & 241 & 225\\
703 & \cellcolor{green!25}687 & \cellcolor{green!25}671 & \cellcolor{green!25}655 & \cellcolor{green!25}639 & \cellcolor{green!25}623 & \cellcolor{green!25}607 & \cellcolor{green!25}591 & \cellcolor{green!25}575 & \cellcolor{green!25}559 & \cellcolor{green!25}543 & \cellcolor{green!25}527 & \cellcolor{green!25}511 & \cellcolor{green!25}495 & \cellcolor{green!25}479 & \cellcolor{green!25}463 & \cellcolor{green!25}447 & \cellcolor{green!25}431 & \cellcolor{green!25}415 & \cellcolor{green!25}399 & \cellcolor{green!25}383 & \cellcolor{green!25}367 & \cellcolor{green!25}351 & 335 & 319 & 303 & 287 & 271 & 255 & 239 & 223 & 207 & 191 & 175\\
653 & \cellcolor{green!25}637 & \cellcolor{green!25}621 & \cellcolor{green!25}605 & \cellcolor{green!25}589 & \cellcolor{green!25}573 & \cellcolor{green!25}557 & \cellcolor{green!25}541 & \cellcolor{green!25}525 & \cellcolor{green!25}509 & \cellcolor{green!25}493 & \cellcolor{green!25}477 & \cellcolor{green!25}461 & \cellcolor{green!25}445 & \cellcolor{green!25}429 & \cellcolor{green!25}413 & \cellcolor{green!25}397 & \cellcolor{green!25}381 & \cellcolor{green!25}365 & \cellcolor{green!25}349 & \cellcolor{green!25}333 & \cellcolor{green!25}317 & \cellcolor{blue!25}301 & 285 & 269 & 253 & 237 & 221 & 205 & 189 & 173 & 157 & 141 & 125\\
603 & \cellcolor{green!25}587 & \cellcolor{green!25}571 & \cellcolor{green!25}555 & \cellcolor{green!25}539 & \cellcolor{green!25}523 & \cellcolor{green!25}507 & \cellcolor{green!25}491 & \cellcolor{green!25}475 & \cellcolor{green!25}459 & \cellcolor{green!25}443 & \cellcolor{green!25}427 & \cellcolor{green!25}411 & \cellcolor{green!25}395 & \cellcolor{green!25}379 & \cellcolor{green!25}363 & \cellcolor{green!25}347 & \cellcolor{green!25}331 & \cellcolor{green!25}315 & \cellcolor{green!25}299 & \cellcolor{red!25}283 & \cellcolor{red!25}267 & \cellcolor{blue!25}251 & 235 & 219 & 203 & 187 & 171 & 155 & 139 & 123 & 107 & 91 & 75\\
553 & \cellcolor{green!25}537 & \cellcolor{green!25}521 & \cellcolor{green!25}505 & \cellcolor{green!25}489 & \cellcolor{green!25}473 & \cellcolor{green!25}457 & \cellcolor{green!25}441 & \cellcolor{green!25}425 & \cellcolor{green!25}409 & \cellcolor{green!25}393 & \cellcolor{green!25}377 & \cellcolor{green!25}361 & \cellcolor{green!25}345 & \cellcolor{green!25}329 & \cellcolor{green!25}313 & \cellcolor{green!25}297 & \cellcolor{green!25}281 & \cellcolor{green!25}265 & \cellcolor{yellow!15}249 & \cellcolor{red!25}233 & \cellcolor{red!25}217 & \cellcolor{blue!25}201 & 185 & 169 & 153 & 137 & 121 & 105 & 89 & 73 & 57 & 41 & 25\\
503 & \cellcolor{green!25}487 & \cellcolor{green!25}471 & \cellcolor{green!25}455 & \cellcolor{green!25}439 & \cellcolor{green!25}423 & \cellcolor{green!25}407 & \cellcolor{green!25}391 & \cellcolor{green!25}375 & \cellcolor{green!25}359 & \cellcolor{green!25}343 & \cellcolor{green!25}327 & \cellcolor{green!25}311 & \cellcolor{green!25}295 & \cellcolor{green!25}279 & \cellcolor{green!25}263 & \cellcolor{green!25}247 & \cellcolor{gray!25}231 & \cellcolor{gray!25}215 & \cellcolor{yellow!15}199 & \cellcolor{red!25}183 & \cellcolor{red!25}167 & \cellcolor{blue!25}151 & 135 & 119 & 103 & 87 & 71 & 55 & 39 & 23 & 7 \\
453 & \cellcolor{green!25}437 & \cellcolor{green!25}421 & \cellcolor{green!25}405 & \cellcolor{green!25}389 & \cellcolor{green!25}373 & \cellcolor{green!25}357 & \cellcolor{green!25}341 & \cellcolor{green!25}325 & \cellcolor{green!25}309 & \cellcolor{green!25}293 & \cellcolor{green!25}277 & \cellcolor{green!25}261 & \cellcolor{green!25}245 & \cellcolor{green!25}229 & \cellcolor{green!25}213 & \cellcolor{orange!55}197 & \cellcolor{gray!25}181 & \cellcolor{gray!55}165 & \cellcolor{yellow!55}149 & \cellcolor{red!25}133 & \cellcolor{red!55}117 & \cellcolor{blue!55}101 & 85  & 69 & 53 & 37 & 21 & 5 \\
403 & 387 & 371 & 355 & 339 & 323 & 307 & 291 & 275 & 259 & 243 & 227 & 211 & 195 & 179 & 163 & 147 & 131 & 115 & 99 & 83 & 67 & 51 & 35 & 19 & 3 \\
350 & \cellcolor{brown!60}334 & \cellcolor{brown!60}318 & 302 & 286 & 270 & 254 & 238 & 222 & 206 & 190 & 174 & 158 & 142 & 126 & 110 & 94 & 78 & 62 & 46 & 30 & 14 \\
300 & \cellcolor{brown!60}284 & \cellcolor{brown!60}268 & 252 & 236 & 220 & 204 & 188 & 172 & 156 & 140 & 124 & 108 & 92 & 76 & 60 & 44 & 28 & 12 \\
250 & \cellcolor{brown!60}234 & \cellcolor{brown!60}218 & 202 & 186 & 170 & 154 & 138 & 122 & 106 & 90 & 74 & 58 & 42 & 26 & 10 \\
200 & \cellcolor{brown!60}184 & \cellcolor{brown!60}168 & 152 & 136 & 120 & 104 & 88 & 72 & 56 & 40 & 24 & 8 \\
150 & \cellcolor{brown!60}134 & \cellcolor{brown!60}118 & 102 & 86 & 70 & 54 & 38 & 22 & 6 \\
100 & \cellcolor{brown!60}84 & \cellcolor{brown}68 & 52 & 36 & 20 & 4 \\
50 & 34 & 18 & 2 \\
0 &  \\
	\end{tabular}
\end{figure}

The figure shows precisely these special gaps of $f$. \fex

\subsection{Image of the transition map and consequences}

\lem\label{lem:transicao-blocos}
If $\varepsilon_1\delta_j>\delta_j'$, then
\[
 \bigcup_{k=j}^t B_k\subset\Lambda_f\setminus\Gamma_f.
\]
\flem
\dem By Lemma~\ref{lem:coeficientes-blocos-verticais} and the definition of $\rho_1$,
\[
 \nu_f(w_{\delta_j})=\rho_1(\delta_j)
 =\delta_j'=\min B_j.
\]
For $k\geq j$, we have $\delta_k'-\delta_j'=(r_{k-1}-r_{j-1})v_0$, with $r_{k-1}-r_{j-1}\geq0$. Therefore,
\[
 \nu_f\bigl(x^{r_{k-1}-r_{j-1}+b}y^a w_{\delta_j}\bigr)
 =\delta_k'+b v_0+a v_1.
\]
As $a,b$ range over the bounds in \eqref{eq:blocos-verticais} corresponding to the index $k$, these values range over $B_k$. Since they are all gaps of $\Gamma_f$, the inclusion is proved. \f

\teo\label{teo:decomposicao-rho-blocos}
Under the assumption of Case~1, for each $i$,
\[
 \rho_1(A_i)=
 \begin{cases}
  \varepsilon_1A_i,&\varepsilon_1\delta_i<\delta_i',\\
  B_i,&\varepsilon_1\delta_i>\delta_i'.
 \end{cases}
\]
The selected block consists of special gaps of $f$, and the other is disjoint from $\Lambda_f\setminus\Gamma_f$. If $j\in\{1,\ldots,t\}$ is the first index for which $\varepsilon_1\delta_j>\delta_j'$, then
\[
 \rho_1(\Lambda_1\setminus\Gamma_1)
 =\left(\dot\bigcup_{i=0}^{j-1}\varepsilon_1A_i\right)
   \mathbin{\dot\cup}
  \left(\dot\bigcup_{i=j}^t B_i\right).
\]
If no such index exists, then $\rho_1(\Lambda_1\setminus\Gamma_1)=\varepsilon_1(\Lambda_1\setminus\Gamma_1)$. In particular, $\varepsilon_1A_0$ is always the selected block at index zero.
\fteo

\dem By Lemma~\ref{lem:coeficientes-blocos-verticais}, for $\delta\in A_i$,
\[
 \rho_1(\delta)
 =\min\{\varepsilon_1\delta_i,\delta_i'\}
       +\varepsilon_1(\delta-\delta_i).
\]
The comparison between the two alternatives is therefore constant on each block, which proves the formula for $\rho_1(A_i)$. Equality of the minima is impossible, since $\varepsilon_1\delta_i\equiv0\pmod{\varepsilon_1}$ and $\delta_i'\equiv\beta_2\not\equiv0\pmod{\varepsilon_1}$.

To justify the exclusion of the unselected block, apply Proposition~\ref{2.2.15}. The family
\[
 \{v_2-\varepsilon_1\alpha:
       \alpha\in\mathbb N^*\setminus\Lambda_1\}
\]
is disjoint from all $B_i=v_2-\varepsilon_1A_i$, since $A_i\subset\Lambda_1$, and from all $\varepsilon_1A_i$, by congruence modulo $\varepsilon_1$. Moreover, the images of the other blocks belong to disjoint pairs $\varepsilon_1A_i\cup B_i$. Consequently,
\[
 (\varepsilon_1A_i\cup B_i)\cap(\Lambda_f\setminus\Gamma_f)
 =\rho_1(A_i).
\]
Hence the unselected block is disjoint from $\Lambda_f\setminus\Gamma_f$.

Finally, Lemma~\ref{lem:coeficientes-blocos-verticais} shows that the values $\varepsilon_1\delta_i-(r_{i-1}+1)v_0$ are strictly increasing. Hence the inequality $\varepsilon_1\delta_i>\delta_i'$ persists after its first occurrence. Since the first of these values is $\lambda<\beta_2$, index zero selects $\varepsilon_1A_0$. The formulas for the full image follow from the partition of $\Lambda_1\setminus\Gamma_1$. \f

The following figure represents the candidate blocks $\varepsilon_1A_i$ and $B_i=v_2-\varepsilon_1A_i$. The green region corresponds to the family $\{v_2-\varepsilon_1\alpha:\alpha\in\mathbb N^*\setminus\Lambda_1\}$. In each pair of blocks, only the one selected by Theorem~\ref{teo:decomposicao-rho-blocos} consists of special gaps. The minima of the first pair are $\min(\varepsilon_1A_0)=\lambda+v_0$ and $\min B_0=\beta_2+v_0$.

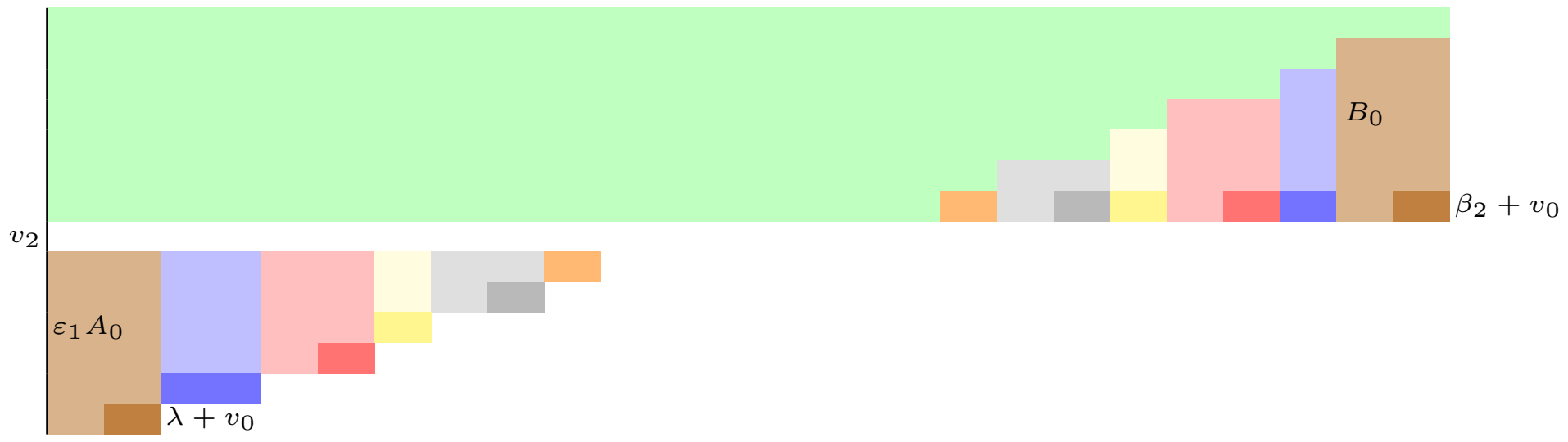
\begin{figure}[H]\centering
	\setlength{\tabcolsep}{1pt}\tiny\noindent
	
	\begin{tabular}{r|ccccccccccccccccccccccccccccccccccccccccccccccc}
		& \cellcolor{green!25}\textcolor{green!25}{737} & \cellcolor{green!25}\textcolor{green!25}{721} & \cellcolor{green!25}\textcolor{green!25}{705} & \cellcolor{green!25}\textcolor{green!25}{689} & \cellcolor{green!25}\textcolor{green!25}{673} & \cellcolor{green!25}\textcolor{green!25}{657} & \cellcolor{green!25}\textcolor{green!25}{641} & \cellcolor{green!25}\textcolor{green!25}{625} & \cellcolor{green!25}\textcolor{green!25}{609} & \cellcolor{green!25}\textcolor{green!25}{593} & \cellcolor{green!25}\textcolor{green!25}{577} & \cellcolor{green!25}\textcolor{green!25}{561} & \cellcolor{green!25}\textcolor{green!25}{545} & \cellcolor{green!25}\textcolor{green!25}{529} & \cellcolor{green!25}\textcolor{green!25}{513} & \cellcolor{green!25}\textcolor{green!25}{497} & \cellcolor{green!25}\textcolor{green!25}{481} & \cellcolor{green!25}\textcolor{green!25}{465} & \cellcolor{green!25}\textcolor{green!25}{449} & \cellcolor{green!25}\textcolor{green!25}{433} & \cellcolor{green!25}\textcolor{green!25}{417} & \cellcolor{green!25}\textcolor{green!25}{401} & \cellcolor{green!25}\textcolor{green!25}{385} & \cellcolor{green!25}\textcolor{green!25}{369} \\
		& \cellcolor{green!25}\textcolor{green!25}{687} & \cellcolor{green!25}\textcolor{green!25}{671} & \cellcolor{green!25}\textcolor{green!25}{655} & \cellcolor{green!25}\textcolor{green!25}{639} & \cellcolor{green!25}\textcolor{green!25}{623} & \cellcolor{green!25}\textcolor{green!25}{607} & \cellcolor{green!25}\textcolor{green!25}{591} & \cellcolor{green!25}\textcolor{green!25}{575} & \cellcolor{green!25}\textcolor{green!25}{559} & \cellcolor{green!25}\textcolor{green!25}{543} & \cellcolor{green!25}\textcolor{green!25}{527} & \cellcolor{green!25}\textcolor{green!25}{511} & \cellcolor{green!25}\textcolor{green!25}{495} & \cellcolor{green!25}\textcolor{green!25}{479} & \cellcolor{green!25}\textcolor{green!25}{463} & \cellcolor{green!25}\textcolor{green!25}{447} & \cellcolor{green!25}\textcolor{green!25}{431} & \cellcolor{green!25}\textcolor{green!25}{415} & \cellcolor{green!25}\textcolor{green!25}{399} & \cellcolor{green!25}\textcolor{green!25}{383} & \cellcolor{green!25}\textcolor{green!25}{367} & \cellcolor{green!25}\textcolor{green!25}{351} & \cellcolor{brown!60}\textcolor{brown!60}{335} & \cellcolor{brown!60}\textcolor{brown!60}{319} \\
		& \cellcolor{green!25}\textcolor{green!25}{637} & \cellcolor{green!25}\textcolor{green!25}{621} & \cellcolor{green!25}\textcolor{green!25}{605} & \cellcolor{green!25}\textcolor{green!25}{589} & \cellcolor{green!25}\textcolor{green!25}{573} & \cellcolor{green!25}\textcolor{green!25}{557} & \cellcolor{green!25}\textcolor{green!25}{541} & \cellcolor{green!25}\textcolor{green!25}{525} & \cellcolor{green!25}\textcolor{green!25}{509} & \cellcolor{green!25}\textcolor{green!25}{493} & \cellcolor{green!25}\textcolor{green!25}{477} & \cellcolor{green!25}\textcolor{green!25}{461} & \cellcolor{green!25}\textcolor{green!25}{445} & \cellcolor{green!25}\textcolor{green!25}{429} & \cellcolor{green!25}\textcolor{green!25}{413} & \cellcolor{green!25}\textcolor{green!25}{397} & \cellcolor{green!25}\textcolor{green!25}{381} & \cellcolor{green!25}\textcolor{green!25}{365} & \cellcolor{green!25}\textcolor{green!25}{349} & \cellcolor{green!25}\textcolor{green!25}{333} & \cellcolor{green!25}\textcolor{green!25}{317} & \cellcolor{blue!25}\textcolor{blue!25}{301} & \cellcolor{brown!60}\textcolor{brown!60}{285} & \cellcolor{brown!60}\textcolor{brown!60}{269}  \\
		& \cellcolor{green!25}\textcolor{green!25}{587} & \cellcolor{green!25}\textcolor{green!25}{571} & \cellcolor{green!25}\textcolor{green!25}{555} & \cellcolor{green!25}\textcolor{green!25}{539} & \cellcolor{green!25}\textcolor{green!25}{523} & \cellcolor{green!25}\textcolor{green!25}{507} & \cellcolor{green!25}\textcolor{green!25}{491} & \cellcolor{green!25}\textcolor{green!25}{475} & \cellcolor{green!25}\textcolor{green!25}{459} & \cellcolor{green!25}\textcolor{green!25}{443} & \cellcolor{green!25}\textcolor{green!25}{427} & \cellcolor{green!25}\textcolor{green!25}{411} & \cellcolor{green!25}\textcolor{green!25}{395} & \cellcolor{green!25}\textcolor{green!25}{379} & \cellcolor{green!25}\textcolor{green!25}{363} & \cellcolor{green!25}\textcolor{green!25}{347} & \cellcolor{green!25}\textcolor{green!25}{331} & \cellcolor{green!25}\textcolor{green!25}{315} & \cellcolor{green!25}\textcolor{green!25}{299} & \cellcolor{red!25}\textcolor{red!25}{283} & \cellcolor{red!25}\textcolor{red!25}{267} & \cellcolor{blue!25}\textcolor{blue!25}{251} & \cellcolor{brown!60}$B_0$ & \cellcolor{brown!60}\textcolor{brown!60}{219}  \\
		& \cellcolor{green!25}\textcolor{green!25}{537} & \cellcolor{green!25}\textcolor{green!25}{521} & \cellcolor{green!25}\textcolor{green!25}{505} & \cellcolor{green!25}\textcolor{green!25}{489} & \cellcolor{green!25}\textcolor{green!25}{473} & \cellcolor{green!25}\textcolor{green!25}{457} & \cellcolor{green!25}\textcolor{green!25}{441} & \cellcolor{green!25}\textcolor{green!25}{425} & \cellcolor{green!25}\textcolor{green!25}{409} & \cellcolor{green!25}\textcolor{green!25}{393} & \cellcolor{green!25}\textcolor{green!25}{377} & \cellcolor{green!25}\textcolor{green!25}{361} & \cellcolor{green!25}\textcolor{green!25}{345} & \cellcolor{green!25}\textcolor{green!25}{329} & \cellcolor{green!25}\textcolor{green!25}{313} & \cellcolor{green!25}\textcolor{green!25}{297} & \cellcolor{green!25}\textcolor{green!25}{281} & \cellcolor{green!25}\textcolor{green!25}{265} & \cellcolor{yellow!15}\textcolor{yellow!15}{249} & \cellcolor{red!25}\textcolor{red!25}{233} & \cellcolor{red!25}\textcolor{red!25}{217} & \cellcolor{blue!25}\textcolor{blue!25}{201} & \cellcolor{brown!60}\textcolor{brown!60}{185} & \cellcolor{brown!60}\textcolor{brown!60}{169}\\
		& \cellcolor{green!25}\textcolor{green!25}{487} & \cellcolor{green!25}\textcolor{green!25}{471} & \cellcolor{green!25}\textcolor{green!25}{455} & \cellcolor{green!25}\textcolor{green!25}{439} & \cellcolor{green!25}\textcolor{green!25}{423} & \cellcolor{green!25}\textcolor{green!25}{407} & \cellcolor{green!25}\textcolor{green!25}{391} & \cellcolor{green!25}\textcolor{green!25}{375} & \cellcolor{green!25}\textcolor{green!25}{359} & \cellcolor{green!25}\textcolor{green!25}{343} & \cellcolor{green!25}\textcolor{green!25}{327} & \cellcolor{green!25}\textcolor{green!25}{311} & \cellcolor{green!25}\textcolor{green!25}{295} & \cellcolor{green!25}\textcolor{green!25}{279} & \cellcolor{green!25}\textcolor{green!25}{263} & \cellcolor{green!25}\textcolor{green!25}{247} & \cellcolor{gray!25}\textcolor{gray!25}{231} & \cellcolor{gray!25}\textcolor{gray!25}{215} & \cellcolor{yellow!15}\textcolor{yellow!15}{199} & \cellcolor{red!25}\textcolor{red!25}{183} & \cellcolor{red!25}\textcolor{red!25}{167} & \cellcolor{blue!25}\textcolor{blue!25}{151} & \cellcolor{brown!60}\textcolor{brown!60}{135} & \cellcolor{brown!60}\textcolor{brown!60}{119} \\
		& \cellcolor{green!25}\textcolor{green!25}{437} & \cellcolor{green!25}\textcolor{green!25}{421} & \cellcolor{green!25}\textcolor{green!25}{405} & \cellcolor{green!25}\textcolor{green!25}{389} & \cellcolor{green!25}\textcolor{green!25}{373} & \cellcolor{green!25}\textcolor{green!25}{357} & \cellcolor{green!25}\textcolor{green!25}{341} & \cellcolor{green!25}\textcolor{green!25}{325} & \cellcolor{green!25}\textcolor{green!25}{309} & \cellcolor{green!25}\textcolor{green!25}{293} & \cellcolor{green!25}\textcolor{green!25}{277} & \cellcolor{green!25}\textcolor{green!25}{261} & \cellcolor{green!25}\textcolor{green!25}{245} & \cellcolor{green!25}\textcolor{green!25}{229} & \cellcolor{green!25}\textcolor{green!25}{213} & \cellcolor{orange!55}\textcolor{orange!55}{197} & \cellcolor{gray!25}\textcolor{gray!25}{181} & \cellcolor{gray!55}\textcolor{gray!55}{165} & \cellcolor{yellow!55}\textcolor{yellow!55}{149} & \cellcolor{red!25}\textcolor{red!25}{133} & \cellcolor{red!55}\textcolor{red!55}{117} & \cellcolor{blue!55}\textcolor{blue!55}{101} & \cellcolor{brown!60}\textcolor{brown!60}{85}  & \cellcolor{brown}\textcolor{brown}{69}&$\beta_2+v_0$  \\
		$v_2$ &  \\
		& \cellcolor{brown!60}\textcolor{brown!60}{334} & \cellcolor{brown!60}\textcolor{brown!60}{318} & \cellcolor{blue!25}\textcolor{blue!25}{302} & \cellcolor{red!25}\textcolor{red!25}{286} & \cellcolor{red!25}\textcolor{red!25}{270} & \cellcolor{yellow!15}\textcolor{yellow!15}{254} & \cellcolor{gray!25}\textcolor{gray!25}{238} & \cellcolor{gray!25}\textcolor{gray!25}{222} & \cellcolor{orange!55}\textcolor{orange!55}{206}  \\
		& \cellcolor{brown!60}\textcolor{brown!60}{284} & \cellcolor{brown!60}\textcolor{brown!60}{268} & \cellcolor{blue!25}\textcolor{blue!25}{252} & \cellcolor{red!25}\textcolor{red!25}{236} & \cellcolor{red!25}\textcolor{red!25}{220} & \cellcolor{yellow!15}\textcolor{yellow!15}{204} & \cellcolor{gray!25}\textcolor{gray!25}{188} & \cellcolor{gray!55}\textcolor{gray!55}{172} \\
		& \multicolumn{2}{l}{\cellcolor{brown!60}$\varepsilon_1 A_0$} & \cellcolor{blue!25}\textcolor{blue!25}{202} & \cellcolor{red!25}\textcolor{red!25}{186} & \cellcolor{red!25}\textcolor{red!25}{170} & \cellcolor{yellow!55}\textcolor{yellow!55}{154}  \\
		& \cellcolor{brown!60}\textcolor{brown!60}{184} & \cellcolor{brown!60}\textcolor{brown!60}{168} & \cellcolor{blue!25}\textcolor{blue!25}{152} & \cellcolor{red!25}\textcolor{red!25}{136} & \cellcolor{red!55}\textcolor{red!55}{120}  \\
		& \cellcolor{brown!60}\textcolor{brown!60}{134} & \cellcolor{brown!60}\textcolor{brown!60}{118} & \cellcolor{blue!55}\textcolor{blue!55}{102}  \\
		& \cellcolor{brown!60}\textcolor{brown!60}{84} & \cellcolor{brown}\textcolor{brown}{68}&$\lambda+v_0$  \\
	\end{tabular}
\end{figure}

\cor\label{cor:criterio-inclusao-escalada}
Under the assumption of Case~1,
\[
 \varepsilon_1(\Lambda_1\setminus\Gamma_1)
 \subset\Lambda_f\setminus\Gamma_f
 \quad\Longleftrightarrow\quad
 \beta_2>\varepsilon_1\delta_t-(r_{t-1}+1)v_0.
\]
\fcor
\dem By Theorem~\ref{teo:decomposicao-rho-blocos}, the inclusion is equivalent to selecting $\varepsilon_1A_i$ for every $i$, that is,
\[
 \beta_2>\varepsilon_1\delta_i-(r_{i-1}+1)v_0
 \qquad(0\leq i\leq t).
\]
Lemma~\ref{lem:coeficientes-blocos-verticais} shows that the largest term on the right-hand side is the one with index $t$. \f

\cor\label{cor:determinacao-torre-binaria}
Let $g\geq2$ and assume $n_k=2$ for $2\leq k\leq g$. Given $\Gamma_f$ and $\Lambda_1\setminus\Gamma_1$, with the minimal generators satisfying the assumption of Case~1, the set $\Lambda_f\setminus\Gamma_f$ is determined.
\fcor
\dem At each stage, the values of the semiroot are considered on its own scale:
\[
 \Gamma_k=\left\langle
  \frac{v_0}{\varepsilon_k},\ldots,
  \frac{v_k}{\varepsilon_k}\right\rangle.
\]
The special gaps of $f_1$ determine the blocks $A_i$. Since $\beta_2=v_2-(n_1-1)v_1$ is determined by $\Gamma_f$, the blocks $B_i$ and the comparisons between the minima are also determined. Theorem~\ref{teo:decomposicao-rho-blocos} thus determines $\rho_1(\Lambda_1\setminus\Gamma_1)$ on the scale of $f$. By Theorem~\ref{t5.2}, applied to $f_2$,
\begin{equation}\label{eq:primeira-etapa-torre-binaria}
\Lambda_2\setminus\Gamma_2=\frac{1}{\varepsilon_2}
       \rho_1(\Lambda_1\setminus\Gamma_1)\mathbin{\dot\cup}
   \left\{\frac{v_2}{\varepsilon_2}-2\alpha:
          \alpha\in\mathbb N^*\setminus\Lambda_1\right\}\mathbin{\dot\cup}
   \left\{\frac{v_2}{\varepsilon_2}+2\eta:
          \eta\in\mathbb N\setminus\Gamma_1\right\}.
\end{equation}
All terms are determined by the data in the statement.

For $2\leq k<g$, Lemma~\ref{lem:propagacao-total-binaria} gives $\rho_k(\delta)=\varepsilon_k\delta$ for every $\delta\in\Lambda_k\setminus\Gamma_k$, including values greater than $v_k/\varepsilon_k$. Since $\varepsilon_k/\varepsilon_{k+1}=2$, Theorem~\ref{t5.2}, applied to $f_{k+1}$, gives
\begin{equation}\label{eq:recorrencia-torre-binaria}
\Lambda_{k+1}\setminus\Gamma_{k+1}
 =2(\Lambda_k\setminus\Gamma_k)\mathbin{\dot\cup}
   \left\{\frac{v_{k+1}}{\varepsilon_{k+1}}-2\alpha:
          \alpha\in\mathbb N^*\setminus\Lambda_k\right\}\mathbin{\dot\cup}
   \left\{\frac{v_{k+1}}{\varepsilon_{k+1}}+2\eta:
          \eta\in\mathbb N\setminus\Gamma_k\right\}.
\end{equation}
Once $\Lambda_k\setminus\Gamma_k$ is known, we recover $\Lambda_k=(\Gamma_k\setminus\{0\})\cup(\Lambda_k\setminus\Gamma_k)$. By induction, the recurrence determines $\Lambda_g\setminus\Gamma_g=\Lambda_f\setminus\Gamma_f$. The assumption of Case~1 is used only at the first stage; \eqref{eq:recorrencia-torre-binaria} is independent of the decomposition of $\Lambda_1\setminus\Gamma_1$. \f

\subsection{A limitation for larger ramification indices}

For $g\geq3$ and $n_{g-1}\geq3$, the inclusion $\varepsilon_{g-1}\delta_{g-1}\in\Lambda_f\setminus\Gamma_f$ may fail for $\delta_{g-1}\in\Lambda_{g-1}\setminus\Gamma_{g-1}$. The following example exhibits this situation when $\delta_{g-1}\geq v_{g-1}/\varepsilon_{g-1}$.

\ex Consider $f_2=(t^{12},t^{27}+t^{30}+t^{31})$.

\begin{figure}[H]\centering
	\setlength{\tabcolsep}{2pt}\footnotesize\noindent
	\begin{tabular}{r|ccccccccccccccccccccccccc}
		305 & \cellcolor{green!25}293 & \cellcolor{green!25}281 & \cellcolor{green!25}269 & \cellcolor{green!25}257 & \cellcolor{green!25}245 & \cellcolor{green!25}233 & \cellcolor{green!25}221 & \cellcolor{green!25}209 & \cellcolor{green!25}197 & \cellcolor{green!25}185 & \cellcolor{green!25}173 & 161 & 149 & 137 & 125 & 113 & 101 & 89 & 77 & 65 & 53 & 41 & 29 & 17 & 5 \\
		278 & \cellcolor{green!25}266 & \cellcolor{green!25}254 & \cellcolor{green!25}242 & \cellcolor{green!25}230 & \cellcolor{green!25}218 & \cellcolor{green!25}206 & \cellcolor{green!25}194 & \cellcolor{green!25}182 & \cellcolor{green!25}170 & \cellcolor{green!25}158 & 146 & 134 & 122 & 110 & 98 & 86 & 74 & 62 & 50 & 38 & 26 & 14 & 2 \\
		251 & \cellcolor{green!25}239 & \cellcolor{green!25}227 & \cellcolor{green!25}215 & \cellcolor{green!25}203 & \cellcolor{green!25}191 & \cellcolor{green!25}179 & \cellcolor{green!25}167 & \cellcolor{green!25}155 & \cellcolor{green!25}143 & 131 & 119 & 107 & 95 & 83 & 71 & 59 & 47 & 35 & 23 & 11 \\
		224 & \cellcolor{green!25}212 & \cellcolor{green!25}200 & \cellcolor{green!25}188 & \cellcolor{green!25}176 & \cellcolor{green!25}164 & \cellcolor{green!25}152 & \cellcolor{green!25}140 & \cellcolor{green!25}128 & 116 & 104 & 92 & 80 & 68 & 56 & 44 & 32 & 20 & 8 \\
		193 & \cellcolor{green!25}181 & \cellcolor{green!25}169 & \cellcolor{green!25}157 & \cellcolor{green!25}145 & \cellcolor{green!25}133 & \cellcolor{green!25}121 & \cellcolor{green!25}109 & \cellcolor{green!25}97 & 85 & 73 & 61 & 49 & 37 & 25 & 13 & 1 \\
		166 & \cellcolor{green!25}154 & \cellcolor{green!25}142 & \cellcolor{green!25}130 & \cellcolor{green!25}118 & \cellcolor{green!25}106 & \cellcolor{green!25}94 & \cellcolor{green!25}82 & 70 & 58 & 46 & 34 & 22 & 10 \\
		139 & \cellcolor{green!25}127 & \cellcolor{green!25}115 & \cellcolor{green!25}103 & \cellcolor{green!25}91 & \cellcolor{green!25}79 & \cellcolor{green!25}67 & \cellcolor{green!25}55 & 43 & 31 & 19 & 7 \\
		112 & 100 & 88 & 76 & 64 & 52 & 40 & 28 & 16 & 4 \\
		81 & \cellcolor{green!25}69 & 57 & 45 & 33 & 21 & 9 \\
		54 & \cellcolor{green!25}42 & 30 & 18 & 6 \\
		27 & 15 & 3 \\
		0 &  \\
	\end{tabular}
\end{figure}

For $f=(t^{24},t^{54}+t^{60}+t^{62}+t^{63})$, no special gap in the class of $2v_2$ remains special. The same holds, among others, for the gaps $2\cdot155$, $2\cdot143$, and $2\cdot182$.

\begin{figure}[H]\centering
        \setlength{\tabcolsep}{0.9pt}\tiny\noindent
	\begin{tabular}{r|ccccccccccccccccccccccccccccccccccccccccccccccccccccc}
		1283 & \cellcolor{green!25}1259 & \cellcolor{green!25}1235 & \cellcolor{green!25}1211 & \cellcolor{green!25}1187 & \cellcolor{green!25}1163 & \cellcolor{green!25}1139 & \cellcolor{green!25}1115 & \cellcolor{green!25}1091 & \cellcolor{green!25}1067 & \cellcolor{green!25}1043 & \cellcolor{green!25}1019 & \cellcolor{green!25}995 & \cellcolor{green!25}971 & \cellcolor{green!25}947 & \cellcolor{green!25}923 & \cellcolor{green!25}899 & \cellcolor{green!25}875 & \cellcolor{green!25}851 & \cellcolor{green!25}827 & \cellcolor{green!25}803 & \cellcolor{green!25}779 & \cellcolor{green!25}755 & \cellcolor{green!25}731 & \cellcolor{green!25}707 & \cellcolor{green!25}683 & \cellcolor{green!25}659 & \cellcolor{green!25}635 & \cellcolor{green!25}611 & \cellcolor{green!25}587 & 563 & 539\\
		1229 & \cellcolor{green!25}1205 & \cellcolor{green!25}1181 & \cellcolor{green!25}1157 & \cellcolor{green!25}1133 & \cellcolor{green!25}1109 & \cellcolor{green!25}1085 & \cellcolor{green!25}1061 & \cellcolor{green!25}1037 & \cellcolor{green!25}1013 & \cellcolor{green!25}989 & \cellcolor{green!25}965 & \cellcolor{green!25}941 & \cellcolor{green!25}917 & \cellcolor{green!25}893 & \cellcolor{green!25}869 & \cellcolor{green!25}845 & \cellcolor{green!25}821 & \cellcolor{green!25}797 & \cellcolor{green!25}773 & \cellcolor{green!25}749 & \cellcolor{green!25}725 & \cellcolor{green!25}701 & \cellcolor{green!25}677 & \cellcolor{green!25}653 & \cellcolor{green!25}629 & \cellcolor{green!25}605 & \cellcolor{green!25}581 & \cellcolor{green!25}557 & \cellcolor{green!25}533 & 509 & 485\\
		1175 & \cellcolor{green!25}1151 & \cellcolor{green!25}1127 & \cellcolor{green!25}1103 & \cellcolor{green!25}1079 & \cellcolor{green!25}1055 & \cellcolor{green!25}1031 & \cellcolor{green!25}1007 & \cellcolor{green!25}983 & \cellcolor{green!25}959 & \cellcolor{green!25}935 & \cellcolor{green!25}911 & \cellcolor{green!25}887 & \cellcolor{green!25}863 & \cellcolor{green!25}839 & \cellcolor{green!25}815 & \cellcolor{green!25}791 & \cellcolor{green!25}767 & \cellcolor{green!25}743 & \cellcolor{green!25}719 & \cellcolor{green!25}695 & \cellcolor{green!25}671 & \cellcolor{green!25}647 & \cellcolor{green!25}623 & \cellcolor{green!25}599 & \cellcolor{green!25}575 & \cellcolor{green!25}551 & \cellcolor{green!25}527 & \cellcolor{green!25}503 & 479 & 455 & 431\\
		1121 & \cellcolor{green!25}1097 & \cellcolor{green!25}1073 & \cellcolor{green!25}1049 & \cellcolor{green!25}1025 & \cellcolor{green!25}1001 & \cellcolor{green!25}977 & \cellcolor{green!25}953 & \cellcolor{green!25}929 & \cellcolor{green!25}905 & \cellcolor{green!25}881 & \cellcolor{green!25}857 & \cellcolor{green!25}833 & \cellcolor{green!25}809 & \cellcolor{green!25}785 & \cellcolor{green!25}761 & \cellcolor{green!25}737 & \cellcolor{green!25}713 & \cellcolor{green!25}689 & \cellcolor{green!25}665 & \cellcolor{green!25}641 & \cellcolor{green!25}617 & \cellcolor{green!25}593 & \cellcolor{green!25}569 & \cellcolor{green!25}545 & \cellcolor{green!25}521 & \cellcolor{green!25}497 & \cellcolor{green!25}473 & 449 & 425 & 401 & 377\\
		1059 & \cellcolor{green!25}1035 & \cellcolor{green!25}1011 & \cellcolor{green!25}987 & \cellcolor{green!25}963 & \cellcolor{green!25}939 & \cellcolor{green!25}915 & \cellcolor{green!25}891 & \cellcolor{green!25}867 & \cellcolor{green!25}843 & \cellcolor{green!25}819 & \cellcolor{green!25}795 & \cellcolor{green!25}771 & \cellcolor{green!25}747 & \cellcolor{green!25}723 & \cellcolor{green!25}699 & \cellcolor{green!25}675 & \cellcolor{green!25}651 & \cellcolor{green!25}627 & \cellcolor{green!25}603 & \cellcolor{green!25}579 & \cellcolor{green!25}555 & \cellcolor{green!25}531 & \cellcolor{green!25}507 & \cellcolor{green!25}483 & \cellcolor{green!25}459 & \cellcolor{green!25}435 & \cellcolor{green!25}411 & \cellcolor{green!25}387 & \cellcolor{green!25}363 & 339 & 315\\
		1005 & \cellcolor{green!25}981 & \cellcolor{green!25}957 & \cellcolor{green!25}933 & \cellcolor{green!25}909 & \cellcolor{green!25}885 & \cellcolor{green!25}861 & \cellcolor{green!25}837 & \cellcolor{green!25}813 & \cellcolor{green!25}789 & \cellcolor{green!25}765 & \cellcolor{green!25}741 & \cellcolor{green!25}717 & \cellcolor{green!25}693 & \cellcolor{green!25}669 & \cellcolor{green!25}645 & \cellcolor{green!25}621 & \cellcolor{green!25}597 & \cellcolor{green!25}573 & \cellcolor{green!25}549 & \cellcolor{green!25}525 & \cellcolor{green!25}501 & \cellcolor{green!25}477 & \cellcolor{green!25}453 & \cellcolor{green!25}429 & \cellcolor{green!25}405 & \cellcolor{green!25}381 & \cellcolor{green!25}357 & \cellcolor{green!25}333 & \cellcolor{green!25}309 & 285 & 261\\
		951 & \cellcolor{green!25}927 & \cellcolor{green!25}903 & \cellcolor{green!25}879 & \cellcolor{green!25}855 & \cellcolor{green!25}831 & \cellcolor{green!25}807 & \cellcolor{green!25}783 & \cellcolor{green!25}759 & \cellcolor{green!25}735 & \cellcolor{green!25}711 & \cellcolor{green!25}687 & \cellcolor{green!25}663 & \cellcolor{green!25}639 & \cellcolor{green!25}615 & \cellcolor{green!25}591 & \cellcolor{green!25}567 & \cellcolor{green!25}543 & \cellcolor{green!25}519 & \cellcolor{green!25}495 & \cellcolor{green!25}471 & \cellcolor{green!25}447 & \cellcolor{green!25}423 & \cellcolor{green!25}399 & \cellcolor{green!25}375 & \cellcolor{green!25}351 & \cellcolor{green!25}327 & \cellcolor{green!25}303 & \cellcolor{green!25}279 & 255 & 231 & 207\\
		897 & \cellcolor{green!25}873 & \cellcolor{green!25}849 & \cellcolor{green!25}825 & \cellcolor{green!25}801 & \cellcolor{green!25}777 & \cellcolor{green!25}753 & \cellcolor{green!25}729 & \cellcolor{green!25}705 & \cellcolor{green!25}681 & \cellcolor{green!25}657 & \cellcolor{green!25}633 & \cellcolor{green!25}609 & \cellcolor{green!25}585 & \cellcolor{green!25}561 & \cellcolor{green!25}537 & \cellcolor{green!25}513 & \cellcolor{green!25}489 & \cellcolor{green!25}465 & \cellcolor{green!25}441 & \cellcolor{green!25}417 & \cellcolor{green!25}393 & \cellcolor{green!25}369 & \cellcolor{green!25}345 & \cellcolor{green!25}321 & \cellcolor{green!25}297 & \cellcolor{green!25}273 & \cellcolor{green!25}249 & 225 & 201 & 177 & 153\\
		835 & \cellcolor{green!25}811 & \cellcolor{green!25}787 & \cellcolor{green!25}763 & \cellcolor{green!25}739 & \cellcolor{green!25}715 & \cellcolor{green!25}691 & \cellcolor{green!25}667 & \cellcolor{green!25}643 & 619 & 595 & 571 & 547 & 523 & 499 & 475 & 451 & 427 & 403 & 379 & 355 & 331 & 307 & 283 & 259 & 235 & 211 & 187 & 163 & 139 & 115 & 91\\
		781 & \cellcolor{green!25}757 & \cellcolor{green!25}733 & \cellcolor{green!25}709 & \cellcolor{green!25}685 & \cellcolor{green!25}661 & \cellcolor{green!25}637 & \cellcolor{green!25}613 & 589 & 565 & 541 & 517 & 493 & 469 & 445 & 421 & 397 & 373 & 349 & 325 & 301 & 277 & 253 & 229 & 205 & 181 & 157 & 133 & 109 & 85 & 61 & 37\\
		727 & \cellcolor{green!25}703 & \cellcolor{green!25}679 & \cellcolor{green!25}655 & \cellcolor{green!25}631 & \cellcolor{green!25}607 & \cellcolor{green!25}583 & \cellcolor{green!25}559 & 535 & 511 & 487 & 463 & 439 & 415 & 391 & 367 & 343 & 319 & 295 & 271 & 247 & 223 & 199 & 175 & 151 & 127 & 103 & 79 & 55 & 31 & 7 \\
		673 & 649 & 625 & 601 & 577 & 553 & 529 & 505 & 481 & 457 & 433 & 409 & 385 & 361 & 337 & 313 & 289 & 265 & 241 & 217 & 193 & 169 & 145 & 121 & 97 & 73 & 49 & 25 & 1 \\
		610 & \cellcolor{green!25}586 & \cellcolor{green!25}562 & \cellcolor{green!25}538 & \cellcolor{green!25}514 & \cellcolor{green!25}490 & \cellcolor{green!25}466 & \cellcolor{green!25}442 & \cellcolor{green!25}418 & 394 & 370 & 346 & 322 & 298 & 274 & 250 & 226 & 202 & 178 & 154 & 130 & 106 & 82 & 58 & 34 & 10 \\
		556 & \cellcolor{green!25}532 & \cellcolor{green!25}508 & \cellcolor{green!25}484 & \cellcolor{green!25}460 & \cellcolor{green!25}436 & \cellcolor{green!25}412 & \cellcolor{green!25}388 & 364 & 340 & 316 & 292 & 268 & 244 & 220 & 196 & 172 & 148 & 124 & 100 & 76 & 52 & 28 & 4 \\
		502 & \cellcolor{green!25}478 & \cellcolor{green!25}454 & \cellcolor{green!25}430 & \cellcolor{green!25}406 & \cellcolor{green!25}382 & \cellcolor{green!25}358 & \cellcolor{green!25}334 & 310 & 286 & 262 & 238 & 214 & 190 & 166 & 142 & 118 & 94 & 70 & 46 & 22 \\
		448 & 424 & 400 & 376 & 352 & 328 & 304 & 280 & 256 & 232 & 208 & 184 & 160 & 136 & 112 & 88 & 64 & 40 & 16 \\
		386 & \cellcolor{green!25}362 & \cellcolor{green!25}338 & \cellcolor{green!25}314 & \cellcolor{green!25}290 & \cellcolor{green!25}266 & \cellcolor{green!25}242 & \cellcolor{green!25}218 & \cellcolor{green!25}194 & 170 & 146 & 122 & 98 & 74 & 50 & 26 & 2 \\
		332 & \cellcolor{green!25}308 & \cellcolor{green!25}284 & \cellcolor{green!25}260 & \cellcolor{green!25}236 & \cellcolor{green!25}212 & \cellcolor{green!25}188 & \cellcolor{green!25}164 & 140 & 116 & 92 & 68 & 44 & 20 \\
		278 & \cellcolor{green!25}254 & \cellcolor{green!25}230 & \cellcolor{green!25}206 & \cellcolor{green!25}182 & \cellcolor{green!25}158 & \cellcolor{green!25}134 & \cellcolor{green!25}110 & 86 & 62 & 38 & 14 \\
		224 & 200 & 176 & 152 & 128 & 104 & 80 & 56 & 32 & 8 \\
		162 & \cellcolor{green!25}138 & 114 & 90 & 66 & 42 & 18 \\
		108 & \cellcolor{green!25}84 & 60 & 36 & 12 \\
		54 & 30 & 6 \\
		0 &  \\
	\end{tabular}
\end{figure}

\fex

\section{Case~2: horizontal decomposition}\label{sec:caso-horizontal}

\subsection{Construction of the horizontal blocks}

Throughout this section, we assume $g\geq2$. The case $\Lambda_1\setminus\Gamma_1=\emptyset$ corresponds to an empty decomposition. We therefore assume that $\Lambda_1\setminus\Gamma_1\neq\emptyset$ and that, in the canonical representations of its minimal generators $\delta_0<\cdots<\delta_t$,
\[
 \delta_i=s_i\frac{v_1}{\varepsilon_1}
          -r_i\frac{v_0}{\varepsilon_1},\qquad
 1\leq s_i<\frac{v_0}{\varepsilon_1},\quad r_i\geq1,
\]
the coefficients satisfy
\[
 s_i>s_j,\quad r_i>r_j\quad(i<j).
\]
Set $m=v_1/\varepsilon_1$ and $n_1=v_0/\varepsilon_1$. Thus, $\Gamma_1=\langle n_1,m\rangle$, and $\mu_1=(n_1-1)(m-1)$ is the conductor computed on the scale of $f_1$. We retain the choice of representatives whose coefficients of $dy$ have maximal value, as in \eqref{eq:coeficiente-maximal}. The extension to the tower of semiroots will be treated in Corollary~\ref{cor:determinacao-torre-binaria-horizontal}.

We adopt the convention $s_{-1}=n_1$ and define, for $0\leq i\leq t$,
\[
 A_i=\big[(\Lambda_1\setminus\Gamma_1)\cap(\delta_i+\Gamma_1)\big]
       \setminus\bigcup_{k=0}^{i-1}A_k,
\]
with the union empty for $i=0$. Since every special gap belongs to $\delta_i+\Gamma_1$ for some $i$, the sets $A_i$ form a partition of $\Lambda_1\setminus\Gamma_1$:
\[
 A_i\cap A_j=\emptyset\quad(i\neq j),\qquad
 \bigcup_{i=0}^t A_i=\Lambda_1\setminus\Gamma_1.
\]

Under the assumption of this case, the blocks admit the description
\[
 A_i=\left\{\delta_i+a_i\frac{v_1}{\varepsilon_1}
                       +b_i\frac{v_0}{\varepsilon_1}:
 \begin{array}{l}
  a_i,b_i\in\mathbb Z_{\geq0},\\
  a_i\leq s_{i-1}-s_i-1,\quad b_i\leq r_i-1
 \end{array}\right\}.
\]
Indeed, the gaps in $\delta_i+\Gamma_1$ are the values $s\frac{v_1}{\varepsilon_1}-r\frac{v_0}{\varepsilon_1}$ with $s_i\leq s\leq n_1-1$ and $1\leq r\leq r_i$. For $i\geq1$, since $s_k\geq s_{i-1}$ and $r_k\geq r_{i-1}>r_i$ for every $k<i$, the union of the preceding blocks removes precisely the elements with $s\geq s_{i-1}$. Therefore, $A_i$ corresponds to the horizontal strip
\[
 s_i\leq s\leq s_{i-1}-1,\qquad 1\leq r\leq r_i,
\]
which also applies to $i=0$ by the convention $s_{-1}=n_1$. In particular,
\[
 \min A_i=\delta_i,\qquad |A_i|=(s_{i-1}-s_i)r_i.
\]

The sets $A_i$ are on the value scale of $f_1$. On the scale of $f$, the corresponding blocks are
\[
 \varepsilon_1 A_i=\{\varepsilon_1\delta:\delta\in A_i\},
 \qquad \min(\varepsilon_1 A_i)=\varepsilon_1\delta_i.
\]

For each $i=0,\ldots,t$, define
\[
 \delta_i'=\beta_2+(n_1-s_{i-1})v_1+v_0
\]
and
\[
 B_i=\left\{\delta_i'+a_i v_1+b_i v_0:
 \begin{array}{l}
  a_i,b_i\in\mathbb Z_{\geq0},\\
  a_i\leq s_{i-1}-s_i-1,\quad b_i\leq r_i-1
 \end{array}\right\}.
\]
Since $\varepsilon_1(A_i-\delta_i)$ consists of the values $a_i v_1+b_i v_0$ subject to the same bounds, we obtain
\[
 B_i=\delta_i'+\varepsilon_1(A_i-\delta_i),
 \qquad \min B_i=\delta_i'.
\]

The relation $v_2=(n_1-1)v_1+\beta_2$ gives a geometric identification of the blocks. For $\delta=\delta_i+a\,v_1/\varepsilon_1+b\,v_0/\varepsilon_1\in A_i$, we have
\[
 v_2-\varepsilon_1\delta
 =\delta_i'+(s_{i-1}-s_i-1-a)v_1+(r_i-1-b)v_0.
\]
Varying $a$ and $b$ gives
\begin{equation}\label{eq:reflexao-blocos-horizontais}
 B_i=v_2-\varepsilon_1A_i.
\end{equation}
In particular, $|B_i|=|A_i|$, and the $B_i$ are pairwise disjoint. Moreover, $\varepsilon_1A_i$ and $B_i$ belong to distinct classes modulo $\varepsilon_1$.

Both blocks consist of gaps of $\Gamma_f$ smaller than $v_2$. Indeed, $\max A_i\leq\mu_1-1$ and $v_2-\varepsilon_1(\mu_1-1)=\beta_2+v_0>0$; below $v_2$, the elements of $\Gamma_f$ belong to $\varepsilon_1\Gamma_1$. Hence $\varepsilon_1A_i$ is disjoint from $\Gamma_f$. In turn, each element of $B_i$ is positive, smaller than $v_2$, and congruent to $\beta_2$ modulo $\varepsilon_1$, and is therefore not divisible by $\varepsilon_1$. The selection of the blocks containing special gaps will be established in Corollary~\ref{lem:selecao-blocos-horizontais}.

\subsection{Example}

\ex\label{ex:blocos-horizontais} Consider $f_1=(t^{14}, t^{25}+t^{31})$. In this case, $\frac{v_0}{\varepsilon_1}=14$ and $\frac{v_1}{\varepsilon_1}=25$, with
\[
 \begin{aligned}
  (\delta_0,\ldots,\delta_4)&=(45,76,107,138,169),\\
  (s_0,\ldots,s_4)&=(13,12,11,10,9),\\
  (r_0,\ldots,r_4)&=(20,16,12,8,4).
 \end{aligned}
\]
Thus, $A_0,\ldots,A_4$ are represented by the brown, blue, red, orange, and gray strips, respectively; since $s_{i-1}-s_i=1$ for each $i$, each block occupies a single row of the table.

In this example,
\[
 (n_1-s_{-1},n_1-s_0,n_1-s_1,n_1-s_2,n_1-s_3)=(0,1,2,3,4).
\]
Consequently,
\[
 \begin{aligned}
  \min B_0&=\beta_2+v_0,\\
  \min B_1&=\beta_2+v_1+v_0,\\
  \min B_2&=\beta_2+2v_1+v_0,\\
  \min B_3&=\beta_2+3v_1+v_0,\\
  \min B_4&=\beta_2+4v_1+v_0.
 \end{aligned}
\]
Since $s_{i-1}-s_i=1$, we have
\[
 B_i=\{\min B_i+b_i v_0:0\leq b_i\leq r_i-1\},
\]
with widths $20,16,12,8$, and $4$, respectively.

\begin{figure}[H]\centering
\setlength{\tabcolsep}{1pt}\small\noindent
\begin{tabular}{r|*{23}{c}}
325 & \cellcolor{brown!25}311 & \cellcolor{brown!25}297 & \cellcolor{brown!25}283 & \cellcolor{brown!25}269 & \cellcolor{brown!25}255 & \cellcolor{brown!25}241 & \cellcolor{brown!25}227 & \cellcolor{brown!25}213 & \cellcolor{brown!25}199 & \cellcolor{brown!25}185 & \cellcolor{brown!25}171 & \cellcolor{brown!25}157 & \cellcolor{brown!25}143 & \cellcolor{brown!25}129 & \cellcolor{brown!25}115 & \cellcolor{brown!25}101 & \cellcolor{brown!25}87 & \cellcolor{brown!25}73 & \cellcolor{brown!25}59 & \cellcolor{brown!50}45 & 31 & 17 & 3 \\
300 & \cellcolor{blue!25}286 & \cellcolor{blue!25}272 & \cellcolor{blue!25}258 & \cellcolor{blue!25}244 & \cellcolor{blue!25}230 & \cellcolor{blue!25}216 & \cellcolor{blue!25}202 & \cellcolor{blue!25}188 & \cellcolor{blue!25}174 & \cellcolor{blue!25}160 & \cellcolor{blue!25}146 & \cellcolor{blue!25}132 & \cellcolor{blue!25}118 & \cellcolor{blue!25}104 & \cellcolor{blue!25}90 & \cellcolor{blue!50}76 & 62 & 48 & 34 & 20 & 6 \\
275 & \cellcolor{red!25}261 & \cellcolor{red!25}247 & \cellcolor{red!25}233 & \cellcolor{red!25}219 & \cellcolor{red!25}205 & \cellcolor{red!25}191 & \cellcolor{red!25}177 & \cellcolor{red!25}163 & \cellcolor{red!25}149 & \cellcolor{red!25}135 & \cellcolor{red!25}121 & \cellcolor{red!50}107 & 93 & 79 & 65 & 51 & 37 & 23 & 9 \\
250 & \cellcolor{orange!25}236 & \cellcolor{orange!25}222 & \cellcolor{orange!25}208 & \cellcolor{orange!25}194 & \cellcolor{orange!25}180 & \cellcolor{orange!25}166 & \cellcolor{orange!25}152 & \cellcolor{orange!50}138 & 124 & 110 & 96 & 82 & 68 & 54 & 40 & 26 & 12 \\
225 & \cellcolor{gray!25}211 & \cellcolor{gray!25}197 & \cellcolor{gray!25}183 & \cellcolor{gray!50}169 & 155 & 141 & 127 & 113 & 99 & 85 & 71 & 57 & 43 & 29 & 15 & 1 \\
200 & 186 & 172 & 158 & 144 & 130 & 116 & 102 & 88 & 74 & 60 & 46 & 32 & 18 & 4 \\
175 & 161 & 147 & 133 & 119 & 105 & 91 & 77 & 63 & 49 & 35 & 21 & 7 \\
150 & 136 & 122 & 108 & 94 & 80 & 66 & 52 & 38 & 24 & 10 \\
125 & 111 & 97 & 83 & 69 & 55 & 41 & 27 & 13 \\
100 & 86 & 72 & 58 & 44 & 30 & 16 & 2 \\
75 & 61 & 47 & 33 & 19 & 5 \\
50 & 36 & 22 & 8 \\
25 & 11 \\
0 &  \\
\end{tabular}
\end{figure}

Consider also the curve
\[
 f=(t^{28},t^{50}+t^{62}+t^{87}).
\]
In this case, $\varepsilon_1=2$, $\beta_2=87$, and $v_2=737$. In the following table, the special gaps of the forms $v_2+\varepsilon_1\delta$, with $\delta\in\mathbb N\setminus\Gamma_1$, and $v_2-\varepsilon_1\delta$, with $\delta\in\mathbb N^*\setminus\Lambda_1$, are highlighted in green. The other colored blocks are $\varepsilon_1A_0,\ldots,\varepsilon_1A_4$, shown in brown, blue, red, orange, and gray, respectively. These blocks are merely the candidates arising from the first alternative of $\rho_1$, so their elements are not necessarily special gaps of $f$.

\begin{figure}[H]\centering
\setlength{\tabcolsep}{0.9pt}\scriptsize\noindent
\begin{tabular}{r|*{26}{c}}
1387 & \cellcolor{green!25}1359 & \cellcolor{green!25}1331 & \cellcolor{green!25}1303 & \cellcolor{green!25}1275 & \cellcolor{green!25}1247 & \cellcolor{green!25}1219 & \cellcolor{green!25}1191 & \cellcolor{green!25}1163 & \cellcolor{green!25}1135 & \cellcolor{green!25}1107 & \cellcolor{green!25}1079 & \cellcolor{green!25}1051 & \cellcolor{green!25}1023 & \cellcolor{green!25}995 & \cellcolor{green!25}967 & \cellcolor{green!25}939 & \cellcolor{green!25}911 & \cellcolor{green!25}883 & \cellcolor{green!25}855 & \cellcolor{green!25}827 & \cellcolor{green!25}799 & \cellcolor{green!25}771 & \cellcolor{green!25}743 & \cellcolor{green!25}715 & 687 & 659 \\
1337 & \cellcolor{green!25}1309 & \cellcolor{green!25}1281 & \cellcolor{green!25}1253 & \cellcolor{green!25}1225 & \cellcolor{green!25}1197 & \cellcolor{green!25}1169 & \cellcolor{green!25}1141 & \cellcolor{green!25}1113 & \cellcolor{green!25}1085 & \cellcolor{green!25}1057 & \cellcolor{green!25}1029 & \cellcolor{green!25}1001 & \cellcolor{green!25}973 & \cellcolor{green!25}945 & \cellcolor{green!25}917 & \cellcolor{green!25}889 & \cellcolor{green!25}861 & \cellcolor{green!25}833 & \cellcolor{green!25}805 & \cellcolor{green!25}777 & \cellcolor{green!25}749 & \cellcolor{green!25}721 & \cellcolor{green!25}693 & \cellcolor{green!25}665 & 637 & 609 \\
1287 & \cellcolor{green!25}1259 & \cellcolor{green!25}1231 & \cellcolor{green!25}1203 & \cellcolor{green!25}1175 & \cellcolor{green!25}1147 & \cellcolor{green!25}1119 & \cellcolor{green!25}1091 & \cellcolor{green!25}1063 & \cellcolor{green!25}1035 & \cellcolor{green!25}1007 & \cellcolor{green!25}979 & \cellcolor{green!25}951 & \cellcolor{green!25}923 & \cellcolor{green!25}895 & \cellcolor{green!25}867 & \cellcolor{green!25}839 & \cellcolor{green!25}811 & \cellcolor{green!25}783 & \cellcolor{green!25}755 & \cellcolor{green!25}727 & \cellcolor{green!25}699 & \cellcolor{green!25}671 & \cellcolor{green!25}643 & \cellcolor{green!25}615 & 587 & 559 \\
1237 & \cellcolor{green!25}1209 & \cellcolor{green!25}1181 & \cellcolor{green!25}1153 & \cellcolor{green!25}1125 & \cellcolor{green!25}1097 & \cellcolor{green!25}1069 & \cellcolor{green!25}1041 & \cellcolor{green!25}1013 & \cellcolor{green!25}985 & \cellcolor{green!25}957 & \cellcolor{green!25}929 & \cellcolor{green!25}901 & \cellcolor{green!25}873 & \cellcolor{green!25}845 & \cellcolor{green!25}817 & \cellcolor{green!25}789 & \cellcolor{green!25}761 & \cellcolor{green!25}733 & \cellcolor{green!25}705 & \cellcolor{green!25}677 & \cellcolor{green!25}649 & \cellcolor{green!25}621 & \cellcolor{green!25}593 & \cellcolor{green!25}565 & 537 & 509 \\
1187 & \cellcolor{green!25}1159 & \cellcolor{green!25}1131 & \cellcolor{green!25}1103 & \cellcolor{green!25}1075 & \cellcolor{green!25}1047 & \cellcolor{green!25}1019 & \cellcolor{green!25}991 & \cellcolor{green!25}963 & \cellcolor{green!25}935 & \cellcolor{green!25}907 & \cellcolor{green!25}879 & \cellcolor{green!25}851 & \cellcolor{green!25}823 & \cellcolor{green!25}795 & \cellcolor{green!25}767 & \cellcolor{green!25}739 & \cellcolor{green!25}711 & \cellcolor{green!25}683 & \cellcolor{green!25}655 & \cellcolor{green!25}627 & \cellcolor{green!25}599 & \cellcolor{green!25}571 & \cellcolor{green!25}543 & \cellcolor{green!25}515 & 487 & 459 \\
1137 & \cellcolor{green!25}1109 & \cellcolor{green!25}1081 & \cellcolor{green!25}1053 & \cellcolor{green!25}1025 & \cellcolor{green!25}997 & \cellcolor{green!25}969 & \cellcolor{green!25}941 & \cellcolor{green!25}913 & \cellcolor{green!25}885 & \cellcolor{green!25}857 & \cellcolor{green!25}829 & \cellcolor{green!25}801 & \cellcolor{green!25}773 & \cellcolor{green!25}745 & \cellcolor{green!25}717 & \cellcolor{green!25}689 & \cellcolor{green!25}661 & \cellcolor{green!25}633 & \cellcolor{green!25}605 & \cellcolor{green!25}577 & \cellcolor{green!25}549 & \cellcolor{green!25}521 & \cellcolor{green!25}493 & \cellcolor{green!25}465 & 437 & 409 \\
1087 & \cellcolor{green!25}1059 & \cellcolor{green!25}1031 & \cellcolor{green!25}1003 & \cellcolor{green!25}975 & \cellcolor{green!25}947 & \cellcolor{green!25}919 & \cellcolor{green!25}891 & \cellcolor{green!25}863 & \cellcolor{green!25}835 & \cellcolor{green!25}807 & \cellcolor{green!25}779 & \cellcolor{green!25}751 & \cellcolor{green!25}723 & \cellcolor{green!25}695 & \cellcolor{green!25}667 & \cellcolor{green!25}639 & \cellcolor{green!25}611 & \cellcolor{green!25}583 & \cellcolor{green!25}555 & \cellcolor{green!25}527 & \cellcolor{green!25}499 & \cellcolor{green!25}471 & \cellcolor{green!25}443 & \cellcolor{green!25}415 & 387 & 359 \\
1037 & \cellcolor{green!25}1009 & \cellcolor{green!25}981 & \cellcolor{green!25}953 & \cellcolor{green!25}925 & \cellcolor{green!25}897 & \cellcolor{green!25}869 & \cellcolor{green!25}841 & \cellcolor{green!25}813 & \cellcolor{green!25}785 & \cellcolor{green!25}757 & \cellcolor{green!25}729 & \cellcolor{green!25}701 & \cellcolor{green!25}673 & \cellcolor{green!25}645 & \cellcolor{green!25}617 & \cellcolor{green!25}589 & \cellcolor{green!25}561 & \cellcolor{green!25}533 & \cellcolor{green!25}505 & \cellcolor{green!25}477 & \cellcolor{green!25}449 & \cellcolor{green!25}421 & \cellcolor{green!25}393 & \cellcolor{green!25}365 & 337 & 309 \\
987 & \cellcolor{green!25}959 & \cellcolor{green!25}931 & \cellcolor{green!25}903 & \cellcolor{green!25}875 & \cellcolor{green!25}847 & \cellcolor{green!25}819 & \cellcolor{green!25}791 & \cellcolor{green!25}763 & \cellcolor{green!25}735 & \cellcolor{green!25}707 & \cellcolor{green!25}679 & \cellcolor{green!25}651 & \cellcolor{green!25}623 & \cellcolor{green!25}595 & \cellcolor{green!25}567 & \cellcolor{green!25}539 & \cellcolor{green!25}511 & \cellcolor{green!25}483 & \cellcolor{green!25}455 & \cellcolor{green!25}427 & 399 & 371 & 343 & 315 & 287 & 259 \\
937 & \cellcolor{green!25}909 & \cellcolor{green!25}881 & \cellcolor{green!25}853 & \cellcolor{green!25}825 & \cellcolor{green!25}797 & \cellcolor{green!25}769 & \cellcolor{green!25}741 & \cellcolor{green!25}713 & \cellcolor{green!25}685 & \cellcolor{green!25}657 & \cellcolor{green!25}629 & \cellcolor{green!25}601 & \cellcolor{green!25}573 & \cellcolor{green!25}545 & \cellcolor{green!25}517 & \cellcolor{green!25}489 & 461 & 433 & 405 & 377 & 349 & 321 & 293 & 265 & 237 & 209 \\
887 & \cellcolor{green!25}859 & \cellcolor{green!25}831 & \cellcolor{green!25}803 & \cellcolor{green!25}775 & \cellcolor{green!25}747 & \cellcolor{green!25}719 & \cellcolor{green!25}691 & \cellcolor{green!25}663 & \cellcolor{green!25}635 & \cellcolor{green!25}607 & \cellcolor{green!25}579 & \cellcolor{green!25}551 & 523 & 495 & 467 & 439 & 411 & 383 & 355 & 327 & 299 & 271 & 243 & 215 & 187 & 159 \\
837 & \cellcolor{green!25}809 & \cellcolor{green!25}781 & \cellcolor{green!25}753 & \cellcolor{green!25}725 & \cellcolor{green!25}697 & \cellcolor{green!25}669 & \cellcolor{green!25}641 & \cellcolor{green!25}613 & 585 & 557 & 529 & 501 & 473 & 445 & 417 & 389 & 361 & 333 & 305 & 277 & 249 & 221 & 193 & 165 & 137 & 109 \\
787 & \cellcolor{green!25}759 & \cellcolor{green!25}731 & \cellcolor{green!25}703 & \cellcolor{green!25}675 & 647 & 619 & 591 & 563 & 535 & 507 & 479 & 451 & 423 & 395 & 367 & 339 & 311 & 283 & 255 & 227 & 199 & 171 & 143 & 115 & 87 & 59 \\
737 & 709 & 681 & 653 & 625 & 597 & 569 & 541 & 513 & 485 & 457 & 429 & 401 & 373 & 345 & 317 & 289 & 261 & 233 & 205 & 177 & 149 & 121 & 93 & 65 & 37 & 9 \\
650 & \cellcolor{brown!25}622 & \cellcolor{brown!25}594 & \cellcolor{brown!25}566 & \cellcolor{brown!25}538 & \cellcolor{brown!25}510 & \cellcolor{brown!25}482 & \cellcolor{brown!25}454 & \cellcolor{brown!25}426 & \cellcolor{brown!25}398 & \cellcolor{brown!25}370 & \cellcolor{brown!25}342 & \cellcolor{brown!25}314 & \cellcolor{brown!25}286 & \cellcolor{brown!25}258 & \cellcolor{brown!25}230 & \cellcolor{brown!25}202 & \cellcolor{brown!25}174 & \cellcolor{brown!25}146 & \cellcolor{brown!25}118 & \cellcolor{brown!50}90 & 62 & 34 & 6 \\
600 & \cellcolor{blue!25}572 & \cellcolor{blue!25}544 & \cellcolor{blue!25}516 & \cellcolor{blue!25}488 & \cellcolor{blue!25}460 & \cellcolor{blue!25}432 & \cellcolor{blue!25}404 & \cellcolor{blue!25}376 & \cellcolor{blue!25}348 & \cellcolor{blue!25}320 & \cellcolor{blue!25}292 & \cellcolor{blue!25}264 & \cellcolor{blue!25}236 & \cellcolor{blue!25}208 & \cellcolor{blue!25}180 & \cellcolor{blue!50}152 & 124 & 96 & 68 & 40 & 12 \\
550 & \cellcolor{red!25}522 & \cellcolor{red!25}494 & \cellcolor{red!25}466 & \cellcolor{red!25}438 & \cellcolor{red!25}410 & \cellcolor{red!25}382 & \cellcolor{red!25}354 & \cellcolor{red!25}326 & \cellcolor{red!25}298 & \cellcolor{red!25}270 & \cellcolor{red!25}242 & \cellcolor{red!50}214 & 186 & 158 & 130 & 102 & 74 & 46 & 18 \\
500 & \cellcolor{orange!25}472 & \cellcolor{orange!25}444 & \cellcolor{orange!25}416 & \cellcolor{orange!25}388 & \cellcolor{orange!25}360 & \cellcolor{orange!25}332 & \cellcolor{orange!25}304 & \cellcolor{orange!50}276 & 248 & 220 & 192 & 164 & 136 & 108 & 80 & 52 & 24 \\
450 & \cellcolor{gray!25}422 & \cellcolor{gray!25}394 & \cellcolor{gray!25}366 & \cellcolor{gray!50}338 & 310 & 282 & 254 & 226 & 198 & 170 & 142 & 114 & 86 & 58 & 30 & 2 \\
400 & 372 & 344 & 316 & 288 & 260 & 232 & 204 & 176 & 148 & 120 & 92 & 64 & 36 & 8 \\
350 & 322 & 294 & 266 & 238 & 210 & 182 & 154 & 126 & 98 & 70 & 42 & 14 \\
300 & 272 & 244 & 216 & 188 & 160 & 132 & 104 & 76 & 48 & 20 \\
250 & 222 & 194 & 166 & 138 & 110 & 82 & 54 & 26 \\
200 & 172 & 144 & 116 & 88 & 60 & 32 & 4 \\
150 & 122 & 94 & 66 & 38 & 10 \\
100 & 72 & 44 & 16 \\
50 & 22 \\
0 & \\
\end{tabular}
\end{figure}

The comparisons between the minimal elements of the blocks are
\[
\begin{aligned}
 \varepsilon_1\delta_0&=90 <115=\delta_0',\\
 \varepsilon_1\delta_1&=152<165=\delta_1',\\
 \varepsilon_1\delta_2&=214<215=\delta_2',\\
 \varepsilon_1\delta_3&=276>265=\delta_3',\\
 \varepsilon_1\delta_4&=338>315=\delta_4'.
\end{aligned}
\]
Theorem~\ref{teo:decomposicao-rho-blocos-horizontais}, proved below, gives
\[
 \rho_1(\Lambda_1\setminus\Gamma_1)
 =\varepsilon_1A_0\mathbin{\dot\cup}\varepsilon_1A_1
  \mathbin{\dot\cup}\varepsilon_1A_2
  \mathbin{\dot\cup}B_3\mathbin{\dot\cup}B_4.
\]
In the following table, the sets $\varepsilon_1A_0,\varepsilon_1A_1,\varepsilon_1A_2,B_3$, and $B_4$ appear in brown, blue, red, orange, and gray, respectively. By Theorem~\ref{t5.2}, these blocks and the gaps highlighted in green constitute all special gaps of $f$.

\begin{figure}[H]\centering
\setlength{\tabcolsep}{0.9pt}\scriptsize\noindent
\begin{tabular}{r|*{26}{c}}
1387 & \cellcolor{green!25}1359 & \cellcolor{green!25}1331 & \cellcolor{green!25}1303 & \cellcolor{green!25}1275 & \cellcolor{green!25}1247 & \cellcolor{green!25}1219 & \cellcolor{green!25}1191 & \cellcolor{green!25}1163 & \cellcolor{green!25}1135 & \cellcolor{green!25}1107 & \cellcolor{green!25}1079 & \cellcolor{green!25}1051 & \cellcolor{green!25}1023 & \cellcolor{green!25}995 & \cellcolor{green!25}967 & \cellcolor{green!25}939 & \cellcolor{green!25}911 & \cellcolor{green!25}883 & \cellcolor{green!25}855 & \cellcolor{green!25}827 & \cellcolor{green!25}799 & \cellcolor{green!25}771 & \cellcolor{green!25}743 & \cellcolor{green!25}715 & 687 & 659 \\
1337 & \cellcolor{green!25}1309 & \cellcolor{green!25}1281 & \cellcolor{green!25}1253 & \cellcolor{green!25}1225 & \cellcolor{green!25}1197 & \cellcolor{green!25}1169 & \cellcolor{green!25}1141 & \cellcolor{green!25}1113 & \cellcolor{green!25}1085 & \cellcolor{green!25}1057 & \cellcolor{green!25}1029 & \cellcolor{green!25}1001 & \cellcolor{green!25}973 & \cellcolor{green!25}945 & \cellcolor{green!25}917 & \cellcolor{green!25}889 & \cellcolor{green!25}861 & \cellcolor{green!25}833 & \cellcolor{green!25}805 & \cellcolor{green!25}777 & \cellcolor{green!25}749 & \cellcolor{green!25}721 & \cellcolor{green!25}693 & \cellcolor{green!25}665 & 637 & 609 \\
1287 & \cellcolor{green!25}1259 & \cellcolor{green!25}1231 & \cellcolor{green!25}1203 & \cellcolor{green!25}1175 & \cellcolor{green!25}1147 & \cellcolor{green!25}1119 & \cellcolor{green!25}1091 & \cellcolor{green!25}1063 & \cellcolor{green!25}1035 & \cellcolor{green!25}1007 & \cellcolor{green!25}979 & \cellcolor{green!25}951 & \cellcolor{green!25}923 & \cellcolor{green!25}895 & \cellcolor{green!25}867 & \cellcolor{green!25}839 & \cellcolor{green!25}811 & \cellcolor{green!25}783 & \cellcolor{green!25}755 & \cellcolor{green!25}727 & \cellcolor{green!25}699 & \cellcolor{green!25}671 & \cellcolor{green!25}643 & \cellcolor{green!25}615 & 587 & 559 \\
1237 & \cellcolor{green!25}1209 & \cellcolor{green!25}1181 & \cellcolor{green!25}1153 & \cellcolor{green!25}1125 & \cellcolor{green!25}1097 & \cellcolor{green!25}1069 & \cellcolor{green!25}1041 & \cellcolor{green!25}1013 & \cellcolor{green!25}985 & \cellcolor{green!25}957 & \cellcolor{green!25}929 & \cellcolor{green!25}901 & \cellcolor{green!25}873 & \cellcolor{green!25}845 & \cellcolor{green!25}817 & \cellcolor{green!25}789 & \cellcolor{green!25}761 & \cellcolor{green!25}733 & \cellcolor{green!25}705 & \cellcolor{green!25}677 & \cellcolor{green!25}649 & \cellcolor{green!25}621 & \cellcolor{green!25}593 & \cellcolor{green!25}565 & 537 & 509 \\
1187 & \cellcolor{green!25}1159 & \cellcolor{green!25}1131 & \cellcolor{green!25}1103 & \cellcolor{green!25}1075 & \cellcolor{green!25}1047 & \cellcolor{green!25}1019 & \cellcolor{green!25}991 & \cellcolor{green!25}963 & \cellcolor{green!25}935 & \cellcolor{green!25}907 & \cellcolor{green!25}879 & \cellcolor{green!25}851 & \cellcolor{green!25}823 & \cellcolor{green!25}795 & \cellcolor{green!25}767 & \cellcolor{green!25}739 & \cellcolor{green!25}711 & \cellcolor{green!25}683 & \cellcolor{green!25}655 & \cellcolor{green!25}627 & \cellcolor{green!25}599 & \cellcolor{green!25}571 & \cellcolor{green!25}543 & \cellcolor{green!25}515 & 487 & 459 \\
1137 & \cellcolor{green!25}1109 & \cellcolor{green!25}1081 & \cellcolor{green!25}1053 & \cellcolor{green!25}1025 & \cellcolor{green!25}997 & \cellcolor{green!25}969 & \cellcolor{green!25}941 & \cellcolor{green!25}913 & \cellcolor{green!25}885 & \cellcolor{green!25}857 & \cellcolor{green!25}829 & \cellcolor{green!25}801 & \cellcolor{green!25}773 & \cellcolor{green!25}745 & \cellcolor{green!25}717 & \cellcolor{green!25}689 & \cellcolor{green!25}661 & \cellcolor{green!25}633 & \cellcolor{green!25}605 & \cellcolor{green!25}577 & \cellcolor{green!25}549 & \cellcolor{green!25}521 & \cellcolor{green!25}493 & \cellcolor{green!25}465 & 437 & 409 \\
1087 & \cellcolor{green!25}1059 & \cellcolor{green!25}1031 & \cellcolor{green!25}1003 & \cellcolor{green!25}975 & \cellcolor{green!25}947 & \cellcolor{green!25}919 & \cellcolor{green!25}891 & \cellcolor{green!25}863 & \cellcolor{green!25}835 & \cellcolor{green!25}807 & \cellcolor{green!25}779 & \cellcolor{green!25}751 & \cellcolor{green!25}723 & \cellcolor{green!25}695 & \cellcolor{green!25}667 & \cellcolor{green!25}639 & \cellcolor{green!25}611 & \cellcolor{green!25}583 & \cellcolor{green!25}555 & \cellcolor{green!25}527 & \cellcolor{green!25}499 & \cellcolor{green!25}471 & \cellcolor{green!25}443 & \cellcolor{green!25}415 & 387 & 359 \\
1037 & \cellcolor{green!25}1009 & \cellcolor{green!25}981 & \cellcolor{green!25}953 & \cellcolor{green!25}925 & \cellcolor{green!25}897 & \cellcolor{green!25}869 & \cellcolor{green!25}841 & \cellcolor{green!25}813 & \cellcolor{green!25}785 & \cellcolor{green!25}757 & \cellcolor{green!25}729 & \cellcolor{green!25}701 & \cellcolor{green!25}673 & \cellcolor{green!25}645 & \cellcolor{green!25}617 & \cellcolor{green!25}589 & \cellcolor{green!25}561 & \cellcolor{green!25}533 & \cellcolor{green!25}505 & \cellcolor{green!25}477 & \cellcolor{green!25}449 & \cellcolor{green!25}421 & \cellcolor{green!25}393 & \cellcolor{green!25}365 & 337 & 309 \\
987 & \cellcolor{green!25}959 & \cellcolor{green!25}931 & \cellcolor{green!25}903 & \cellcolor{green!25}875 & \cellcolor{green!25}847 & \cellcolor{green!25}819 & \cellcolor{green!25}791 & \cellcolor{green!25}763 & \cellcolor{green!25}735 & \cellcolor{green!25}707 & \cellcolor{green!25}679 & \cellcolor{green!25}651 & \cellcolor{green!25}623 & \cellcolor{green!25}595 & \cellcolor{green!25}567 & \cellcolor{green!25}539 & \cellcolor{green!25}511 & \cellcolor{green!25}483 & \cellcolor{green!25}455 & \cellcolor{green!25}427 & \cellcolor{gray!25}399 & \cellcolor{gray!25}371 & \cellcolor{gray!25}343 & \cellcolor{gray!50}315 & 287 & 259 \\
937 & \cellcolor{green!25}909 & \cellcolor{green!25}881 & \cellcolor{green!25}853 & \cellcolor{green!25}825 & \cellcolor{green!25}797 & \cellcolor{green!25}769 & \cellcolor{green!25}741 & \cellcolor{green!25}713 & \cellcolor{green!25}685 & \cellcolor{green!25}657 & \cellcolor{green!25}629 & \cellcolor{green!25}601 & \cellcolor{green!25}573 & \cellcolor{green!25}545 & \cellcolor{green!25}517 & \cellcolor{green!25}489 & \cellcolor{orange!25}461 & \cellcolor{orange!25}433 & \cellcolor{orange!25}405 & \cellcolor{orange!25}377 & \cellcolor{orange!25}349 & \cellcolor{orange!25}321 & \cellcolor{orange!25}293 & \cellcolor{orange!50}265 & 237 & 209 \\
887 & \cellcolor{green!25}859 & \cellcolor{green!25}831 & \cellcolor{green!25}803 & \cellcolor{green!25}775 & \cellcolor{green!25}747 & \cellcolor{green!25}719 & \cellcolor{green!25}691 & \cellcolor{green!25}663 & \cellcolor{green!25}635 & \cellcolor{green!25}607 & \cellcolor{green!25}579 & \cellcolor{green!25}551 & 523 & 495 & 467 & 439 & 411 & 383 & 355 & 327 & 299 & 271 & 243 & 215 & 187 & 159 \\
837 & \cellcolor{green!25}809 & \cellcolor{green!25}781 & \cellcolor{green!25}753 & \cellcolor{green!25}725 & \cellcolor{green!25}697 & \cellcolor{green!25}669 & \cellcolor{green!25}641 & \cellcolor{green!25}613 & 585 & 557 & 529 & 501 & 473 & 445 & 417 & 389 & 361 & 333 & 305 & 277 & 249 & 221 & 193 & 165 & 137 & 109 \\
787 & \cellcolor{green!25}759 & \cellcolor{green!25}731 & \cellcolor{green!25}703 & \cellcolor{green!25}675 & 647 & 619 & 591 & 563 & 535 & 507 & 479 & 451 & 423 & 395 & 367 & 339 & 311 & 283 & 255 & 227 & 199 & 171 & 143 & 115 & 87 & 59 \\
737 & 709 & 681 & 653 & 625 & 597 & 569 & 541 & 513 & 485 & 457 & 429 & 401 & 373 & 345 & 317 & 289 & 261 & 233 & 205 & 177 & 149 & 121 & 93 & 65 & 37 & 9 \\
650 & \cellcolor{brown!25}622 & \cellcolor{brown!25}594 & \cellcolor{brown!25}566 & \cellcolor{brown!25}538 & \cellcolor{brown!25}510 & \cellcolor{brown!25}482 & \cellcolor{brown!25}454 & \cellcolor{brown!25}426 & \cellcolor{brown!25}398 & \cellcolor{brown!25}370 & \cellcolor{brown!25}342 & \cellcolor{brown!25}314 & \cellcolor{brown!25}286 & \cellcolor{brown!25}258 & \cellcolor{brown!25}230 & \cellcolor{brown!25}202 & \cellcolor{brown!25}174 & \cellcolor{brown!25}146 & \cellcolor{brown!25}118 & \cellcolor{brown!50}90 & 62 & 34 & 6 \\
600 & \cellcolor{blue!25}572 & \cellcolor{blue!25}544 & \cellcolor{blue!25}516 & \cellcolor{blue!25}488 & \cellcolor{blue!25}460 & \cellcolor{blue!25}432 & \cellcolor{blue!25}404 & \cellcolor{blue!25}376 & \cellcolor{blue!25}348 & \cellcolor{blue!25}320 & \cellcolor{blue!25}292 & \cellcolor{blue!25}264 & \cellcolor{blue!25}236 & \cellcolor{blue!25}208 & \cellcolor{blue!25}180 & \cellcolor{blue!50}152 & 124 & 96 & 68 & 40 & 12 \\
550 & \cellcolor{red!25}522 & \cellcolor{red!25}494 & \cellcolor{red!25}466 & \cellcolor{red!25}438 & \cellcolor{red!25}410 & \cellcolor{red!25}382 & \cellcolor{red!25}354 & \cellcolor{red!25}326 & \cellcolor{red!25}298 & \cellcolor{red!25}270 & \cellcolor{red!25}242 & \cellcolor{red!50}214 & 186 & 158 & 130 & 102 & 74 & 46 & 18 \\
500 & 472 & 444 & 416 & 388 & 360 & 332 & 304 & 276 & 248 & 220 & 192 & 164 & 136 & 108 & 80 & 52 & 24 \\
450 & 422 & 394 & 366 & 338 & 310 & 282 & 254 & 226 & 198 & 170 & 142 & 114 & 86 & 58 & 30 & 2 \\
400 & 372 & 344 & 316 & 288 & 260 & 232 & 204 & 176 & 148 & 120 & 92 & 64 & 36 & 8 \\
350 & 322 & 294 & 266 & 238 & 210 & 182 & 154 & 126 & 98 & 70 & 42 & 14 \\
300 & 272 & 244 & 216 & 188 & 160 & 132 & 104 & 76 & 48 & 20 \\
250 & 222 & 194 & 166 & 138 & 110 & 82 & 54 & 26 \\
200 & 172 & 144 & 116 & 88 & 60 & 32 & 4 \\
150 & 122 & 94 & 66 & 38 & 10 \\
100 & 72 & 44 & 16 \\
50 & 22 \\
0 & \\
\end{tabular}
\end{figure}

\fex

\subsection{Maximal coefficients, the image of the transition map, and consequences}

\lem\label{lem:theta-blocos-horizontais}
Under the assumptions of Case~2, set
\[
 \vartheta_i=n_1+(n_1-s_{i-1})m.
\]
For $0\leq i\leq t$ and $\delta\in A_i$, we have
\begin{equation}\label{eq:theta-blocos-horizontais}
 \nu_1(B_{\delta})=\vartheta_i+\delta-\delta_i.
\end{equation}
Moreover, the sequence
\[
 T_i=\varepsilon_1(\delta_i-\vartheta_i)
     =\varepsilon_1\delta_i-(n_1-s_{i-1})v_1-v_0,
 \qquad i=0,\ldots,t,
\]
is strictly increasing.
\flem

\dem By the description of the horizontal blocks,
\[
 M_i=\max A_i=(s_{i-1}-1)m-n_1.
\]
The minima $\delta_i$ are strictly increasing, and the maxima $M_i$ are strictly decreasing. Lemma~\ref{lem:coeficientes-particao-retangular} gives
\[
 \nu_1(B_\delta)=K-M_i+\delta-\delta_i
 =\vartheta_i+\delta-\delta_i,
\]
since $K-M_i=n_1+(n_1-s_{i-1})m=\vartheta_i$.

Finally, for $0\leq i<t$, consider
\[
 \eta_i=\delta_{i+1}+(s_i-s_{i+1})m
       =\delta_i+(r_i-r_{i+1})n_1\in A_i.
\]
By the identities and inequalities already proved,
\[
 \vartheta_i+(r_i-r_{i+1})n_1
 =\nu_1(B_{\eta_i})
 \geq\vartheta_{i+1}+(s_i-s_{i+1})m.
\]
Consequently,
\[
 \delta_{i+1}-\vartheta_{i+1}
 \geq\delta_i-\vartheta_i.
\]
Equality would imply
\[
 (r_i-r_{i+1})n_1=(s_{i-1}-s_{i+1})m.
\]
This is impossible, since $\gcd(n_1,m)=1$ and $0<s_{i-1}-s_{i+1}<n_1$. Therefore, $T_{i+1}>T_i$.\f

\lem\label{lem:transicao-blocos-horizontais}
Under the assumption of Case~2, if
\[
 \varepsilon_1\delta_j>\delta_j'
 =\beta_2+(n_1-s_{j-1})v_1+v_0,
\]
then $B_k\subset\Lambda_f\setminus\Gamma_f$ for every $k\geq j$.
\flem

\dem By Lemma~\ref{lem:theta-blocos-horizontais},
\[
 \beta_2+\varepsilon_1\nu_1(B_{\delta_j})
 =\beta_2+(n_1-s_{j-1})v_1+v_0=\delta_j'.
\]
The definition of $\rho_1$ therefore gives
\[
 \rho_1(\delta_j)=\delta_j'\in\Lambda_f\setminus\Gamma_f.
\]
Choose a differential $\omega$ with $\nu_f(\omega)=\delta_j'$. For $k\geq j$, we have
\[
 \delta_k'-\delta_j'=(s_{j-1}-s_{k-1})v_1,
 \qquad s_{j-1}-s_{k-1}\geq0.
\]
Thus, $\omega_k=y^{s_{j-1}-s_{k-1}}\omega$ satisfies $\nu_f(\omega_k)=\delta_k'$. For each $\delta_k'+av_1+bv_0\in B_k$, it follows that
\[
 \nu_f(x^by^a\omega_k)=\delta_k'+av_1+bv_0.
\]
All these values are gaps of $\Gamma_f$, by the identity $B_k=v_2-\varepsilon_1A_k$ proved in the definition of the blocks. Hence $B_k\subset\Lambda_f\setminus\Gamma_f$.\f

\teo\label{teo:decomposicao-rho-blocos-horizontais}
Let $f_1$ be a $1$-semiroot of $f$, with $g\geq2$, and suppose that the assumption of Case~2 holds. For each $i=0,\ldots,t$ and $\delta\in A_i$,
\[
 \rho_1(\delta)
 =\min\{\varepsilon_1\delta_i,\delta_i'\}
   +\varepsilon_1(\delta-\delta_i).
\]
In particular,
\[
 \rho_1(A_i)=
 \begin{cases}
  \varepsilon_1A_i,&\varepsilon_1\delta_i<\delta_i',\\
  B_i,&\varepsilon_1\delta_i>\delta_i'.
 \end{cases}
\]
If $j$ is the first index such that $\varepsilon_1\delta_j>\delta_j'$, then
\[
 \rho_1(\Lambda_1\setminus\Gamma_1)
 =\left(\dot\bigcup_{i=0}^{j-1}\varepsilon_1A_i\right)
  \mathbin{\dot\cup}
  \left(\dot\bigcup_{i=j}^{t}B_i\right).
\]
If no such index exists, then
\[
 \rho_1(\Lambda_1\setminus\Gamma_1)
 =\varepsilon_1(\Lambda_1\setminus\Gamma_1).
\]
\fteo

\dem By Lemma~\ref{lem:theta-blocos-horizontais}, for $\delta\in A_i$,
\[
 \beta_2+\varepsilon_1\nu_1(B_{\delta})
 =\delta_i'+\varepsilon_1(\delta-\delta_i).
\]
Substituting this equality into the definition of $\rho_1$, we obtain
\[
 \begin{aligned}
 \rho_1(\delta)
 &=\min\{\varepsilon_1\delta,
          \delta_i'+\varepsilon_1(\delta-\delta_i)\}\\
 &=\min\{\varepsilon_1\delta_i,\delta_i'\}
   +\varepsilon_1(\delta-\delta_i).
 \end{aligned}
\]
The two values being compared are distinct: the first is divisible by $\varepsilon_1$, and the second is congruent to $\beta_2$ modulo $\varepsilon_1$. This proves the description of $\rho_1(A_i)$.

The values $T_i$ defined in Lemma~\ref{lem:theta-blocos-horizontais} satisfy $T_0<\cdots<T_t$. Therefore, if $j$ is the first index with $\beta_2<T_j$, then $\beta_2<T_k$ for every $k\geq j$, whereas $\beta_2>T_i$ for $i<j$. Applying the description of $\rho_1(A_i)$ to the partition $\Lambda_1\setminus\Gamma_1=\dot\bigcup_i A_i$ yields the equalities in the statement. The unions are disjoint because the candidate blocks are disjoint.\f

\cor\label{lem:selecao-blocos-horizontais}
Under the assumptions of the preceding theorem, exactly one of the blocks $\varepsilon_1A_i$ and $B_i$ consists of special gaps of $f$, and the other contains no special gaps. More precisely,
\[
 \begin{aligned}
 \varepsilon_1\delta_i<\delta_i'
 &\ \Longrightarrow\quad
 \rho_1(A_i)=\varepsilon_1A_i,
 \quad B_i\cap(\Lambda_f\setminus\Gamma_f)=\emptyset,\\
 \varepsilon_1\delta_i>\delta_i'
 &\ \Longrightarrow\quad
 \rho_1(A_i)=B_i,
 \quad\varepsilon_1A_i\cap(\Lambda_f\setminus\Gamma_f)=\emptyset.
 \end{aligned}
\]
\fcor

\dem Set $C_i=\varepsilon_1A_i\mathbin{\dot\cup}B_i$. The sets $C_i$ are pairwise disjoint and are contained in $\mathbb N\setminus\Gamma_f$, below $v_2$. No element of $C_i$ belongs to the family
\[
 D=\{v_2-\varepsilon_1\alpha:
          \alpha\in\mathbb N^*\setminus\Lambda_1\}
\]
in Proposition~\ref{2.2.15}. For $\varepsilon_1A_i$, this follows from the residue classes modulo $\varepsilon_1$; for $B_i$, it follows from $B_i=v_2-\varepsilon_1A_i$ and $A_i\subset\Lambda_1$. Consequently,
\[
 C_i\cap(\Lambda_f\setminus\Gamma_f)
 =C_i\cap\rho_1(\Lambda_1\setminus\Gamma_1)=\rho_1(A_i),
\]
where the last equality uses the inclusion $\rho_1(A_k)\subset C_k$, already proved in the preceding theorem. The description of $\rho_1(A_i)$ in that theorem completes the proof.\f

\cor\label{cor:criterio-inclusao-escalada-horizontal}
Under the assumptions of Case~2,
\[
 \varepsilon_1(\Lambda_1\setminus\Gamma_1)
 \subset\Lambda_f\setminus\Gamma_f
\]
if and only if
\[
 \beta_2>
 \varepsilon_1\delta_t-(n_1-s_{t-1})v_1-v_0.
\]
\fcor

\dem By Corollary~\ref{lem:selecao-blocos-horizontais}, the inclusion is equivalent to $\varepsilon_1\delta_i<\delta_i'$ for every $i$, that is, $\beta_2>T_i$ for every $i$. Since $T_0<\cdots<T_t$ by Lemma~\ref{lem:theta-blocos-horizontais}, this condition is equivalent to $\beta_2>T_t$. Equality $\beta_2=T_i$ is impossible, since $T_i$ is divisible by $\varepsilon_1$, whereas $\beta_2$ is not.\f

\cor\label{cor:determinacao-torre-binaria-horizontal}
Suppose the semigroup $\Gamma_f=\langle v_0,\ldots,v_g\rangle$ is known, with $g\geq2$ and $n_k=2$ for $2\leq k\leq g$. If $f_1$ satisfies the assumption of Case~2, then $\Lambda_1\setminus\Gamma_1$ determines $\Lambda_f\setminus\Gamma_f$.
\fcor

\dem The characteristic exponents, in particular $\beta_2=v_2-(n_1-1)v_1$, are determined by $\Gamma_f$. Theorem~\ref{teo:decomposicao-rho-blocos-horizontais} determines $\rho_1(\Lambda_1\setminus\Gamma_1)$. Dividing this image by $\varepsilon_2$ and applying \eqref{eq:primeira-etapa-torre-binaria} gives $\Lambda_2\setminus\Gamma_2$ on the scale of $f_2$. The subsequent stages are determined by the recurrence~\eqref{eq:recorrencia-torre-binaria}, whose proof uses only $n_k=2$, without any assumption on the blocks of $f_1$. \f

\ob Under the assumptions of Corollary~\ref{cor:determinacao-torre-binaria-horizontal}, for $2\leq k\leq g$, the description by congruence classes also allows us to recover the gaps of the upper family from the lower family. Indeed, set $m_k=v_0/\varepsilon_{k-1}$ and take $\eta\in\mathbb N\setminus\Gamma_{k-1}$. Since $\eta$ is not a multiple of $m_k$, there exists $a\in\{1,\ldots,m_k-1\}$ with $a\equiv-\eta\pmod{m_k}$. We have $a\notin\Lambda_{k-1}$, since $\min\Lambda_{k-1}=m_k$. Moreover, $v_k/\varepsilon_k>v_0/\varepsilon_k=2m_k$, so that $v_k/\varepsilon_k-2a>0$, and
\[
 \left(\frac{v_k}{\varepsilon_k}+2\eta\right)
 -\left(\frac{v_k}{\varepsilon_k}-2a\right)
 =2(\eta+a)\in\frac{v_0}{\varepsilon_k}\mathbb N.
\]
Thus, a differential of value $v_k/\varepsilon_k+2\eta$ is obtained by multiplying a differential of value $v_k/\varepsilon_k-2a<v_k/\varepsilon_k$ by a power of $x$. This is a particular case of the covering of $\Lambda_k$ by its Apéry set; see Lemma~3.4 of \citet{de2022value}.\fob

\section{The general case: alternations between decompositions}
\label{sec:caso-geral}

In this section, the block decomposition is extended to sequences of minimal generators that alternate between Cases~1 and~2, without imposing a monotonicity condition on the entire sequence. We consider a branch $f$ with $g\geq2$ and its $1$-semiroot $f_1$. We retain the notation
\[
 n_1=\frac{v_0}{\varepsilon_1},\qquad
 m=\frac{v_1}{\varepsilon_1},\qquad
 \Gamma_1=\langle n_1,m\rangle.
\]
Thus, $\gcd(n_1,m)=1$, $\mu_1=(n_1-1)(m-1)$ is the conductor of $\Gamma_1$, and $v_2=(n_1-1)v_1+\beta_2$. If $\Lambda_1\setminus\Gamma_1=\emptyset$, then $\rho_1(\Lambda_1\setminus\Gamma_1)=\emptyset$, and the inclusion $\varepsilon_1(\Lambda_1\setminus\Gamma_1)\subset\Lambda_f\setminus\Gamma_f$ is immediate. We therefore assume $\Lambda_1\setminus\Gamma_1\neq\emptyset$. Following the convention of Section~\ref{sec:preliminares}, $B_\delta$ denotes the coefficient of $dy$ in a differential of value $\delta$ for which $\nu_1(B_\delta)$ is maximal. The valuations $\nu_1$ are computed on the scale of $f_1$.

\subsection{A two-level block decomposition}

Consider the minimal generators of $\Lambda_1$ that do not belong to $\Gamma_1$, ordered by increasing value. For each generator, its normal representation is $sm-rn_1$, with $0<s<n_1$ and $0<r<m$. The ordering of the coordinates $s$ agrees with that of the coordinates $r$. Indeed, if two distinct generators satisfied $s<s'$ and $r\geq r'$, then
\[
 (s'm-r'n_1)-(sm-rn_1)
 =(s'-s)m+(r-r')n_1\in\Gamma_1\setminus\{0\},
\]
contradicting the minimality of the generator $s'm-r'n_1$. Equality of two coordinates also contradicts minimality. Therefore, between consecutive generators, both coordinates increase (Case~1) or both decrease (Case~2).

We group the sequence into maximal consecutive segments in which $s$ and $r$ are increasing. These groups are indexed by $i=0,\ldots,q$, and the generators in each group by $j=1,\ldots,\ell_i$, so that
\[
 \delta_{i,1}<\cdots<\delta_{i,\ell_i},\qquad
 \delta_{i,\ell_i}<\delta_{i+1,1}\quad(0\leq i<q).
\]
Thus,
\[
 \delta_{i,j}=s_{i,j}m-r_{i,j}n_1,\qquad
 0<s_{i,j}<n_1,\quad 0<r_{i,j}<m.
\]
Each step within a group is of Case~1; the step to the next group is of Case~2. This indexing will be used for both generators and blocks.

Define recursively
\begin{equation}\label{eq:blocos-recursivos-gerais}
 \begin{aligned}
 A_{i,j}
 &=\bigl[(\Lambda_1\setminus\Gamma_1)
          \cap(\delta_{i,j}+\Gamma_1)\bigr]
       \setminus\bigcup_{\delta_{h,l}<\delta_{i,j}}A_{h,l},\\
 A_i&=\mathop{\dot\bigcup}_{j=1}^{\ell_i}A_{i,j}.
 \end{aligned}
\end{equation}
Since every special gap belongs to $\delta_{i,j}+\Gamma_1$ for some minimal generator $\delta_{i,j}$,
\[
 \Lambda_1\setminus\Gamma_1=\mathop{\dot\bigcup}_{i=0}^q
          \mathop{\dot\bigcup}_{j=1}^{\ell_i}A_{i,j}
    =\mathop{\dot\bigcup}_{i=0}^q A_i.
\]
The subblocks of each $A_i$ are added through vertical steps of Case~1; a horizontal step of Case~2 starts a new main block.

The bounds of each subblock depend on all preceding generators. To make them explicit, set
\begin{equation}\label{eq:limites-gerais}
 \begin{aligned}
 \underline r_{i,j}
 &=\max\bigl(\{0\}\cup
   \{r_{h,l}:\delta_{h,l}<\delta_{i,j},\ s_{h,l}<s_{i,j}\}\bigr),\\
 \overline s_{i,j}
 &=\min\bigl(\{n_1\}\cup
   \{s_{h,l}:\delta_{h,l}<\delta_{i,j},\ s_{h,l}>s_{i,j}\}\bigr).
 \end{aligned}
\end{equation}
We have $\underline r_{i,j}<r_{i,j}$ and $s_{i,j}<\overline s_{i,j}$. Definition~\eqref{eq:blocos-recursivos-gerais} is equivalent to the rectangular description
\begin{equation}\label{eq:retangulos-gerais}
 A_{i,j}=
 \left\{\delta_{i,j}+am+bn_1:
 \begin{array}{l}
 a,b\in\mathbb Z_{\geq0},\\
 a<\overline s_{i,j}-s_{i,j},\quad
 b<r_{i,j}-\underline r_{i,j}
 \end{array}\right\}.
\end{equation}
Indeed, for a positive gap $sm-rn_1$ in normal form,
\[
 sm-rn_1\in\delta_{h,l}+\Gamma_1
 \quad\Longleftrightarrow\quad
 s\geq s_{h,l},\quad r\leq r_{h,l},
\]
by the uniqueness of the normal form in $\langle n_1,m\rangle$. Among the gaps generated by $\delta_{i,j}$, the preceding generators with $s_{h,l}<s_{i,j}$ remove exactly the region $r\leq\underline r_{i,j}$; those with $s_{h,l}>s_{i,j}$ remove exactly the region $s\geq\overline s_{i,j}$. These observations prove \eqref{eq:retangulos-gerais}. In particular,
\begin{equation}\label{eq:extremos-blocos-gerais}
 \begin{aligned}
 \min A_{i,j}&=\delta_{i,j},\\
 \max A_{i,j}
   &=(\overline s_{i,j}-1)m-(\underline r_{i,j}+1)n_1=\delta_{i,j}+(\overline s_{i,j}-s_{i,j}-1)m+(r_{i,j}-\underline r_{i,j}-1)n_1,\\
 \sharp A_{i,j}
   &=(\overline s_{i,j}-s_{i,j})(r_{i,j}-\underline r_{i,j}).
 \end{aligned}
\end{equation}
The bounds \eqref{eq:limites-gerais} are not reset when the main block changes. This dependence on preceding generators is necessary to allow multiple alternations.

The corresponding blocks on the scale of $f$ are defined by
\begin{equation}\label{eq:blocos-b-gerais}
 \begin{aligned}
 \delta_{i,j}'&=\beta_2+(n_1-\overline s_{i,j})v_1
                          +(\underline r_{i,j}+1)v_0,\\
 B_{i,j}&=\delta_{i,j}'+\varepsilon_1(A_{i,j}-\delta_{i,j}),\\
 B_i&=\mathop{\dot\bigcup}_{j=1}^{\ell_i}B_{i,j}.
 \end{aligned}
\end{equation}
The definition of the blocks $B_{i,j}$ does not presuppose that their elements belong to $\Lambda_f$. The conditions for this inclusion will be established in Theorem~\ref{teo:decomposicao-rho-blocos-gerais}.

For $\delta=\delta_{i,j}+am+bn_1\in A_{i,j}$, we have
\[
\begin{aligned}
 v_2-\varepsilon_1\delta&=\delta_{i,j}'+(\overline s_{i,j}-s_{i,j}-1-a)v_1+(r_{i,j}-\underline r_{i,j}-1-b)v_0\\
 &=\delta_{i,j}'+\varepsilon_1\Bigl((\overline s_{i,j}-s_{i,j}-1-a)m+(r_{i,j}-\underline r_{i,j}-1-b)n_1\Bigr)\\
 &=\delta_{i,j}'+\varepsilon_1(\max A_{i,j}-\delta).
\end{aligned}
\]
Reflecting the two index intervals in \eqref{eq:retangulos-gerais} therefore gives
\begin{equation}\label{eq:reflexao-blocos-gerais}
 B_{i,j}=v_2-\varepsilon_1A_{i,j},\qquad
 \delta_{i,j}'=v_2-\varepsilon_1\max A_{i,j}.
\end{equation}
Since the $A_{i,j}$ form a partition, each of the families $\{\varepsilon_1A_{i,j}\}$ and $\{B_{i,j}\}$ consists of pairwise disjoint sets. The two families are also disjoint from each other, since their elements belong to the classes $0$ and $\beta_2$ modulo $\varepsilon_1$, respectively. For each pair $(i,j)$, we have $\sharp(\varepsilon_1A_{i,j})=\sharp B_{i,j}$. Moreover, all their elements are gaps of $\Gamma_f$ smaller than $v_2$. Indeed, below $v_2$ the elements of $\Gamma_f$ belong to $\varepsilon_1\Gamma_1$, whereas $\varepsilon_1A_{i,j}$ consists of gaps on this scale and $B_{i,j}\equiv\beta_2\not\equiv0\pmod{\varepsilon_1}$. Positivity and the upper bound follow from $A_{i,j}\subset[1,\mu_1-1]$ and $\varepsilon_1(\mu_1-1)<v_2$.

\subsection{Examples}

\ex\label{ex:alternancia-casos}
Consider
\[
 f_1=(t^8,t^{17}+t^{23}+t^{27})
\]
and the passage to
\[
 f=(t^{16},t^{34}+t^{46}+t^{54}+t^{55}).
\]
In this case,
\[
 \frac{v_0}{\varepsilon_1}=8,\qquad
 \frac{v_1}{\varepsilon_1}=17,\qquad
 n_1=8,\qquad
 \varepsilon_1=2,\qquad
 \beta_2=55.
\]
The minimal generators of $\Lambda_1$ that do not belong to $\Gamma_1$ are
\[
 \begin{aligned}
  \delta_{0,1}=31&=7\cdot17-11\cdot8,
      &&(s_{0,1},r_{0,1})=(7,11),\\
  \delta_{1,1}=52&=4\cdot17-2\cdot8,
      &&(s_{1,1},r_{1,1})=(4,2),\\
  \delta_{1,2}=70&=6\cdot17-4\cdot8,
      &&(s_{1,2},r_{1,2})=(6,4).
 \end{aligned}
\]
Therefore, the step from $(7,11)$ to $(4,2)$ is of Case~2, whereas the step from $(4,2)$ to $(6,4)$ is of Case~1. There are two main blocks:
\[
 A_0=A_{0,1},\qquad
 A_1=A_{1,1}\mathbin{\dot\cup}A_{1,2}.
\]
The explicit decomposition is
\[
 \begin{aligned}
 A_0=A_{0,1}
  &=\{31+8b:0\leq b\leq10\},\\
 A_{1,1}
  &=\{52+17a+8b:0\leq a\leq2,\ 0\leq b\leq1\},\\
 A_{1,2}
  &=\{70+8b:0\leq b\leq1\}.
 \end{aligned}
\]
Thus, $A_0$ is the horizontal row generated by $31$, and $A_1=(\Lambda_1\setminus\Gamma_1)\setminus A_0$ is subdivided vertically by the generators $52$ and $70$. In the following table, the columns are indexed by $r$ and the rows by $s$.
\[
\begin{array}{c|ccccccccccc}
s\backslash r
 &1&2&3&4&5&6&7&8&9&10&11\\ \hline
7
 &\cellcolor{brown!25}111&\cellcolor{brown!25}103
 &\cellcolor{brown!25}95&\cellcolor{brown!25}87
 &\cellcolor{brown!25}79&\cellcolor{brown!25}71
 &\cellcolor{brown!25}63&\cellcolor{brown!25}55
 &\cellcolor{brown!25}47&\cellcolor{brown!25}39
 &\cellcolor{brown!25}31\\
6
 &\cellcolor{blue!25}94&\cellcolor{blue!25}86
 &\cellcolor{red!25}78&\cellcolor{red!25}70\\
5
 &\cellcolor{blue!25}77&\cellcolor{blue!25}69\\
4
 &\cellcolor{blue!25}60&\cellcolor{blue!25}52
\end{array}
\]
\begin{center}
 \colorbox{brown!25}{$A_0$}\qquad
 \colorbox{blue!25}{$A_{1,1}$}\qquad
 \colorbox{red!25}{$A_{1,2}$}.
\end{center}

The accumulated bounds and the valuations of the maximal coefficients,
given by Lemma~\ref{lem:theta-blocos-gerais} as $\nu_1\left(B_{\delta_{i, j}}\right)=\frac{\delta_{i, j}^{\prime}-\beta_2}{\varepsilon_1}=\frac{\left(n_1-\bar{s}_{i, j}\right) v_1+\left(\underline{r}_{i, j}+1\right) v_0}{\varepsilon_1}=\left(n_1-\bar{s}_{i, j}\right) m+\left(\underline{r}_{i, j}+1\right) n_1$, are
\[
 \begin{array}{c|c|c|c}
 (i,j)&\delta_{i,j}&(\overline s_{i,j},\underline r_{i,j})
       &\nu_1(B_{\delta_{i,j}})\\ \hline
 (0,1)&31&(8,0)&8\\
 (1,1)&52&(7,0)&25\\
 (1,2)&70&(7,2)&41
 \end{array}
\]
Definition \eqref{eq:blocos-b-gerais} gives
\[
 \begin{aligned}
 \delta_{0,1}'&=55+16=71,
 &B_0=B_{0,1}
   &=\{71+16b:0\leq b\leq10\},\\
 \delta_{1,1}'&=55+34+16=105,
 &B_{1,1}
   &=\{105+34a+16b:0\leq a\leq2,\ 0\leq b\leq1\},\\
 \delta_{1,2}'&=55+34+3\cdot16=137,
 &B_{1,2}
   &=\{137+16b:0\leq b\leq1\}.
 \end{aligned}
\]
The comparisons between the minima are
\[
 \begin{array}{c|c|c}
 \text{subblock}&\varepsilon_1\delta_{i,j}&\delta_{i,j}'\\ \hline
 A_0&62&71\\
 A_{1,1}&104&105\\
 A_{1,2}&140&137
 \end{array}
\]
and show that the first alternative is selected in $A_0$ and
$A_{1,1}$, whereas the second is selected in $A_{1,2}$.
Theorem~\ref{teo:decomposicao-rho-blocos-gerais} gives, in this example,
\[
 \rho_1(\Lambda_1\setminus\Gamma_1)
 =2A_0\mathbin{\dot\cup}2A_{1,1}
  \mathbin{\dot\cup}B_{1,2}.
\]

\begin{figure}[H]
\centering
\caption{Differential values of the branch $f$. In brown, $2A_0$; in blue,
$2A_{1,1}$; in red, $B_{1,2}$, which replaces $2A_{1,2}$.}
\begingroup
\setlength{\tabcolsep}{1pt}
\renewcommand{\arraystretch}{1.05}
\resizebox{\textwidth}{!}{%
\begin{tabular}{r|*{33}{c}}
531 & \cellcolor{green!25}515  &  \cellcolor{green!25}499  &  \cellcolor{green!25}483  &  \cellcolor{green!25}467  &  \cellcolor{green!25}451  &  \cellcolor{green!25}435  &  \cellcolor{green!25}419  &  \cellcolor{green!25}403  &  \cellcolor{green!25}387  &  \cellcolor{green!25}371  &  \cellcolor{green!25}355  &  \cellcolor{green!25}339  &  \cellcolor{green!25}323  &  \cellcolor{green!25}307  &  \cellcolor{green!25}291  &  \cellcolor{green!25}275  &  259  &  243  &  227  &  211  &  195  &  179  &  163  &  147  &  131  &  115  &  99  &  83  &  67  &  51  &  35  &  19  &  3 \\
497 & \cellcolor{green!25}481  &  \cellcolor{green!25}465  &  \cellcolor{green!25}449  &  \cellcolor{green!25}433  &  \cellcolor{green!25}417  &  \cellcolor{green!25}401  &  \cellcolor{green!25}385  &  \cellcolor{green!25}369  &  \cellcolor{green!25}353  &  \cellcolor{green!25}337  &  \cellcolor{green!25}321  &  \cellcolor{green!25}305  &  \cellcolor{green!25}289  &  \cellcolor{green!25}273  &  \cellcolor{green!25}257  &  \cellcolor{green!25}241  &  225  &  209  &  193  &  177  &  161  &  145  &  129  &  113  &  97  &  81  &  65  &  49  &  33  &  17  &  1 \\
463 & \cellcolor{green!25}447  &  \cellcolor{green!25}431  &  \cellcolor{green!25}415  &  \cellcolor{green!25}399  &  \cellcolor{green!25}383  &  \cellcolor{green!25}367  &  \cellcolor{green!25}351  &  \cellcolor{green!25}335  &  \cellcolor{green!25}319  &  \cellcolor{green!25}303  &  \cellcolor{green!25}287  &  \cellcolor{green!25}271  &  \cellcolor{green!25}255  &  \cellcolor{green!25}239  &  \cellcolor{green!25}223  &  \cellcolor{green!25}207  &  191  &  175  &  159  &  143  &  127  &  111  &  95  &  79  &  63  &  47  &  31  &  15 \\
429 & \cellcolor{green!25}413  &  \cellcolor{green!25}397  &  \cellcolor{green!25}381  &  \cellcolor{green!25}365  &  \cellcolor{green!25}349  &  \cellcolor{green!25}333  &  \cellcolor{green!25}317  &  \cellcolor{green!25}301  &  \cellcolor{green!25}285  &  \cellcolor{green!25}269  &  \cellcolor{green!25}253  &  \cellcolor{green!25}237  &  \cellcolor{green!25}221  &  \cellcolor{green!25}205  &  189  &  173  &  157  &  141  &  125  &  109  &  93  &  77  &  61  &  45  &  29  &  13 \\
395 & \cellcolor{green!25}379  &  \cellcolor{green!25}363  &  \cellcolor{green!25}347  &  \cellcolor{green!25}331  &  \cellcolor{green!25}315  &  \cellcolor{green!25}299  &  \cellcolor{green!25}283  &  \cellcolor{green!25}267  &  \cellcolor{green!25}251  &  \cellcolor{green!25}235  &  \cellcolor{green!25}219  &  \cellcolor{green!25}203  &  \cellcolor{green!25}187  &  \cellcolor{green!25}171  &  155  &  139  &  123  &  107  &  91  &  75  &  59  &  43  &  27  &  11 \\
361 & \cellcolor{green!25}345  &  \cellcolor{green!25}329  &  \cellcolor{green!25}313  &  \cellcolor{green!25}297  &  \cellcolor{green!25}281  &  \cellcolor{green!25}265  &  \cellcolor{green!25}249  &  \cellcolor{green!25}233  &  \cellcolor{green!25}217  &  \cellcolor{green!25}201  &  \cellcolor{green!25}185  &  \cellcolor{green!25}169  &  \cellcolor{red!25}153  &  \cellcolor{red!50}137  &  121  &  105  &  89  &  73  &  57  &  41  &  25  &  9 \\
327 & \cellcolor{green!25}311  &  \cellcolor{green!25}295  &  \cellcolor{green!25}279  &  \cellcolor{green!25}263  &  \cellcolor{green!25}247  &  231  &  215  &  199  &  183  &  167  &  151  &  135  &  119  &  103  &  87  &  71  &  55  &  39  &  23  &  7 \\
293 & 277  &  261  &  245  &  229  &  213  &  197  &  181  &  165  &  149  &  133  &  117  &  101  &  85  &  69  &  53  &  37  &  21  &  5 \\
238 & \cellcolor{brown!25}222  &  \cellcolor{brown!25}206  &  \cellcolor{brown!25}190  &  \cellcolor{brown!25}174  &  \cellcolor{brown!25}158  &  \cellcolor{brown!25}142  &  \cellcolor{brown!25}126  &  \cellcolor{brown!25}110  &  \cellcolor{brown!25}94  &  \cellcolor{brown!25}78  &  \cellcolor{brown!50}62  &  46  &  30  &  14 \\
204 & \cellcolor{blue!25}188  &  \cellcolor{blue!25}172  &  156  &  140  &  124  &  108  &  92  &  76  &  60  &  44  &  28  &  12 \\
170 & \cellcolor{blue!25}154  &  \cellcolor{blue!25}138  &  122  &  106  &  90  &  74  &  58  &  42  &  26  &  10 \\
136 & \cellcolor{blue!25}120  &  \cellcolor{blue!50}104  &  88  &  72  &  56  &  40  &  24  &  8 \\
102 & 86  &  70  &  54  &  38  &  22  &  6 \\
68 & 52  &  36  &  20  &  4 \\
34 & 18  &  2 \\
0 &  \\
\end{tabular}%
}
\endgroup
\label{fig:alternancia-beta55}
\end{figure}

Comparing $f$ with the branch
\[
 f^+=(t^{16},t^{34}+t^{46}+t^{54}+t^{59})
\]
distinguishes the two regimes. The two branches have the same
semiroot $f_1$, but their third generators are $293$ and $297$,
respectively. The conditions
$\varepsilon_1\delta_{i,j}<\delta_{i,j}'$ are equivalent,
in the order of the minimal generators, to $\beta_2>46$, $\beta_2>54$, and $\beta_2>58$.
For $f^+$, with $\beta_2=59$, all three blocks remain in
the first alternative. For $f$, with $\beta_2=55$, only the
last subblock is replaced.
Green identifies the complementary families of
Theorem~\ref{t5.2}; the other colors distinguish the subblocks.

\begin{figure}[H]
\centering
\caption{Special gaps of $f^+$, with $\beta_2=59$.
In brown, $2A_0$; in blue, $2A_{1,1}$; in red, $2A_{1,2}$.}
\begingroup
\setlength{\tabcolsep}{1pt}
\renewcommand{\arraystretch}{1.05}
\resizebox{\textwidth}{!}{%
\begin{tabular}{r|*{33}{c}}
535 & \cellcolor{green!25}519  &  \cellcolor{green!25}503  &  \cellcolor{green!25}487  &  \cellcolor{green!25}471  &  \cellcolor{green!25}455  &  \cellcolor{green!25}439  &  \cellcolor{green!25}423  &  \cellcolor{green!25}407  &  \cellcolor{green!25}391  &  \cellcolor{green!25}375  &  \cellcolor{green!25}359  &  \cellcolor{green!25}343  &  \cellcolor{green!25}327  &  \cellcolor{green!25}311  &  \cellcolor{green!25}295  &  \cellcolor{green!25}279  &  263  &  247  &  231  &  215  &  199  &  183  &  167  &  151  &  135  &  119  &  103  &  87  &  71  &  55  &  39  &  23  &  7 \\
501 & \cellcolor{green!25}485  &  \cellcolor{green!25}469  &  \cellcolor{green!25}453  &  \cellcolor{green!25}437  &  \cellcolor{green!25}421  &  \cellcolor{green!25}405  &  \cellcolor{green!25}389  &  \cellcolor{green!25}373  &  \cellcolor{green!25}357  &  \cellcolor{green!25}341  &  \cellcolor{green!25}325  &  \cellcolor{green!25}309  &  \cellcolor{green!25}293  &  \cellcolor{green!25}277  &  \cellcolor{green!25}261  &  \cellcolor{green!25}245  &  229  &  213  &  197  &  181  &  165  &  149  &  133  &  117  &  101  &  85  &  69  &  53  &  37  &  21  &  5 \\
467 & \cellcolor{green!25}451  &  \cellcolor{green!25}435  &  \cellcolor{green!25}419  &  \cellcolor{green!25}403  &  \cellcolor{green!25}387  &  \cellcolor{green!25}371  &  \cellcolor{green!25}355  &  \cellcolor{green!25}339  &  \cellcolor{green!25}323  &  \cellcolor{green!25}307  &  \cellcolor{green!25}291  &  \cellcolor{green!25}275  &  \cellcolor{green!25}259  &  \cellcolor{green!25}243  &  \cellcolor{green!25}227  &  \cellcolor{green!25}211  &  195  &  179  &  163  &  147  &  131  &  115  &  99  &  83  &  67  &  51  &  35  &  19  &  3 \\
433 & \cellcolor{green!25}417  &  \cellcolor{green!25}401  &  \cellcolor{green!25}385  &  \cellcolor{green!25}369  &  \cellcolor{green!25}353  &  \cellcolor{green!25}337  &  \cellcolor{green!25}321  &  \cellcolor{green!25}305  &  \cellcolor{green!25}289  &  \cellcolor{green!25}273  &  \cellcolor{green!25}257  &  \cellcolor{green!25}241  &  \cellcolor{green!25}225  &  \cellcolor{green!25}209  &  193  &  177  &  161  &  145  &  129  &  113  &  97  &  81  &  65  &  49  &  33  &  17  &  1 \\
399 & \cellcolor{green!25}383  &  \cellcolor{green!25}367  &  \cellcolor{green!25}351  &  \cellcolor{green!25}335  &  \cellcolor{green!25}319  &  \cellcolor{green!25}303  &  \cellcolor{green!25}287  &  \cellcolor{green!25}271  &  \cellcolor{green!25}255  &  \cellcolor{green!25}239  &  \cellcolor{green!25}223  &  \cellcolor{green!25}207  &  \cellcolor{green!25}191  &  \cellcolor{green!25}175  &  159  &  143  &  127  &  111  &  95  &  79  &  63  &  47  &  31  &  15 \\
365 & \cellcolor{green!25}349  &  \cellcolor{green!25}333  &  \cellcolor{green!25}317  &  \cellcolor{green!25}301  &  \cellcolor{green!25}285  &  \cellcolor{green!25}269  &  \cellcolor{green!25}253  &  \cellcolor{green!25}237  &  \cellcolor{green!25}221  &  \cellcolor{green!25}205  &  \cellcolor{green!25}189  &  \cellcolor{green!25}173  &  157  &  141  &  125  &  109  &  93  &  77  &  61  &  45  &  29  &  13 \\
331 & \cellcolor{green!25}315  &  \cellcolor{green!25}299  &  \cellcolor{green!25}283  &  \cellcolor{green!25}267  &  \cellcolor{green!25}251  &  235  &  219  &  203  &  187  &  171  &  155  &  139  &  123  &  107  &  91  &  75  &  59  &  43  &  27  &  11 \\
297 & 281  &  265  &  249  &  233  &  217  &  201  &  185  &  169  &  153  &  137  &  121  &  105  &  89  &  73  &  57  &  41  &  25  &  9 \\
238 & \cellcolor{brown!25}222  &  \cellcolor{brown!25}206  &  \cellcolor{brown!25}190  &  \cellcolor{brown!25}174  &  \cellcolor{brown!25}158  &  \cellcolor{brown!25}142  &  \cellcolor{brown!25}126  &  \cellcolor{brown!25}110  &  \cellcolor{brown!25}94  &  \cellcolor{brown!25}78  &  \cellcolor{brown!50}62  &  46  &  30  &  14 \\
204 & \cellcolor{blue!25}188  &  \cellcolor{blue!25}172  &  \cellcolor{red!25}156  &  \cellcolor{red!50}140  &  124  &  108  &  92  &  76  &  60  &  44  &  28  &  12 \\
170 & \cellcolor{blue!25}154  &  \cellcolor{blue!25}138  &  122  &  106  &  90  &  74  &  58  &  42  &  26  &  10 \\
136 & \cellcolor{blue!25}120  &  \cellcolor{blue!50}104  &  88  &  72  &  56  &  40  &  24  &  8 \\
102 & 86  &  70  &  54  &  38  &  22  &  6 \\
68 & 52  &  36  &  20  &  4 \\
34 & 18  &  2 \\
0 &  \\
\end{tabular}%
}
\endgroup
\label{fig:alternancia-beta59}
\end{figure}

\fex

\ex\label{ex:multiplas-alternancias}
Consider the semiroot
\[
 f_1=(t^7,t^{16}+t^{18}+t^{19}),\qquad
 \Gamma_1=\langle7,16\rangle,\qquad \mu_1=90.
\]
Computing the module of differentials gives the minimal generators
\[
 \delta_{0,1}=25,\quad\delta_{0,2}=33,\quad
 \delta_{1,1}=50,\quad\delta_{1,2}=59,
\]
with the respective pairs
\[
 \begin{aligned}
 (s_{0,1},r_{0,1})&=(2,1),&
 (s_{0,2},r_{0,2})&=(6,9),\\
 (s_{1,1},r_{1,1})&=(4,2),&
 (s_{1,2},r_{1,2})&=(5,3).
 \end{aligned}
\]
The three transitions correspond, successively, to Cases $1$, $2$, and $1$.
Thus, there are two alternations and two main blocks, both
subdivided:
\[
 \begin{aligned}
 A_0&=A_{0,1}\mathbin{\dot\cup}A_{0,2},
 & A_{0,1}&=\{25+16a:0\leq a\leq4\},
 & A_{0,2}&=\{33+7b:0\leq b\leq7\},\\
 A_1&=A_{1,1}\mathbin{\dot\cup}A_{1,2},
 & A_{1,1}&=\{50,66\},
 & A_{1,2}&=\{59\}.
 \end{aligned}
\]
Figure~\ref{fig:multiplas-alternancias-f1} shows the decomposition
of the special gaps of $f_1$ into these four subblocks.

\begin{figure}[H]
\centering
\caption{Gaps of $\Gamma_1=\langle7,16\rangle$. The colors distinguish the four subblocks of
$\Lambda_1\setminus\Gamma_1$.}
\begingroup
\normalfont
\setlength{\tabcolsep}{3.5pt}
\renewcommand{\arraystretch}{1.12}
\begin{tabular}{r|*{13}{c}}
96 & \cellcolor{brown!25}89 & \cellcolor{blue!25}82 & \cellcolor{blue!25}75 & \cellcolor{blue!25}68 & \cellcolor{blue!25}61 & \cellcolor{blue!25}54 & \cellcolor{blue!25}47 & \cellcolor{blue!25}40 & \cellcolor{blue!50}{33} & 26 & 19 & 12 & 5 \\
80 & \cellcolor{brown!25}73 & \cellcolor{red!25}66 & \cellcolor{orange!50}{59} & 52 & 45 & 38 & 31 & 24 & 17 & 10 & 3 \\
64 & \cellcolor{brown!25}57 & \cellcolor{red!50}{50} & 43 & 36 & 29 & 22 & 15 & 8 & 1 \\
48 & \cellcolor{brown!25}41 & 34 & 27 & 20 & 13 & 6 \\
32 & \cellcolor{brown!50}{25} & 18 & 11 & 4 \\
16 & 9 & 2 \\
0 &  \\
\end{tabular}
\par\smallskip
{\small
 \colorbox{brown!25}{$A_{0,1}$}\quad
 \colorbox{blue!25}{$A_{0,2}$}\quad
 \colorbox{red!25}{$A_{1,1}$}\quad
 \colorbox{orange!25}{$A_{1,2}$}
}
\endgroup
\label{fig:multiplas-alternancias-f1}
\end{figure}

The data for each subblock are given below; the last column
uses the formula in Lemma~\ref{lem:theta-blocos-gerais}.
\[
 \begin{array}{c|c|c|c|c}
 (i,j)&\delta_{i,j}&\overline s_{i,j}&\underline r_{i,j}
       &\nu_1(B_{\delta_{i,j}})\\ \hline
 \cellcolor{brown!25}(0,1)&25&7&0&7\\
 \cellcolor{blue!25}(0,2)&33&7&1&14\\
 \cellcolor{red!25}(1,1)&50&6&1&30\\
 \cellcolor{orange!25}(1,2)&59&6&2&37
 \end{array}
\]
In particular, upon entering $A_1$, the lower bound remains
$\underline r_{1,1}=1$, and the upper bound becomes
$\overline s_{1,1}=6$, determined by $\delta_{0,2}$ rather than
$\delta_{0,1}$.
This dependence shows why accumulated bounds are necessary.

For
\[
 f=(t^{14},t^{32}+t^{36}+t^{38}+t^{39}),
\]
we have $\varepsilon_1=2$, $\beta_2=39$, and $v_2=231$.
In the order of the minimal generators, the values $\delta_{i,j}'$ are
$53$, $67$, $99$, and $113$, whereas the values
$2\delta_{i,j}$ are $50$, $66$, $100$, and $118$.
The description in Theorem~\ref{teo:decomposicao-rho-blocos-gerais} gives
\[
 \rho_1(\Lambda_1\setminus\Gamma_1)
 =2A_{0,1}\mathbin{\dot\cup}2A_{0,2}
  \mathbin{\dot\cup}B_{1,1}\mathbin{\dot\cup}B_{1,2},
 \qquad B_{1,1}=\{99,131\},\quad B_{1,2}=\{113\}.
\]
In this example, the second alternative is already selected in
the first subblock of $A_1$ and propagates to the next subblock.
In Figure~\ref{fig:multiplas-alternancias-f}, the colors of the
subblocks are retained in passing to $f$: brown and blue identify
$2A_{0,1}$ and $2A_{0,2}$, whereas red and orange identify
$B_{1,1}$ and $B_{1,2}$, respectively. The two complementary
families of Theorem~\ref{t5.2} are shown in green.

\begin{figure}[H]
\centering
\caption{Gaps of $\Gamma_f=\langle14,32,231\rangle$. The colored cells correspond to the
special gaps of $f$.}
\begingroup
\normalfont
\setlength{\tabcolsep}{1pt}
\renewcommand{\arraystretch}{1.12}
\resizebox{\textwidth}{!}{%
\begin{tabular}{r|*{30}{c}}
423 & \cellcolor{green!25}409 & \cellcolor{green!25}395 & \cellcolor{green!25}381 & \cellcolor{green!25}367 & \cellcolor{green!25}353 & \cellcolor{green!25}339 & \cellcolor{green!25}325 & \cellcolor{green!25}311 & \cellcolor{green!25}297 & \cellcolor{green!25}283 & \cellcolor{green!25}269 & \cellcolor{green!25}255 & \cellcolor{green!25}241 & \cellcolor{green!25}227 & \cellcolor{green!25}213 & 199 & 185 & 171 & 157 & 143 & 129 & 115 & 101 & 87 & 73 & 59 & 45 & 31 & 17 & 3 \\
391 & \cellcolor{green!25}377 & \cellcolor{green!25}363 & \cellcolor{green!25}349 & \cellcolor{green!25}335 & \cellcolor{green!25}321 & \cellcolor{green!25}307 & \cellcolor{green!25}293 & \cellcolor{green!25}279 & \cellcolor{green!25}265 & \cellcolor{green!25}251 & \cellcolor{green!25}237 & \cellcolor{green!25}223 & \cellcolor{green!25}209 & \cellcolor{green!25}195 & 181 & 167 & 153 & 139 & 125 & 111 & 97 & 83 & 69 & 55 & 41 & 27 & 13 \\
359 & \cellcolor{green!25}345 & \cellcolor{green!25}331 & \cellcolor{green!25}317 & \cellcolor{green!25}303 & \cellcolor{green!25}289 & \cellcolor{green!25}275 & \cellcolor{green!25}261 & \cellcolor{green!25}247 & \cellcolor{green!25}233 & \cellcolor{green!25}219 & \cellcolor{green!25}205 & \cellcolor{green!25}191 & \cellcolor{green!25}177 & \cellcolor{green!25}163 & 149 & 135 & 121 & 107 & 93 & 79 & 65 & 51 & 37 & 23 & 9 \\
327 & \cellcolor{green!25}313 & \cellcolor{green!25}299 & \cellcolor{green!25}285 & \cellcolor{green!25}271 & \cellcolor{green!25}257 & \cellcolor{green!25}243 & \cellcolor{green!25}229 & \cellcolor{green!25}215 & \cellcolor{green!25}201 & \cellcolor{green!25}187 & \cellcolor{green!25}173 & \cellcolor{green!25}159 & \cellcolor{green!25}145 & \cellcolor{red!25}131 & 117 & 103 & 89 & 75 & 61 & 47 & 33 & 19 & 5 \\
295 & \cellcolor{green!25}281 & \cellcolor{green!25}267 & \cellcolor{green!25}253 & \cellcolor{green!25}239 & \cellcolor{green!25}225 & \cellcolor{green!25}211 & \cellcolor{green!25}197 & \cellcolor{green!25}183 & \cellcolor{green!25}169 & \cellcolor{green!25}155 & \cellcolor{green!25}141 & \cellcolor{green!25}127 & \cellcolor{orange!50}113 & \cellcolor{red!50}99 & 85 & 71 & 57 & 43 & 29 & 15 & 1 \\
263 & \cellcolor{green!25}249 & \cellcolor{green!25}235 & \cellcolor{green!25}221 & \cellcolor{green!25}207 & \cellcolor{green!25}193 & \cellcolor{green!25}179 & 165 & 151 & 137 & 123 & 109 & 95 & 81 & 67 & 53 & 39 & 25 & 11 \\
231 & 217 & 203 & 189 & 175 & 161 & 147 & 133 & 119 & 105 & 91 & 77 & 63 & 49 & 35 & 21 & 7 \\
192 & \cellcolor{brown!25}178 & \cellcolor{blue!25}164 & \cellcolor{blue!25}150 & \cellcolor{blue!25}136 & \cellcolor{blue!25}122 & \cellcolor{blue!25}108 & \cellcolor{blue!25}94 & \cellcolor{blue!25}80 & \cellcolor{blue!50}66 & 52 & 38 & 24 & 10 \\
160 & \cellcolor{brown!25}146 & 132 & 118 & 104 & 90 & 76 & 62 & 48 & 34 & 20 & 6 \\
128 & \cellcolor{brown!25}114 & 100 & 86 & 72 & 58 & 44 & 30 & 16 & 2 \\
96 & \cellcolor{brown!25}82 & 68 & 54 & 40 & 26 & 12 \\
64 & \cellcolor{brown!50}50 & 36 & 22 & 8 \\
32 & 18 & 4 \\
0 &  \\
\end{tabular}%
}
\par\smallskip
{\small
 \colorbox{brown!25}{$2A_{0,1}$}\quad
 \colorbox{blue!25}{$2A_{0,2}$}\quad
 \colorbox{red!25}{$B_{1,1}$}\quad
 \colorbox{orange!25}{$B_{1,2}$}\quad
 \colorbox{green!25}{other special gaps}
}
\endgroup
\label{fig:multiplas-alternancias-f}
\end{figure}
\fex

\subsection{Accumulated bounds and coefficient valuations}

Delorme's theorem on increasing semimodules, in the formulation of
Theorem 7.13 (4) of \citet{cano2023analytic}, provides the
property of semimodules of differential values of branches with
one Puiseux pair that determines the evolution of the accumulated bounds.

\lem\label{lem:limites-acumulados-gerais}
If $\delta_{i,j}$ and $\delta_{h,l}$ are consecutive minimal
generators, in this order, then
\begin{equation}\label{eq:recorrencia-limites-gerais}
 (\overline s_{h,l},\underline r_{h,l})=
 \begin{cases}
  (\overline s_{i,j},r_{i,j}),
       &s_{h,l}>s_{i,j}\quad\text{\rm(Case~1)},\\
  (s_{i,j},\underline r_{i,j}),
       &s_{h,l}<s_{i,j}\quad\text{\rm(Case~2)}.
 \end{cases}
\end{equation}
Moreover,
\begin{equation}\label{eq:incrementos-minimos-gerais}
 \delta_{0,1}'=\beta_2+v_0,\qquad
 \delta_{h,l}'-\delta_{i,j}'
       \in\langle v_0,v_1\rangle\setminus\{0\}.
\end{equation}
In the order of the minimal generators, both the values $\delta_{i,j}'$
and the differences $\varepsilon_1\delta_{i,j}-\delta_{i,j}'$
form strictly increasing sequences.
\flem

\dem In this proof, the generators are indexed consecutively by
$\delta_0<\cdots<\delta_t$; the same index is used for their coordinates,
their accumulated bounds, and the values $\delta_k'$. For each $k$, define
the least common value
\[
 u_k=\min\left[
 \left(\Gamma_1\cup\bigcup_{h<k}(\delta_h+\Gamma_1)\right)
          \cap(\delta_k+\Gamma_1)\right].
\]
The rectangular description \eqref{eq:retangulos-gerais} gives
\begin{equation}\label{eq:eixos-gerais}
 \begin{aligned}
 u_k
 &=\delta_k+\min\{(\overline s_k-s_k)m,\,
                          (r_k-\underline r_k)n_1\}\\
 &=\min\{\overline s_km-r_kn_1,\,
                          s_km-\underline r_kn_1\}.
 \end{aligned}
\end{equation}
Indeed, reaching either boundary of the rectangle yields a
value in $\Gamma_1$ or in a previously generated semimodule;
when both increments remain strictly below the boundaries,
the value stays in the subblock of $\delta_k$.
Theorem 7.13 (4) of \citet{cano2023analytic}, applied to the basis
$n_1,m,\delta_0,\ldots,\delta_t$ of $\Lambda_1$, states that
\begin{equation}\label{eq:crescimento-delorme-geral}
 \delta_{k+1}>u_k\qquad(0\leq k<t).
\end{equation}

We prove the following property by induction on $k$, establishing
at each step the recurrence
\eqref{eq:recorrencia-limites-gerais} in passing from $k$ to $k+1$ as well:
\begin{equation}\label{eq:encaixe-acumulado-geral}
 s_h<\overline s_k,\qquad r_h>\underline r_k
 \qquad(k\leq h\leq t).
\end{equation}
For $k=0$, we have $\overline s_0=n_1$ and $\underline r_0=0$,
and the property follows from the normal representations of the generators.
Fix $k<t$ and assume, as the induction hypothesis, that
$s_h<\overline s_k$ and $r_h>\underline r_k$ for every $h\geq k$.
We prove the recurrence in passing from $k$ to $k+1$ and
\eqref{eq:encaixe-acumulado-geral} with $k$ replaced by $k+1$.
Set
\[
 X=s_km-\underline r_kn_1,\qquad
 Y=\overline s_km-r_kn_1.
\]
These values are distinct: equality would imply
$(\overline s_k-s_k)m=(r_k-\underline r_k)n_1$,
which is impossible, since $\gcd(n_1,m)=1$ and
$0<\overline s_k-s_k<n_1$.

If $X<Y$, then $u_k=X$. A subsequent generator with $s_h<s_k$
would also satisfy $r_h<r_k$. By the induction hypothesis,
$r_h>\underline r_k$, so that
\[
 \delta_h=s_hm-r_hn_1<X=u_k<\delta_{k+1}\leq\delta_h,
\]
a contradiction. Thus, all subsequent generators satisfy
$s_h>s_k$ and $r_h>r_k$. Since
$s_k<s_{k+1}<\overline s_k$ and no previous coordinate
$s_h$ lies strictly between $s_k$ and $\overline s_k$,
definition \eqref{eq:limites-gerais} gives
\[
 \overline s_{k+1}=\overline s_k,\qquad
 \underline r_{k+1}=r_k.
\]
This establishes Case~1 of \eqref{eq:recorrencia-limites-gerais}
and property \eqref{eq:encaixe-acumulado-geral} for $k+1$.

If $Y<X$, then $u_k=Y$. A subsequent generator with $s_h>s_k$
would satisfy $r_h>r_k$. By the induction hypothesis,
$s_h<\overline s_k$, so that
\[
 \delta_h=s_hm-r_hn_1<Y=u_k<\delta_{k+1}\leq\delta_h,
\]
again a contradiction. Therefore, all subsequent generators
satisfy $s_h<s_k$ and $r_h<r_k$.
For every previous index $h<k$ with $s_h<s_k$, we have
$r_h\leq\underline r_k<r_{k+1}$ and hence $s_h<s_{k+1}$.
The previous indices with $s_h>s_k$ remain to the right of
$s_{k+1}$. Thus, \eqref{eq:limites-gerais} implies
\[
 \overline s_{k+1}=s_k,\qquad
 \underline r_{k+1}=\underline r_k,
\]
which completes the inductive step.

By \eqref{eq:blocos-b-gerais}, $\delta_0'=\beta_2+v_0$ and
\[
 \delta_{k+1}'-\delta_k'=
 \begin{cases}
  (r_k-\underline r_k)v_0,&X<Y,\\
  (\overline s_k-s_k)v_1,&Y<X.
 \end{cases}
\]
These increments belong to
$\langle v_0,v_1\rangle\setminus\{0\}$ and equal
$\varepsilon_1(u_k-\delta_k)$. Finally,
\[
 \begin{aligned}
 &(\varepsilon_1\delta_{k+1}-\delta_{k+1}')
       -(\varepsilon_1\delta_k-\delta_k')\\
 &\qquad=\varepsilon_1(\delta_{k+1}-u_k)>0,
 \end{aligned}
\]
by \eqref{eq:crescimento-delorme-geral}. \f

The recurrence describes how the two decompositions combine:
a Case~1 transition preserves the upper bound and increases
the lower bound; a Case~2 transition preserves the lower bound
and decreases the upper bound. No additional restriction is
imposed on the number of alternations.

\lem\label{lem:theta-blocos-gerais}
For each subblock $A_{i,j}$ and every $\delta\in A_{i,j}$,
\begin{equation}\label{eq:theta-blocos-gerais}
 \beta_2+\varepsilon_1\nu_1(B_\delta)
 =\delta_{i,j}'+\varepsilon_1(\delta-\delta_{i,j}).
\end{equation}
In particular,
$\beta_2+\varepsilon_1\nu_1(B_{\delta_{i,j}})=\delta_{i,j}'$.
\flem

\dem Order the subblocks by their respective minima
$\delta_{i,j}$. Description~\eqref{eq:retangulos-gerais}
shows that each is rectangular. Moreover, the maxima
\[
 M_{i,j}=(\overline s_{i,j}-1)m
              -(\underline r_{i,j}+1)n_1
\]
are strictly decreasing in this order, since
\[
 \delta_{i,j}'=\beta_2+\varepsilon_1(K-M_{i,j}),
 \qquad K=(n_1-1)m,
\]
and the $\delta_{i,j}'$ are strictly increasing by
Lemma~\ref{lem:limites-acumulados-gerais}.
Lemma~\ref{lem:coeficientes-particao-retangular} therefore gives
\[
 \nu_1(B_\delta)=K-M_{i,j}+\delta-\delta_{i,j}
 \qquad(\delta\in A_{i,j}).
\]
Multiplying by $\varepsilon_1$ and adding $\beta_2$
yields \eqref{eq:theta-blocos-gerais}. \f

\subsection{Transition and general decomposition of the image}

\lem\label{lem:transicao-blocos-gerais}
If $\rho_1(\delta_{i,j})=\min(\varepsilon_1A_{i,j})$, then
$\varepsilon_1A_{i,j}\subset\Lambda_f\setminus\Gamma_f$.
If $\rho_1(\delta_{i,j})=\min B_{i,j}$, then
\[
 B_{h,l}\subset\Lambda_f\setminus\Gamma_f
 \qquad\text{for every }\delta_{h,l}\geq\delta_{i,j}.
\]
The second conclusion covers all subsequent subblocks,
including those belonging to other main blocks.
\flem

\dem If $\rho_1(\delta_{i,j})=\min(\varepsilon_1A_{i,j})$, choose a differential $\omega$
with $\nu_f(\omega)=\varepsilon_1\delta_{i,j}$.
For each $\delta_{i,j}+am+bn_1\in A_{i,j}$,
\[
 \nu_f(x^by^a\omega)
 =\varepsilon_1\delta_{i,j}+av_1+bv_0
 =\varepsilon_1(\delta_{i,j}+am+bn_1).
\]
These values are gaps of $\Gamma_f$, as shown
after \eqref{eq:reflexao-blocos-gerais}.

If $\rho_1(\delta_{i,j})=\min B_{i,j}$, choose $\omega$ with
$\nu_f(\omega)=\delta_{i,j}'$.
By Lemma~\ref{lem:limites-acumulados-gerais},
$\delta_{h,l}'-\delta_{i,j}'\in\langle v_0,v_1\rangle$
whenever $\delta_{h,l}\geq\delta_{i,j}$.
Write this difference as $av_1+bv_0$, with $a,b\geq0$.
Then
\[
 \nu_f(x^by^a\omega)=\delta_{h,l}'.
\]
Multiplying this differential by the monomials corresponding to
the rectangle $A_{h,l}-\delta_{h,l}$ yields all values of
$B_{h,l}$. The identity $B_{h,l}=v_2-\varepsilon_1A_{h,l}$
ensures that these are gaps of $\Gamma_f$ and hence special gaps.
The construction follows from the additivity of the valuation
and remains valid independently of maximizing the value of
the coefficient of $dy$ in the multiplied differentials. \f

\teo\label{teo:decomposicao-rho-blocos-gerais}
Let $f_1$ be a $1$-semiroot of a branch $f$ with $g\geq2$.
For each subblock $A_{i,j}$ and every $\delta\in A_{i,j}$,
\begin{equation}\label{eq:rho-blocos-gerais}
 \rho_1(\delta)=
 \min\{\varepsilon_1\delta_{i,j},\delta_{i,j}'\}
       +\varepsilon_1(\delta-\delta_{i,j}).
\end{equation}
Consequently,
\[
 \rho_1(A_{i,j})=
 \begin{cases}
  \varepsilon_1A_{i,j},&\varepsilon_1\delta_{i,j}<\delta_{i,j}',\\
  B_{i,j},&\varepsilon_1\delta_{i,j}>\delta_{i,j}'.
 \end{cases}
\]
If $(I,J)$ is the first pair, in the order of the minimal generators,
such that $\rho_1(\delta_{I,J})=\min B_{I,J}$, then
\begin{equation}\label{eq:corte-hierarquico-geral}
 \begin{aligned}
 \rho_1(\Lambda_1\setminus\Gamma_1)={}&
 \left(\mathop{\dot\bigcup}_{i=0}^{I-1}\varepsilon_1A_i\right)
 \mathbin{\dot\cup}
 \left(\mathop{\dot\bigcup}_{j=1}^{J-1}
                              \varepsilon_1A_{I,j}\right)\\
 &\mathbin{\dot\cup}
 \left(\mathop{\dot\bigcup}_{j=J}^{\ell_I}B_{I,j}\right)
 \mathbin{\dot\cup}
 \left(\mathop{\dot\bigcup}_{i=I+1}^{q}B_i\right).
 \end{aligned}
\end{equation}
If there is no such pair, then
$\rho_1(\Lambda_1\setminus\Gamma_1)
=\varepsilon_1(\Lambda_1\setminus\Gamma_1)$.
Finally, for each $(i,j)$,
\begin{equation}\label{eq:exclusao-blocos-gerais}
 \bigl(\varepsilon_1A_{i,j}\mathbin{\dot\cup}B_{i,j}\bigr)
       \cap(\Lambda_f\setminus\Gamma_f)
 =\rho_1(A_{i,j}).
\end{equation}
Thus, the unselected block contains no special gaps of $f$.
\fteo

\dem Substituting the identity in Lemma~\ref{lem:theta-blocos-gerais}
into the definition of $\rho_1$ gives
\eqref{eq:rho-blocos-gerais}. The two terms being compared are
distinct, since they belong to the classes $0$ and $\beta_2$ modulo
$\varepsilon_1$. This yields the formula for $\rho_1(A_{i,j})$.

By Lemma~\ref{lem:limites-acumulados-gerais}, the differences
$\varepsilon_1\delta_{i,j}-\delta_{i,j}'$ are strictly
increasing in the order of the generators. Thus, the first
positive difference determines the second alternative for
all subsequent subblocks. Applying the formula for $\rho_1(A_{i,j})$
to the partition of $\Lambda_1\setminus\Gamma_1$ gives
\eqref{eq:corte-hierarquico-geral}; in the absence of a positive
difference, it gives
$\rho_1(\Lambda_1\setminus\Gamma_1)
=\varepsilon_1(\Lambda_1\setminus\Gamma_1)$.
The unions are disjoint by the geometric properties of the
blocks. Lemma~\ref{lem:transicao-blocos-gerais} explicitly
realizes the propagation of the $B_{i,j}$ by multiplication
of differentials.

By the identity $B_{i,j}=v_2-\varepsilon_1A_{i,j}$ and the inclusion
$A_{i,j}\subset\Lambda_1$, the block $B_{i,j}$ is disjoint from
\[
 \left\{v_2-\varepsilon_1\delta;
 \delta\in\mathbb N^*\setminus\Lambda_1\right\},
\]
which occurs in Proposition~\ref{2.2.15}; the same holds for
$\varepsilon_1A_{i,j}$, by the residue classes modulo $\varepsilon_1$.
Since both lie below $v_2$, the disjointness of the blocks and
the formula already obtained for $\rho_1(A_{i,j})$ give
\[
 \bigl(\varepsilon_1A_{i,j}\mathbin{\dot\cup}B_{i,j}\bigr)
       \cap(\Lambda_f\setminus\Gamma_f)=\rho_1(A_{i,j}),
\]
which is \eqref{eq:exclusao-blocos-gerais}. \f

Formula \eqref{eq:corte-hierarquico-geral} formalizes the order
of the checks: within each main block, the minima of the subblocks
are examined in the order $j=1,\ldots,\ell_i$. As long as
$\rho_1(\delta_{i,j})=\min(\varepsilon_1A_{i,j})$, one proceeds
to the next minimum. The first equality
$\rho_1(\delta_{i,j})=\min B_{i,j}$ determines the second
alternative for all subsequent subblocks, even after further
alternations between Cases~1 and~2.

\cor\label{cor:criterio-inclusao-escalada-geral}
Under the hypotheses of the preceding theorem, with
$\Lambda_1\setminus\Gamma_1\neq\emptyset$,
\[
 \varepsilon_1(\Lambda_1\setminus\Gamma_1)
 \subset\Lambda_f\setminus\Gamma_f
 \quad\Longleftrightarrow\quad
 \varepsilon_1\delta_{q,\ell_q}<\delta_{q,\ell_q}'.
\]
Equivalently,
\begin{equation}\label{eq:criterio-inclusao-geral}
 \beta_2>
 \varepsilon_1\delta_{q,\ell_q}
 -(n_1-\overline s_{q,\ell_q})v_1
 -(\underline r_{q,\ell_q}+1)v_0.
\end{equation}
\fcor

\dem By identity \eqref{eq:exclusao-blocos-gerais}, the inclusion
is equivalent to selecting $\varepsilon_1A_{i,j}$ for every $(i,j)$,
that is, to $\varepsilon_1\delta_{i,j}<\delta_{i,j}'$ for
all subblocks. By the strict monotonicity in
Lemma~\ref{lem:limites-acumulados-gerais}, it suffices that
this inequality hold for the last generator
$\delta_{q,\ell_q}$. Substituting the expression for
$\delta_{q,\ell_q}'$ into \eqref{eq:blocos-b-gerais} gives
\eqref{eq:criterio-inclusao-geral}. Equality between the minima
is impossible, since $\varepsilon_1\delta_{i,j}$ is divisible by
$\varepsilon_1$, whereas $\delta_{i,j}'$ is not. \f

In the monotone cases, the formulas of the preceding sections are recovered.
Using, only for this comparison, the indexing
$\delta_0<\cdots<\delta_t$ of
Sections~\ref{sec:caso-vertical} and~\ref{sec:caso-horizontal}, if every transition
is of Case~1, we have $\overline s_i=n_1$ and
$\underline r_i=r_{i-1}$, with $r_{-1}=0$. Thus,
\[
 \delta_i'=\beta_2+(r_{i-1}+1)v_0,
\]
and the inclusion criterion is
$\beta_2>\varepsilon_1\delta_t-(r_{t-1}+1)v_0$.
If every transition is of Case~2, then
$\overline s_i=s_{i-1}$ and $\underline r_i=0$, with $s_{-1}=n_1$,
so that
\[
 \delta_i'=\beta_2+(n_1-s_{i-1})v_1+v_0,
\]
and the criterion is
$\beta_2>\varepsilon_1\delta_t-(n_1-s_{t-1})v_1-v_0$.

\subsection{Determination along a binary tower}

\cor\label{cor:determinacao-torre-binaria-geral}
Suppose that the semigroup
$\Gamma_f=\langle v_0,\ldots,v_g\rangle$ of a branch $f$ is known, with
$g\geq2$ and $n_k=2$ for $2\leq k\leq g$.
Then $\Lambda_1\setminus\Gamma_1$ determines
$\Lambda_f\setminus\Gamma_f$, independently of the alternations
between Cases~1 and~2 in the sequence of minimal generators.
\fcor

\dem The semigroup $\Gamma_f$ determines the semigroups
$\Gamma_k$ and the characteristic exponents, in particular
$\beta_2=v_2-(n_1-1)v_1$.
If $\Lambda_1\setminus\Gamma_1\neq\emptyset$,
Theorem~\ref{teo:decomposicao-rho-blocos-gerais} determines
$\rho_1(\Lambda_1\setminus\Gamma_1)$; if the set is empty,
its image is also empty.
Formula~\eqref{eq:primeira-etapa-torre-binaria}, which includes
division by $\varepsilon_2$, determines
$\Lambda_2\setminus\Gamma_2$ on the scale of $f_2$.
Lemma~\ref{lem:propagacao-total-binaria} and
recurrence~\eqref{eq:recorrencia-torre-binaria} successively
determine the special gaps of the remaining semiroots,
up to $\Lambda_g\setminus\Gamma_g=\Lambda_f\setminus\Gamma_f$.
None of these subsequent steps uses the initial block decomposition. \f

\section*{Use of artificial intelligence}
OpenAI ChatGPT/Codex was used to assist with the English translation,
bibliographic editing, and \LaTeX{} formatting of this manuscript.

\bibliographystyle{plainnat}
\bibliography{references}
\end{document}